\documentclass[matprg]{svjour3}
\usepackage{subcaption}
\usepackage{amsmath,amsxtra,amsfonts,amscd,amssymb,bm}
\usepackage{algorithm}
\usepackage{algpseudocode}
\usepackage[top=2.5cm,bottom=2.5cm,right=2.5cm,left=2.5cm]{geometry}
\usepackage[mathscr]{eucal}
\usepackage{color}
\usepackage{cite}
\numberwithin{equation}{section}
\usepackage{hyperref}

\def\wh{\widehat}

\def\sl{\textstyle{\sum}_{\ell=1}^k}

\def\O{\mathcal{O}}

\def\cG{{\cal G}}
\def\cF{{\cal F}}

\newcommand{\tsum}{\textstyle\sum}
\newcommand{\bbe}{\mathbb{E}}

\def\eqref#1{(\ref{#1})}
\usepackage{epsfig}
\usepackage{calc}
\usepackage{amstext}
\usepackage{amsmath}
\usepackage{multicol}
\usepackage{amsfonts}
\usepackage{amssymb}
\usepackage{mathrsfs}
\usepackage{bm}
\usepackage{xcolor}
\usepackage{comment}
\usepackage[scientific-notation=true]{siunitx}

\newcommand{\calF}{\mathcal{F}}
\newcommand{\bbec}[1]{\bbe_{\lceil #1 \rceil}}

\newcommand{\bbr}{\Bbb{R}}
\newcommand{\beq}{\begin{equation}}
\newcommand{\eeq}{\end{equation}}
\newcommand{\beqa}{\begin{eqnarray}}
\newcommand{\eeqa}{\end{eqnarray}}
\newcommand{\beqas}{\begin{eqnarray*}}
\newcommand{\eeqas}{\end{eqnarray*}}

\newcommand{\newOmega}{\overline \Omega}

\usepackage{algorithm}
\usepackage{algpseudocode}
\def\argmin{\mathop{\mathrm{argmin}}}
\def\argmax{\mathop{\mathrm{argmax}}}

\def\bbe{{\mathbb{E}}}

\newcommand{\nn}{\nonumber}

\def\qed{\ \hfill$\square$\par\smallskip}
\newcommand{\aic}[2]{{\color{black}#2}}
\newcommand{\revision}[2]{{\color{black}#2}}
\newcommand{\tli}[2]{{\color{black}#2}}
\newcommand{\rf}[1]{(\ref{#1})}
\renewcommand{\top}{{T}}
\def\inter{\hbox{\rm  int\,}}
\def\cL{{\cal L}}
\def\cH{{\cal H}}

\def\cK{{\cal K}}

\def\euc{E}
\def\feaReg{X}

\def\xmid{y}
\def\xavg{x}

\def\stx{{y}}
\def\sstx{\overline{y}}

\title{Accelerated stochastic approximation with state-dependent noise
}
\author{
  Sasila Ilandarideva$^{1}$
  \and
  Anatoli Juditsky$^{1}$
\and
  Guanghui Lan$^{2}$
\and
  Tianjiao Li$^{2}$
   \thanks{SI and AJ were partially supported by Multidisciplinary Institute in Artificial intelligence MIAI {@} Grenoble Alpes ANR-19-P3IA-0003. GL and TL were partially supported by Division of Mathematical Sciences grant DMS-1953199 and Air Force Office of Scientific Research grant FA9550-22-1-0447.}
   }

  \institute{
  Corresponding Author: Tianjiao Li. Coauthors are listed according to the alphabetic order.
  \vspace{0.08in}\\
  $^{1}$ LJK, University Grenoble Alpes, 38401 Domaine Universitaire de Saint-Martin-d'Hères, France\\
  Email: {sasila.ilandarideva@univ-grenoble-alpes.fr; anatoli.juditsky@univ-grenoble-alpes.fr}
  \vspace{0.06in}\\
  $^{2}$ H. Milton Stewart School of Industrial and Systems
    Engineering, Georgia Institute of Technology, Atlanta, GA, USA\\
Email: {george.lan@isye.gatech.edu; tli432@gatech.edu}}
\date{the date of receipt and acceptance should be inserted later}

\begin{document}

\maketitle

\begin{abstract}
%\aic{
%This paper considers a class of stochastic smooth convex optimization problems under a fairly general assumption on the noise of stochastic gradients. Different from the classical setting with uniformly bounded noise, herein we assume that the variance of stochastic gradients depends on the sub-optimality of the search points. These problems arise naturally from a few different applications, including the well-known generalized linear regression in statistics. However, none of the existing stochastic optimization methods are optimal in terms of the dependence on accuracy, problem parameters, and mini-batch size
%for solving this class of problems. To address this issue, we study and compare two accelerated stochastic approximation methods with a certain duality relationship. We show that the first method, known as the stochastic accelerated gradient descent (SAGD), can achieve a nearly optimal convergence up to the dependence on the mini-batch size.
%While the second method, named stochastic gradient extrapolation (SGE), achieves the optimal convergence rate, thus attaining the optimal iteration and sample complexities simultaneously for the first time in the literature. We then generalize the latter method to solve problems satisfying certain quadratic growth conditions, and specialize it for the classic sparse recovery problem in statistical learning. Finally, we report some numerical experiments to demonstrate the performance of our proposed algorithms in high-dimensional settings. }
{
We consider a class of stochastic smooth convex optimization problems under rather general assumptions on the noise in the stochastic gradient observation. As opposed to the classical problem setting in which the variance of noise is assumed to be uniformly bounded, herein we assume that the variance of stochastic gradients is related to the ``sub-optimality'' of the approximate solutions delivered by the algorithm. Such problems naturally arise in a variety of applications, in particular, in the well-known generalized linear regression problem in statistics. However, to the best of our knowledge, none of the existing stochastic approximation algorithms for solving this class of problems attain optimality in terms of the dependence on accuracy, problem parameters, and mini-batch size.

We discuss two non-Euclidean accelerated stochastic approximation routines---stochastic accelerated gradient descent (SAGD) and stochastic gradient extrapolation (SGE)---which carry a particular duality relationship. We show that both SAGD and SGE, under appropriate conditions, achieve the optimal convergence rate, attaining the optimal iteration and sample complexities simultaneously. However, corresponding assumptions for the SGE algorithm are more general; they allow, for instance, for efficient application of the SGE to statistical estimation problems under heavy tail noises and discontinuous score functions. We also discuss the application of the SGE to problems satisfying quadratic growth conditions, and show how it can be used to recover sparse solutions. Finally, we report on some simulation experiments to illustrate numerical performance of our proposed algorithms in high-dimensional settings.}
 \end{abstract}


 \section{Introduction} \label{sec_intro}
This paper focuses on the stochastic
optimization problem given by
\beq \label{cp-ASGD}
f^* := \min\limits_{x \in \feaReg}  f(x)
\eeq
where $\feaReg$ is a closed convex subset of a Euclidean space $\euc$
and $f: \feaReg \to \bbr$ is a smooth convex function with Lipschitz continuous gradient, i.e.,
for some $L \ge 0$,
\beq \label{smoothness_ASGD}
0 \le f(y) - f(x) - \langle \nabla f(x), y- x\rangle
\le \tfrac{L}{2} \|y- x\|^2, \ \ \ \forall x, y \in \feaReg.
\eeq
%Here $g(x)$ denotes the gradient of $f$ at $x$.
We assume throughout the paper that the problem in \rf{cp-ASGD} is solvable, i.e., 
the set  $\feaReg^*$ of optimal solutions
is nonempty.

We consider the stochastic setting where only stochastic first-order information about $f$ is available for solving problem \eqref{cp-ASGD}.  Specifically, at the current search point $x_t \in \feaReg$,  a stochastic oracle (SO) generates the stochastic operator $\cG(x_t, \xi_t)$, where $\xi_t \in \Xi$ denotes a random variable, whose probability distribution is supported on $\Xi$. We suppose that $\xi_t$ is independent of $x_0, ..., x_t$, and $\{\xi_t\}_{t\geq 0}$ are mutually independent; we
also assume that  $\cG(x_t, \xi_t)$ is an unbiased estimator of $g(x_t)=\nabla f(x_t)$ satisfying
\begin{align}\label{assump:bias_SO}
    \bbe_{\xi_t} [\cG(x_t, \xi_t)] = g(x_t)
\end{align}
(expectation w.r.t. to the distribution $\xi_t$).

Stochastic approximation (SA) and stochastic mirror descent (SMD) methods are routinely used to solve stochastic optimization problems; see, e.g., \cite{nemirovskii1979complexity,polyak1990new,polyak1992acceleration, nemirovski2009robust}.
More specifically, it was shown in \cite{nemirovski2009robust} that SMD can achieve the optimal sample complexity for general nonsmooth optimization and saddle point problems. For smooth stochastic optimization problems, Lan \cite{lan2012optimal,LanBook2020} introduced an accelerated stochastic approximation (AC-SA), also known as stochastic accelerated gradient descent (SAGD), which was obtained by replacing exact gradients with their unbiased estimators in the celebrated accelerated gradient methods \cite{nesterov1983method} (see also
\cite{YudinNemirovski77,NemirovskiYudin83-smooth,nemyud:83,NemirovskiNester85}
for early developments).
It was shown in \cite{lan2012optimal} that AC-SA achieves the optimal sample complexity for smooth, nonsmooth and stochastic convex optimization (see \cite{ghadimi2012optimal, ghadimi2013optimal} for generalization to the strongly convex settings). It should be noted that the original analysis of AC-SA in \cite{lan2012optimal} was carried out under {\em uniformly bounded variance condition}
% For instance, \cite{nemirovski2009robust} demonstrated the effectiveness of SA methods for non-smooth stochastic convex optimization; a stochastic accelerated gradient descent (SAGD) algorithm  for the smooth stochastic convex optimization, a stochastic variant of the celebrated Nesterov's accelerated gradient decent \cite{nesterov1983method, LanBook2020}, was proposed in \cite{lan2012optimal} which also established the optimal convergence rate under the {\em bounded variance condition}
\beq\label{assump_uniform_noise}
\bbe_{\xi_t}[\|\cG({x}, \xi_t) - g({x})\|_*^2] \le \sigma^2
\eeq
where $\|\cdot\|_*$ is the norm conjugate to $\|\cdot\|$.
However, it has been observed recently (see, e.g., \cite{srebro2010smoothness, juditsky2023sparse}) that this uniformly bounded variance condition is not necessarily satisfied in some important applications in which the variance of the stochastic gradient depends on the search point $x_t$.
% Recently, there has been an increasing interest in using SA type methods for solving smooth stochastic convex optimization with {\em state-dependent noise;} see, e.g., \cite{srebro2010smoothness, juditsky2023sparse}.
As a motivation, consider
the fundamental to Statistical Learning problem of parameter estimation in the {\em generalized linear regression} (GLR) model in which one aims to estimate an unknown parameter vector $x^*\in \feaReg\subset \bbr^n$ given observations $(\phi_t,\eta_t)$,
\begin{align}\label{prob_GLM}
    \eta_t = u(\phi^\top_t x^*) +  \zeta_t,\quad t=1,2,...,
\end{align}
where, in generalized linear models terminology, $\eta_t \in
\bbr$ are responses,  $\phi_t \in \bbr^n$ are random regressors, $\zeta_t \in \bbr$ are zero-mean random noises which are assumed to be mutually independent and independent of $\phi_t$, and $u: \bbr\rightarrow \bbr$ is the (generally nonlinear) ``activation function''.
Then it follows directly that
\begin{equation} \label{eq:opt_GLR}
\bbe[\phi_t(u(\phi^\top_t x^*) - \eta_t)] =
\bbe[\phi_t \, \zeta_t] =0.
\end{equation}
Thus, the problem of recovery of $x^*$ from observations $\eta_t$ and $\phi_t$ may be formulated as
 a stochastic optimization problem.
 Specifically,
 when denoting $v: \bbr \to \bbr$ the primitive of $u$, i.e., $v'(t) = u(t)$ and assuming that $x^*\in\inter \feaReg$, \eqref{eq:opt_GLR} may be seen as
 as the optimality condition
 %for
 %this problem. More specifically,
%and denote .
%This setting implies that $\bbe[\phi(u(\phi^\top x^*) - \eta)] = \mathbf{0}$, thus the problem can be solved by considering
for the problem
\beq\label{prob_GLM_opt}
    \min_{x \in \feaReg} \left\{f(x) := \bbe[v(\phi^\top x) - \phi^\top x \eta]\right\}.
\eeq
%where $v'(t) = u(t)$.
Clearly, the gradient of $f$ is given by $g(x) = \bbe[\phi(u(\phi^\top x) - \eta)]$, and one of its unbiased stochastic estimator is  $\cG(x, \overbrace{(\phi, \zeta)}^{=:\xi}) = \phi(u(\phi^\top x) - \eta)$.
Under mild assumptions, one can show (see  Section~\ref{sec:glr}) that the noise
\[\cG({x}, (\phi_t, \zeta_t))-g({x})=\phi_t[u(\phi_t^\top {x})-u(\phi_t^\top x^*)]-\phi_t\zeta_t + g(x^*) - g({x})\]
(\revision{}{here $g(x^*)=0$}) of the gradient observation $\cG({x},(\phi_t, \zeta_t))$ at ${x}$ satisfies the condition
\beq\label{cond:glr}
\bbe_{\xi_t}\big[\|\cG({x}, \xi_t) - g({x})\|_*^2\big] \le \sigma^2({x})=\cL [f({x}) - f^*]+\sigma_*^2 \tag{SN}
\eeq
for some $\cL, \sigma_*>0$. In what follows, with some terminology abuse, we refer to \rf{cond:glr} and similar conditions as {\em state-dependent noise} assumptions. More generally, if compared to the ``uniform noise'' condition \rf{assump_uniform_noise}, \rf{cond:glr} may be seen as a refined  assumption on the structure of the stochastic oracle. Furthermore, one can easily verify that in various settings of the GLR problem, condition \rf{cond:glr} holds, while condition \eqref{assump_uniform_noise} of uniformly bounded noise variance  is violated (e.g., in the case of unbounded $\feaReg$ and linear $u$). {Similar conditions have been recently introduced in the context of reduced variance stochastic algorithms for finite sum minimization, see, e.g.,  \cite{bietti2017stochastic,gower2019sgd,gower2021sgd, khaled2022better} and references therein.} Various stochastic optimization problems in which similar state-dependent noise assumptions apply can also be found in recent literature on
machine learning \cite{woodworth2021even,juditsky2023sparse,ilandarideva2023stochastic} and
reinforcement learning, see, e.g., \cite{ kotsalis2022simple1, kotsalis2022simple, li2021accelerated, li2022stochastic}.

% In spite of these efforts, some important questions about efficient and robust first-order methods for stochastic optimization with state-dependent noise remained unanswered.
Motivated by these aforementioned statistical applications, Juditsky et al. \cite{juditsky2023sparse} proposed an SMD algorithm that exploits state-dependent noise assumption to attain optimal convergence rates in the situation of ``dominating stochastic error'', i.e., when the amplitude of the error of the gradient observation is comparable to the amplitude of the gradient $\nabla f$ of the problem objective. In the similar setting, the authors of \cite{woodworth2021even} have recently established sharp lower complexity bounds for stochastic optimization under state-dependent noise in the Euclidean setting for both general and strongly convex situations.
However, it is well-known that SMD is suboptimal in the so-call mini-batch setting, where the noise of the gradient estimator is reduced by using a batch of samples. This setting has been widely used for applications of stochastic optimization algorithms especially under a distributed computing environment. To achieve the accelerated convergence, the application of the classical SAGD algorithm of \cite{lan2012optimal} to the state-dependent noise setting was the subject of  \cite{cotter2011better, liu2018mass,woodworth2021even}.
The authors of \cite{woodworth2021even} proved (expected) optimal accuracy bound $\mathcal{O}\big(\cH/k^2\big)$ after $k$ iterations for the ``standard'' SAGD under condition of uniform (for all $\xi$) $\cH$-Lipschitz continuity on the stochastic operator $\cG(\cdot, \xi)$.
However, this assumption impose significant limitations on the form of the stochastic operator in statistical learning applications and is violated in the simple case of unbounded (e.g., Gaussian) regressors $\phi_t$, etc. To conclude, in spite of these efforts, {to the best of our knowledge, the question of building an optimal stochastic approximation routine} of general smooth convex optimization under state-dependent noise assumption \eqref{cond:glr} (when only $\nabla f$ rather than $\cal G(\cdot, \xi)$ is Lipschitz continuous) {has not received a complete answer.}

\subsection{Contributions and organization} Given the state of affairs, this paper focuses on designing accelerated algorithms and providing sharp analysis for the general stochastic optimization problem \eqref{cp-ASGD} with state-dependent noise.
Our contribution is threefold.
\begin{enumerate}\item
We analyze the convergence rates of the generic (non-Euclidean) SAGD for solving stochastic optimization with state-dependent noise. We show that under condition~\eqref{cond:glr}, SAGD attains a convergence rate
\begin{align*}
\mathcal{O}\left(\tfrac{L R^2}{k^2} + \tfrac{\cL R^2}{km} + \tfrac{\sqrt{L \cL}R^2}{k \sqrt{m}}  + \sqrt{\tfrac{\sigma_*^2 R^2}{km}}\right)
\end{align*}
where $k$ is the number of iterations, $R$ is the initial distance to the optimal solution, and $m$ is the batch size. The terms in the above bound are optimal, except for the third term which has a sub-optimal dependence on the batch size $m$.
% the optimal sample complexity of $\mathcal{O}\left(\sqrt{\tfrac{L R^2}{\epsilon}} + \tfrac{\cL R^2 }{\epsilon} + \tfrac{R^2\sigma_*^2}{\epsilon^2}\right)$ when the batch size $m=1$.
As a consequence, in order to achieve the optimal iteration complexity $\mathcal{O}\big(\sqrt{L R^2/\epsilon}\big)$, SAGD requires a larger batch size, resulting in a sub-optimal sample complexity\footnote{In what follows we refer to the total number $N=N(\epsilon)$ of calls to the stochastic oracle which are necessary for the approximate solution $\wh x_k$ after $k=k(\epsilon)$ iterations and $N$ oracle calls to attain the expected (in)accuracy $\epsilon$, i.e., \begin{equation}\label{e-acc}\bbe[f(\wh x_k)]-f^*\leq \epsilon\end{equation}  as {\em sample} (or {\em information}) {\em $\epsilon$-complexity} of the method. We also call {\em iteration $\epsilon$-complexity} the minimal iteration count $k$ such that
\rf{e-acc} holds.} \[
\mathcal{O} \left(\sqrt{\tfrac{LR^2}{\epsilon}}  + {\tfrac{\sqrt{L}\cL  R^3}{ \epsilon^{3/2}}} + {\tfrac{R^2 \sigma_*^2}{ \epsilon^2}}\right).
\]
However, imposing the condition of boundness of the second moment of the Lipschitz constant of the gradient observation   $\cG(\cdot, \xi)$ allows to improve the second term in the above sample complexity bound to $\mathcal{O}\left(1/\epsilon\right)$. The corresponding iteration complexity bound $\mathcal{O}\big(\sqrt{L R^2/\epsilon}\big)$ is an improvement w.r.t. the bound $\mathcal{O}\big(\sqrt{\cH R^2/\epsilon}\big)$ of \cite{woodworth2021even}, and our assumption is more general than the ``uniform'' Lipschitz continuity assumption in \cite{woodworth2021even}. Moreover, unlike \cite{woodworth2021even}, our analysis does not require the feasible region $X$ to be a bounded set.

\item
% {To achieve the optimal convergence rate without imposing further assumptions, we consider an alternative accelerated SA method named stochastic gradient extrapolation method (SGE).}
{Under state-dependent noise assumption \rf{cond:glr} we analyze an alternative accelerated SA algorithm---stochastic gradient extrapolation
method (SGE).} The gradient extrapolation method was introduced in \cite{lan2018random} by exchanging the primal and the dual variables in a game interpretation of
Nesterov's accelerated gradient method (see \cite{lan2018random} and Chapter 3 and 4 of \cite{LanBook2020}). SGE uses the same sequence of points for both
gradient estimations and output solutions. {This appears to be a significant advantage of the SGE over the SAGD in the present setting and allows for direct ``compensation'' of the state-dependent noise term} by the suboptimality gap of approximate solutions.
As a result, SGE achieves the optimal convergence rate after $k$ iterations
\begin{align*}
    \mathcal{O}\left(\tfrac{L R^2}{k^2} + \tfrac{\cL R^2}{km} + \sqrt{\tfrac{\sigma_*^2 R^2}{km}}\right).
\end{align*}
Consequently, it attains the optimal
iteration complexity $\mathcal{O}\big(\sqrt{L R^2/\epsilon}\big)$
along with the optimal sample complexity $\mathcal{O}\big(\sqrt{L R^2/\epsilon} + \cL R^2 /{\epsilon} + R^2\sigma_*^2/\epsilon^2\big)$, as supported by lower bounds in \cite{nemyud:83, woodworth2021even}.
\item We propose a multi-stage algorithm with restarts for solving problems satisfying the quadratic growth condition stating that for some $\mu > 0$ and $x^*\in \feaReg$,\footnote{We suppose for convenience that in this case the optimal solution $x^*$ is unique. Note that \eqref{eq:general_quadractic_growth_1} can be seen as a relaxation of the strong convexity assumption, i.e., for any $x, y \in \feaReg$,
% for some $\mu > 0$,
$f(y) - f(x) - \langle \nabla f(x), y- x\rangle \geq \mu \|y-x\|^2/2.
$}
\begin{align}\label{eq:general_quadractic_growth_1}
    f(x) - f^*\geq \tfrac{\mu}{2}\|x-x^*\|^2, ~~\forall x \in \feaReg.
\end{align}
We show that the proposed procedure attains the optimal iteration complexity $\mathcal{O}\big( \sqrt{L/\mu}\log(1/{\epsilon})\big)$ and the optimal sample complexity $\mathcal{O}\big( \sqrt{L/\mu}\log(1/{\epsilon})  + \cL/\mu \log(1/\epsilon) + \sigma_*^2/(\mu \epsilon)\big)$ simultaneously.
Furthermore, we specify the multi-stage SGE to solve the sparse recovery problem. This is done by incorporating hard-thresholding of the approximate solution at the end of each algorithm stage to enforce the sparsity. The convergence results match the quadratic growth setting up to a multiplicative factor of the sparsity level $s$, with only logarithmic dependence on the dimension $n$. To the best of our knowledge, the corresponding convergence guarantees are new to the optimization literature and lead to extremely efficient algorithms for sparse recovery problems in the distributed setting.
\end{enumerate}
Remaining sections of the paper are organized as follows. In Section~\ref{sec_assump}, we formalize the general problem statement and introduce the mini-batch setup. In Section~\ref{sec_SAGD}, we study the SAGD method for solving the general convex problem with state-dependent noise.
Section~\ref{sec_SGE} introduces SGE and provides convergence guarantees in the general convex problem. We introduce  the multi-stage SGE for solving problems satisfying quadratic growth condition in Section~\ref{sec:quadratic_growth}. Section~\ref{sec_sparse} further extends the multi-stage SGE to be applied to the sparse recovery problem. Finally, in Section \ref{sec:numeric}, we present some results of a preliminary simulation study illustrating the numerical performance of the proposed algorithms in the high-dimensional setting of sparse recovery. Proofs of the statements are postponed \revision{till}{until} the appendix.

\subsection{Notation} For any $n \geq 1$, we use $[n]$ to denote the set of integers $\{1, ..., n \}$. For $x \in \bbr$, we let $(x)_+ = \max(x,0)$ and $(x)_- = \max(-x,0)$. In what follows, $\euc$ is a finite-dimensional real-vector (Euclidean) space. Given a norm $\|\cdot\|$ on $\euc$, the associated dual norm $\|\cdot\|_*$ is defined as $\|z\|_* := \sup\{\langle x, z\rangle : \|x\|\leq 1\}$.  We define $\omega:\euc\rightarrow \bbr$, the {\em distance generating function,} which is a continuously differentiable strongly convex function with modulus $1$. Without loss of generality, we assume that $\omega(x) \geq \omega(0) = 0$ and for some $\Omega \geq 1$,
\begin{align}\label{V_to_norm_general}
    \omega(x) \leq \tfrac{\Omega}{2} \|x\|^2,\quad \forall x \in \euc.
\end{align}
Ideally, we want $\Omega$ to be ``not too large''.\footnote{\textcolor{black}{Note that condition~\eqref{V_to_norm_general} is only required in Sections~\ref{sec:quadratic_growth} and \ref{sec_sparse}.}} Meanwhile, a desired distance generating function should be ``prox-friendly'', i.e., for any $a \in \euc$, the minimization problem
\begin{align*}
    \min_{x \in \feaReg} \{\langle a, x \rangle + \omega(x)\}
\end{align*}
can be easily solved. Note that when $\|\cdot\|$ is the Euclidean norm, \revision{}{we can set $\omega(x)=\|x\|_2^2/{2}$, and the corresponding $\Omega = 1$}. Another ``standard'' choice is $\|\cdot\| = \|\cdot\|_1$, $\|\cdot\|_*=\|\cdot\|_\infty$, and one can choose the distance generating function $\omega$ (cf. \cite{nesterov2013first})
\begin{align*}
    \omega(x) = \tfrac{1}{2} e \log n \cdot n^{(p-1)(2-p)/p}\|x\|_p^2, \quad p = 1 + \tfrac{1}{\log n};
\end{align*}
the corresponding $\Omega$ satisfies $\Omega\leq e^2 \log n$ in this case.

Given an initialization $x_0\in \feaReg$, we define the $x_0$-associated Bregman's divergence of $x,y\in \feaReg$ as
\begin{align*}
    V_{x_0} (x,y) = \omega(y-x_0) - \omega(x-x_0) - \langle \nabla \omega(x-x_0),y-x\rangle.
\end{align*}
Clearly, for any $y, x, x_0\in \feaReg$, we have
\begin{align}\label{bound_bregman}
     V_{x_0}(x_0,y) \leq \tfrac{\Omega}{2}\|y-x_0\|^2 \quad \text{and} \quad V_{x_0}(x, y) \geq \tfrac{1}{2}\|x-y\|^2.
\end{align}
We use the shorthand notation $V(x,y)$ for $V_{x_0}(x,y)$ when $x_0$ is clear in the content.

Unless stated otherwise, all relations between random variables are assumed to hold almost surely.
\section{Problem statement}\label{sec_assump}
We summarize below the setting of the stochastic optimization problem and our assumptions.
\subsection{Assumptions}
The problem under consideration is a stochastic optimization problem \eqref{cp-ASGD} with a convex and smooth objective function as given in \eqref{smoothness_ASGD}.
% Additionally, we assume the quadratic growth condition \eqref{eq:general_quadractic_growth} with $\mu \geq 0$, and the problem reduces to the general convex problem when $\mu = 0$.

We assume that stochastic oracle $\cG(\cdot, \cdot)$ is unbiased, i.e., satisfies \eqref{assump:bias_SO}. Furthermore, we consider the following ``state-dependent'' oracle noise assumption.

\begin{itemize}
    \item{[{\em State-dependent variance}]} We assume that for all $x\in X$ and some $\cL<\infty$,
    % {\color{red} Do we use the fact that $L\leq \cL$? Do we need this assumption?} \textcolor{red}{\textbf{Not really, I removed this one and change \eqref{complexity:SAGD_2_0} accordingly}}
\beq \label{assump:variance}
\bbe_{\xi_t}\big[\|\cG(x, \xi_t) - g({x})\|_*^2\big] \le \revision{}{\sigma^2}(x):=  \cL[f({x}) - f^* - \langle g(x^*), {x} - x^*\rangle] + \sigma_*^2, \quad t \in \mathbb{Z}_+, \tag{SN}
\eeq
for some $x^* \in \feaReg^*$. 
% and $\sigma_t(\cdot)$ is a real (nonrandom) function on $\feaReg$.
\end{itemize}
Assumption~\eqref{assump:variance} can be further weakened to \revision{}{$\sigma^2(x):= \cL[f(x) - f^*] + \sigma_*^2$}. In the case of unconditional minimizer $x^*\in \inter X$, the latter condition is clearly equivalent to \rf{assump:variance}.  In the case of $g(x^*)\neq 0$, utilizing the relaxed assumption results in convergence guarantees which depend explicitly on $f(x_0) - f^*$. We use \eqref{assump:variance} for the sake of convenience, the term $\langle g(x^*), x - x^*\rangle$ in the right-hand side when combined with smoothness of $f$ allow us to upper bound the variance of $\cG(x_0, \xi_t)$ at $x_0$ by the term proportional to $\|x_0-x^*\|^2$; see, e.g., \eqref{eq:GEM_theorem_result_general}.

\begin{figure}[ht]
    \centering
    \begin{tabular}{cc}
    \includegraphics[width=0.4\textwidth]{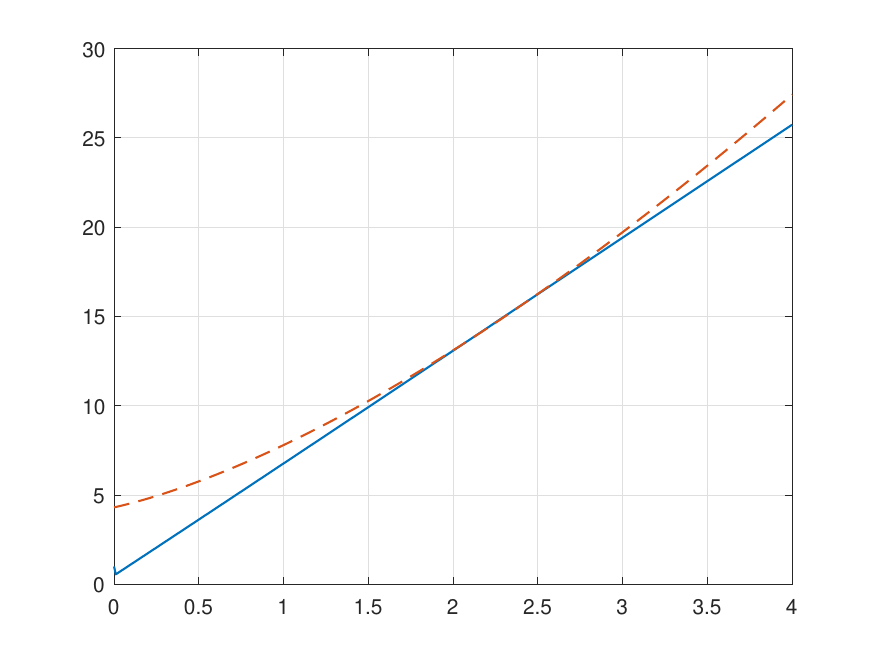}&\includegraphics[width=0.4\textwidth]{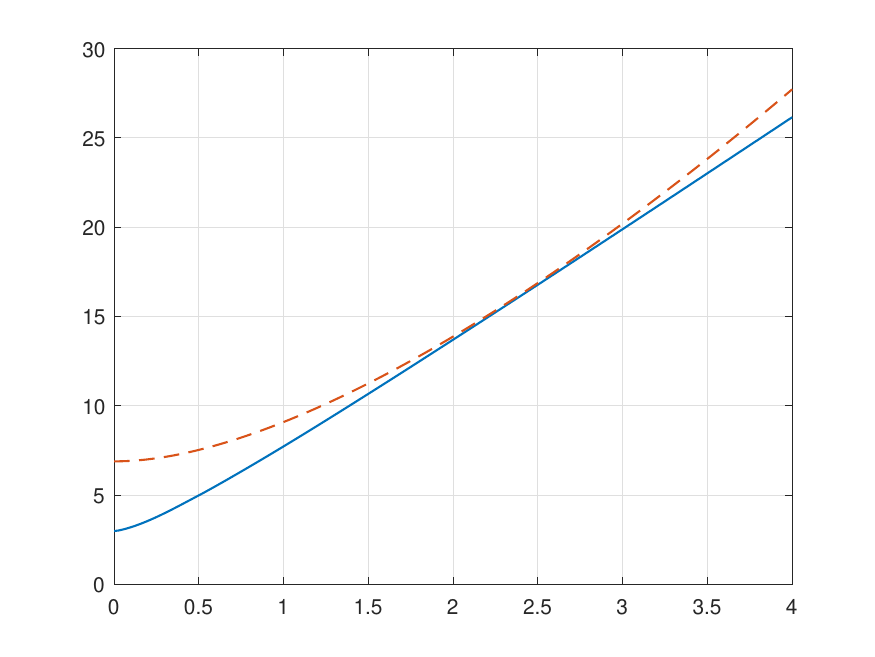}
    \end{tabular}
    \caption{Left plot: variance of the stochastic oracle $\cG_1(x,\xi)$ as function of $x$ (solid line) and upper bound $4\sigma^2(0)+3[f(x)-f(0)]$ (dashed line);
     right plot: variance of  $\cG_2(x,\xi)$ as function of $x$ (solid line) and upper bound $2.3\sigma^2(0)+3[f(x)-f(0)]$ (dashed line).}
    \label{fig:2plots}
\end{figure}

When proving the convergence rates for the stochastic variant of Nesterov's accelerated gradient descent method (SAGD; see Section~\ref{sec_SAGD}), we also consider the following assumption:
\begin{itemize}
    \item{[{\em Lipschitz continuous stochastic gradient}]} For each $\xi\in \Xi$, there exists a $\cK(\xi) > 0$, such that $\bbe_\xi[\cK(\xi)^2] < \infty$ and
\beq \label{assump:variance_1}
\|\cG(x, \xi) - \cG(y, \xi)\|_* \leq \cK(\xi) \|x-y\|, \quad \forall x, y \in X. \tag{LP}
\eeq
\end{itemize}
This assumption relaxes the assumption in \cite{woodworth2021even} that assumes $\cG(\cdot, \xi)$ is $\cH$-Lipschitz continuous for all $\xi \in \Xi$. For instance, for GLR model with Gaussian/sub-Gaussian regressors $\phi_t$, \eqref{assump:variance_1} holds, but the $\cH$-uniform Lipschitz continuous condition is violated. Nevertheless, \eqref{assump:variance_1} is not a necessary condition of Assumption~\eqref{assump:variance} since the latter one may hold in various situations of interest where $\cG(\cdot,\xi)$ is not Lipschitz (and even not continuous).
% ({\bf GL: What do we mean by smooth here?}).
In Figure \ref{fig:2plots} we present the plot of the expectation in the left-hand side of \rf{assump:variance} as a function of $x$ for two choices of scalar discontinuous gradient observation $\cG_1(x,\phi)=\phi u(\phi x)$ and $\cG_2(x,[\phi,\zeta])=\phi u(\phi x+\zeta)$ where $\phi$ and $\zeta$ are independent r.v. with Student $t_4$ distribution and $u(t)=(\tfrac{1}{2}+\sqrt{|t|})\mathrm{sign}(t)$.

Given the limitations of Assumption~\eqref{assump:variance_1}, we will only use it partially in Section~\ref{sec_SAGD} in order to improve the convergence rates of SAGD. In the following sections, we will propose an alternative accelerated algorithm called SGE that does not rely on Assumption~\eqref{assump:variance_1} but attains stronger convergence guarantees; see Section~\ref{sec_SGE} for more details.

\subsection{Mini-batch setup}
We consider the mini-batch approach widely used in practice. Specifically, we assume that at each search point \aic{$u_t$,}{} the stochastic oracle is called repeatedly, thus generating $m_t$ i.i.d. samples $\{\xi_{t,i}\}_{i=1}^{m_t}$, $m_t$ being the number of oracle calls.
{We define  filtration $\calF_t=\sigma\big(u_0,\xi_{0,1},...,\xi_{0,m_0},\xi_{1,1},...,\xi_{1,m_1},......,\xi_{t,1},...,\xi_{t,m_t}\big)$ and use the shorthand notation $\bbec{t}$ for the conditional expectation with respect to $\calF_t$.}

Given an $\cF_{t-1}$-measurable search point $u_{t}\in X$, we compute the unbiased estimator $G_t=G_t(u_t)$ of $g(u_t)$,\footnote{Here $u_t$ is a general place holder for a $\calF_{t-1}$-measurable search point. Note that, with a slight ambiguity of notation, search points $z_t$ and $\xavg_t$ in Algorithms~\ref{alg:SAGD} and \ref{adaptive_GEM} are $
\calF_{t}$-measurable. }
% We assume that mini-batch is allowed, where the mini-batch operator is calculated by
\begin{align*}
G_t = \frac{1}{m_t} \sum_{i=1}^{m_t} \cG( u_t, \xi_{t,i}).
\end{align*}
%Observe that
%\begin{itemize}
%\item random variables in $\{\xi_{t,i}\}_{i=1}^{m_t}$ are $\calF_t$-measurable.
%\item search points with index $t$ which are deterministic functions of $G_\tau (u_\tau),\tau\leq t-1$, %e.g., $x_t$ and $z_t$ in Algorithm~\ref{adaptive_GEM},
%are $\calF_{t-1}$-measurable.\footnote{Here $u_t$ is a general place holder for a $\calF_{t-1}$-measurable search point. With a slight ambiguity of notation, search points $z_t$ and $\xavg_t$ in Algorithms~\ref{alg:SAGD} and \ref{adaptive_GEM} are $
%\calF_{t}$-measurable. }
%\end{itemize}

% Here, $m_t$ denotes the mini-batch size, and ${\xi_{t,i}}_{i=1}^{m_t}$ represents the samples obtained from the stochastic oracle.
%We define  filtration $\calF_t=\sigma\big(u_0,\xi_{0,1},...,\xi_{0,m_0},\xi_{1,1},...,\xi_{1,m_1},......,\xi_{t,1},...,\xi_{t,m_t}\big)$, so that
%\begin{itemize}
%\item random variables in $\{\xi_{t,i}\}_{i=1}^{m_t}$ are $\calF_t$-measurable.
%\item search points with index $t$ which are deterministic functions of $G_\tau (u_\tau),\tau\leq t-1$, %e.g., $x_t$ and $z_t$ in Algorithm~\ref{adaptive_GEM},
%are $\calF_{t-1}$-measurable.\footnote{Here $u_t$ is a general place holder for a $\calF_{t-1}$-measurable search point. With a slight ambiguity of notation, search points $z_t$ and $\xavg_t$ in Algorithms~\ref{alg:SAGD} and \ref{adaptive_GEM} are $
%\calF_{t}$-measurable. }
%\end{itemize}
%We use the shorthand notation $\bbec{t}$ to denote the conditional expectation with respect to the filtration $\calF_t$.

Based on the state-dependent noise Assumption~\eqref{assump:variance}, we have the following characterization of the properties of the mini-batch estimator:
{\color{black}
\begin{lemma}\label{assump:expectation_revised}
Denote
\beq\label{vstar}
V_{x, x_0}^*(y)= \max_{z \in E} \left\{\langle y, z - x \rangle - V_{x_0}(x,z) \right\}
\eeq (we use the shorthand notation $V_{x}^*(y)$ when $x_0$ is clear in content). The mini-batch estimator $G_t$ satisfies for any \aic{$x_0 \in E$, $\calF_{t-1}$-measurable $x$,}{$x_0,x$, and $u\in X$} and $\gamma \in \bbr$,
\beq\label{eq:assum_minibatch1_revised}
%\bbe[\|\delta_t\|_*^2]\le
\aic{\bbec{t-1}
\big[V_{x, x_0}^*\big(c\cdot[G_t (u_t) - g (u_t)]\big)\big] \le \tfrac{c^2}{2m_t} \cdot \bbec{t-1}[\|\cG_t(u_t, \xi_{t, 1}) - g(u_t)\|_*^2],
}{
\bbe_{\{\xi_{t,i}\}_{i=1}^{m_t}}\big[V_{x, x_0}^*\big(\gamma[G_t (u) - g (u)]\big)\big] \le \tfrac{\gamma^2}{2m_t} \bbe_{\xi_{t, 1}}[\|\cG_t(u, \xi_{t, 1}) - g(u)\|_*^2],}
\eeq
Consequently, when Assumption \eqref{assump:variance} holds,
\beq\label{eq:assum_minibatch_revised}
\aic{
\bbec{t-1}\big[V_{x, x_0}^*\big(c[G_t (u_t) - g (u_t)]\big)\big] \le \tfrac{c^2}{2m_t}  \left\{\cL[f(u_t) - f^* - \langle g(x^*), u_t - x^*\rangle] + \sigma_*^2\right\}.
}{
\bbe_{\{\xi_{t,i}\}_{i=1}^{m_t}}\big[V_{x, x_0}^*\big(\gamma[G_t (u) - g (u)]\big)\big] \leq\tfrac{\gamma^2}{2m_t}  \left\{\cL[f(u) - f^* - \langle g(x^*), u - x^*\rangle] + \sigma_*^2\right\}.
}
\eeq
% , for some $x^*\in \feaReg^*$, where $\Omega$ is defined in \eqref{V_to_norm_general}.
\end{lemma}}

In the following section, we present accelerated algorithms equipped with mini-batches \revision{in order to achieve}{that achieve} the optimal iteration complexity for stochastic optimization with state-dependent noise.

\section{SAGD for general convex problem with state-dependent noise}\label{sec_SAGD}
Algorithm 1 describes a mini-batch variant of the standard stochastic accelerated gradient descent (SAGD) method. The SAGD method, also called accelerated stochastic approximation (AC-SA), maintains three sequences of points. Specifically, $\{z_t\}$ is the sequence of ``prox-centers'' for the prox-mapping updates \eqref{eq:ad_SAGD2}, and $\xmid_t$ and $\xavg_t$ are weighted averages of the past $z_t$'s;  $\{\xmid_t\}$ is the sequence of search points where stochastic gradients are estimated using mini-batches, and points $\{\xavg_t\}$ is the trajectory of approximate solutions (outputs) at each iteration.

\begin{algorithm}[H]
	\caption{Stochastic accelerated gradient descent method (SAGD)}
	\begin{algorithmic}
		\State {\bf Input:} initial point $z_0= \xavg_0$, nonnegative nonrandom parameters $\{\beta_t\}$ and $\{\eta_t\}$, and batch size $\{m_t\}$.
		% \State $G_0=G_{-1} = \tfrac{1}{m_0} \tsum_{i=1}^{m_0} \cG(z_0, \xi_{0,i})$.
		\For {$t =1, 2, \ldots,$}
		%\State
\begin{subequations}
\label{eq:ad_SAGD0}
\begin{align}
            \xmid_t &= (1-\beta_t)\xavg_{t-1}+\beta_t z_{t-1} \label{eq:ad_SAGD1}\\
            G_t &= \tfrac{1}{m_t} \tsum_{i=1}^{m_t} \cG(\xmid_t, \xi_{t,i}). \nn\\
			z_t &= \argmin_{z \in \feaReg} \{ \langle G_t, z \rangle+\eta_t V(z_{t-1}, z)
% 			+ \lambda_t V_{z_0}(z_0, x)
			\}, \label{eq:ad_SAGD2}\\
			\xavg_t &= (1-\beta_t) \xavg_{t-1} + \beta_t z_t,\label{eq:ad_SAGD3}
		\end{align}
\end{subequations}		\EndFor
	\end{algorithmic} \label{alg:SAGD}
\end{algorithm}

We start with the following characterization of the approximate solution $x_k$ of Algorithm \ref{alg:SAGD}.

\begin{theorem}\label{thm:SAGD_2}
Suppose Assumption~\eqref{assump:variance} is satisfied. Let the algorithmic parameters $\beta_t$ and $\eta_t$ satisfy for some $\theta_t \geq 0$,
\begin{subequations}
\label{SAGD_cond0}
\begin{align}
    \theta_t \beta_t \eta_t &\leq \theta_{t-1} \beta_{t-1} \eta_{t-1}, ~t=2,..., k \label{SAGD_cond_2}\\
    \eta_t &> L\beta_t, ~t=1,..., k. \label{SAGD_cond_3}
\end{align}
\end{subequations}
Furthermore, suppose that $\beta_1 = 1$ and
\begin{align}
    \theta_t(1-\beta_t)(1+ r_t \cL) \leq \theta_{t-1}, ~t = 2, ..., k. \label{SAGD_cond_1}
\end{align}
Then
\begin{align}
    \theta_k \bbe[f(\xavg_k) &- f^*] + \theta_k \beta_k \eta_k \bbe[V(z_k, x^*)]\nn\\
    & \leq \theta_1\beta_1\eta_1 V(z_0,x^*) + \sum_{t=1}^k \theta_t r_t \beta_t L \cL \bbe[V(z_{t-1}, x^*)] + \sum_{t=1}^k \theta_t r_t \sigma_*^2 \label{eq:thm_SAGD_1}
\end{align}
where
%\tli{$r_t := \tfrac{\beta_t \revision{\Omega}{\newOmega}}{2(\eta_t - L \beta_t) m_t}$\revision{}{ and $\newOmega$ is defined in Lemma~\ref{assump:expectation}}}
$r_t := \left[2(\eta_t - L \beta_t) m_t\right]^{-1}\beta_t$.
\end{theorem}
{\color{black}\begin{corollary}\label{cor:sagd2}
In the premise of Theorem \ref{thm:SAGD_2}, {let $k\geq 2$,}  and let $\theta_t = (t+1)(t+2)$, $\beta_t = \tfrac{3}{t+2}$, $m_t = m$, and $\eta_t = \tfrac{\eta}{t+1}$ with
$$\eta = \max\left\{4L, \tfrac{6 (k-1) \cL}{m}, \sqrt{\tfrac{9(k+{1})^2  L \cL}{m}}, \tfrac{ \sigma_*}{D}\sqrt{\tfrac{2(k+2)^3 }{3 m}}\right\}$$
{where $D>0$ is such that $V(z_0,x^*)\leq D ^2$.}
Then $\bbe[V(z_t,x^*)] \leq 3 D ^2$ for all $1\leq t\leq k$, and
\begin{align}\label{eq:thm_SAGD_2}
    \bbe[f(\xavg_k) -f^*] \leq \tfrac{12LD ^2}{(k+1)(k+2)}+ \tfrac{6  \cL D ^2}{(k+2) m} + \tfrac{18D ^2\sqrt{ L \cL}}{(k+1)\sqrt{m}} + \tfrac{4\sigma_*D\sqrt{2 }}{\sqrt{ (k+1)m}}.
\end{align}
\end{corollary}}
\paragraph{Remarks.}
Bound \rf{eq:thm_SAGD_2} of  Corollary \ref{cor:sagd2} allows us to establish the complexity bounds for the SAGD algorithm when solving a general convex problem with state-dependent noise. Let us first consider the case when $m=1$.
% \revision{In this case,  we can take $\Omega/m = 1$}{In this case, we have $\newOmega=1$,}
% \textcolor{red}{(It is clear that when $m=1$, Lemma 1 holds with $\Omega=1$.)}
The total number of iterations/oracle calls used by SAGD to find an $\epsilon$-optimal solution, i.e., $\wh x\in X$ such that $\bbe[f(\wh x) -f^*] \leq \epsilon$, is bounded by
\begin{align}\label{complexity:SAGD_2_0}
    \mathcal{O} \left(\sqrt{\tfrac{LD ^2}{\epsilon}} + \tfrac{\cL D ^2}{\epsilon} + \tfrac{\sqrt{L \cL}D^2 }{\epsilon}+ \tfrac{D ^2 \sigma_*^2}{\epsilon^2} \right).
\end{align}
%In this complexity bound, the three terms represent the deterministic error, the state-dependent stochastic error, and the state-independent stochastic error, respectively.
\revision{In the usual setting with $L =\mathcal{O}(\cL)$}{Considering the setting with $L =\mathcal{O}(\cL)$, as in the context of \cite{woodworth2021even}}, this upper bound matches the optimal sample complexity under Assumption~\eqref{assump:variance}, supported by \textcolor{black}{the lower bound in Theorem 4 of \cite{woodworth2021even}}.

In the mini-batch setting,
% \revision{}{(where $\newOmega = \Omega$)},
the bound in \eqref{eq:thm_SAGD_2} means that in order to achieve the optimal iteration complexity of $\mathcal{O}\left(\sqrt{LD ^2/\epsilon} \right)$, the batch size
$
m \geq \max\left\{1, \tfrac{k \cL}{L}, \tfrac{k^2  \cL}{L} , \tfrac{k^3  \sigma_*^2}{D ^2 L^2} \right\}=\max\left\{1, \tfrac{k^2  \cL}{L} , \tfrac{k^3 \sigma_*^2}{D ^2 L^2} \right\}
$
is needed. Consequently, the total sample complexity is bounded by
\tli{}{
\begin{align}\label{complexity:SAGD_2}
\mathcal{O} \left(\sqrt{\tfrac{LD ^2}{\epsilon}} + \tfrac{\sqrt{L}\cL  D ^3}{\epsilon^{3/2}} + \tfrac{D ^2 \sigma_*^2}{\epsilon^2} \right)
\end{align}}
in this situation.
When comparing \eqref{complexity:SAGD_2} to \eqref{complexity:SAGD_2_0}, we observe that the second term in \eqref{complexity:SAGD_2} is sub-optimal. This implies that, based on our analysis, SAGD does not simultaneously achieve the optimal iteration complexity and sample complexity under Assumption~\eqref{assump:variance}.

{The analysis on the convergence of SAGD in the state-dependent noise setting reveals the ``bottleneck'': in the recursion \rf{eq:ad_SAGD0}}, different points $y_t$ and $x_t$ are used for gradient estimations and output solutions. Improving the convergence rates in {Corollary \ref{cor:sagd2}} requires better control of the objective value at the points of stochastic gradient estimation. {This can be \revision{achived}{achieved} by imposing the Lipschitz regularity assumption in \eqref{assump:variance_1} on stochastic gradients.}
\tli{}{
\begin{lemma}\label{assump:expectation_2}
Suppose Assumptions~\eqref{assump:variance} and \eqref{assump:variance_1} hold. Let $\{\xavg_t\}$, $\{\xmid_t\}$ and $\{z_t\}$ be generated by Algorithm~\ref{alg:SAGD}.  Then for any $\gamma \in \bbr$,
\begin{align*}
    \bbec{t-1}\big[V^*_{z_{t-1}}(\gamma[G_t - g(\xmid_t)])\big] \leq {\tfrac{\gamma^2}{m_t}}\left\{{\bar \cK}^2\beta_t^2\|\xavg_{t-1} - z_{t-1} \|^2 + 3\cL[f(\xavg_{t-1}) - f^* - \langle g(x^*), \xavg_{t-1} - x^*\rangle] + 3\sigma_*^2\right\},
\end{align*}
where $V^*_x(\cdot)$ is as defined in \rf{vstar} and
\begin{align}
    {\bar \cK}:=\left(3\bbe_{\xi_{t,1}}[\cK(\xi_{t,1})^2] + 3L^2\right)^{1/2}.\label{K_constant}
\end{align}
% ({\bf GL: Please list the definition of $\cK$ in a separate line, and label it. Please refer to this definition later (e.g., in Theorem 2 or the discussion. Otherwise, it is difficult to find this definition}).
\end{lemma}}
We have the following analog of Theorem \ref{thm:SAGD_2} under Assumptions~\eqref{assump:variance} and \eqref{assump:variance_1}.
\begin{theorem}\label{thm:SAGD_1}
Let Assumptions~\eqref{assump:variance} and~\eqref{assump:variance_1} hold. Let $\beta_t,\eta_t$ satisfy \rf{SAGD_cond0} for some $\theta_t\geq 0$, and let also  $\beta_1 = 1$ and
\begin{align}\label{eq:SAGD_5}
    \theta_t (1-\beta_t + 3r_t \cL) \leq \theta_{t-1}, ~t=2,...,k.
\end{align}
Then
\begin{align}
\theta_k &\bbe[f(\xavg_k) - f^*] + \theta_k \beta_k \eta_k \bbe[V(z_k, x^*)]\nn\\
    & \leq \theta_1\beta_1\eta_1 V(z_0,x^*) +  \tfrac{3\theta_1 r_1 \cL L}{2} \|z_0-x^*\|^2+ \sum_{t=1}^k \theta_t r_t \beta_t^2 \revision{\cK^2}{\bar \cK^2} \bbe\|x_{t-1} -z_{t-1}\|^2 + \sum_{t=1}^k 3\theta_t r_t \sigma_*^2
\label{eq:thm_SAGD_3}
\end{align}
where \revision{$\cK$}{$\bar \cK$} is defined in \eqref{K_constant} and $r_t := \tfrac{\beta_t }{2(\eta_t - L \beta_t) m_t}$.
\end{theorem}
\begin{corollary}\label{cor:sagd22}
In the premise of Theorem \ref{thm:SAGD_1}, suppose that $k\geq 2$, $V(z_0,x^*)\leq D ^2$, and let $\beta_t = \tfrac{3}{t+2}$, $m_t = m$ and $\eta_t = \tfrac{\eta}{t+1}$ with
$$\tli{}{\eta = \max\left\{4L, \tfrac{18\tli{\revision{\Omega}{\newOmega}}{} (k+1) \cL}{m}, 12\sqrt{\tfrac{k \tli{\revision{\Omega}{\newOmega}}{}{\bar \cK^2}}{m}}, \tfrac{ \sigma_*}{D}\sqrt{\tfrac{2(k+2)^3 }{ m}}\right\}.}$$
Then $\bbe[V(z_t,z_0)] \leq 3 D ^2$, $t=1,...,k$, and
\begin{align}\label{eq:thm_SAGD_4}
\tli{}{
    \bbe[f(\xavg_k) -f^*] \leq \tfrac{13LD ^2}{(k+1)(k+2)}+ \tfrac{54 \tli{\revision{\Omega}{\newOmega}}{} \cL D ^2}{(k+2) m} + \tfrac{72 {\bar \cK} D ^2}{(k+2)\sqrt{(k+1)m}}+ \tfrac{4\sigma_*D\sqrt{6}}{\sqrt{(k+1)m}}.}
\end{align}
\end{corollary}
\paragraph{Remarks.}
Let $m=1$\tli{(and \revision{$\Omega/m = 1$}{$\overline{\Omega}=1$} in Lemmas~\ref{assump:expectation} and \ref{assump:expectation_2})}{}. By \eqref{eq:thm_SAGD_4}, the iteration/sample complexity of the SAGD is bounded by
\begin{align*}
\mathcal{O} \left\{\sqrt{\tfrac{LD ^2}{\epsilon}} + \tfrac{\cL D ^2}{\epsilon} + \left(\tfrac{\revision{\cK}{\bar \cK} D ^2}{\epsilon}\right)^{\frac{2}{3}} + \tfrac{D ^2 \sigma_*^2}{\epsilon^2} \right\}.
\end{align*}
% {\crd Except for the term $\mathcal{O}\left(\left(\cK D ^2/\epsilon\right)^{\frac{2}{3}}\right)$, other three terms are optimal.
% What does this mean? I suggest that we remove this and formulate differently. See the "black" sentence below. }
\textcolor{black}{When $\epsilon = \mathcal{O}(L^3 D^2/\revision{\cK^2}{\bar \cK^2})$, this sample complexity is optimal.}

In the mini-batch setting\tli{(when \revision{$\Omega/m = \Omega$}{$\overline{\Omega}=\Omega$} in Lemmas~\ref{assump:expectation} and \ref{assump:expectation_2})}{}, in order to obtain the optimal iteration complexity $\mathcal{O}\left(\sqrt{LD ^2/\epsilon} \right)$, SAGD batch size $m$ should be of the order of
{\begin{align*}
\max\left\{ 1, \tfrac{k  \cL}{L}, \tfrac{k \revision{\cK^2}{\bar \cK^2}}{L^2}, \tfrac{k^3  \sigma_*^2}{L^2 D ^2}\right\}.
\end{align*}}
As a result, the corresponding sample complexity becomes
{\color{black}
\begin{align}\label{complexity:SAGD_1}
    \mathcal{O} \left\{\sqrt{\tfrac{LD ^2}{\epsilon}} + \tfrac{\cL D ^2}{\epsilon} + \tfrac{\bar \cK^2 D ^2}{L \epsilon} + \tfrac{ D ^2 \sigma_*^2}{\epsilon^2} \right\}.
\end{align}}
When $\revision{\cK^2}{\bar \cK^2}=\mathcal{O}(L\cL)$, this complexity bound is optimal (cf. Theorem 4 of \cite{woodworth2021even}). The result of Corollary \ref{cor:sagd22} refines the corresponding statement of \cite{woodworth2021even} in three aspects. First, the corresponding iteration complexity bound $\mathcal{O}\left(\sqrt{LD ^2/\epsilon} \right)$  is stated in terms of the Lipschitz constant of the expected gradient. Second, it relies upon Assumption \rf{assump:variance_1} which is significantly weaker than the assumption of uniform Lipschitz continuity of the stochastic gradient observation $\cG(x,\cdot)$ used in \cite{woodworth2021even}. Third, unlike \cite{woodworth2021even}, our analysis does not require the feasible region $X$ to be bounded.
%
%
%In the next section, we consider an alternative accelerated algorithm that attains the optimal convergence guarantees in terms of both iteration complexity and sample complexity without Assumption~\eqref{assump:variance_1}.

\section{SGE for general convex problem with state-dependent noise}\label{sec_SGE}
In this section, we discuss an alternative acceleration scheme  to solve the general smooth convex problem with state-dependent noise which we refer to as stochastic gradient extrapolation (SGE).
% We introduce the SGE algorithm in Section~\ref{sec: SGE_general_convex}, and we study the convergence guarantee of SGE in expectation and in high probability in Section~\ref{sec:SGE_expectation} and \ref{sec_hpb}, respectively.
% \subsection{The stochastic gradient extrapolation method}\label{sec: SGE_general_convex}
SGE (Algorithm~\ref{adaptive_GEM}) is a variant of the gradient extrapolation method proposed in \cite{lan2018random}.
We consider here the mini-batch version of the routine.
\begin{algorithm}[H]
	\caption{Stochastic gradient extrapolation method (SGE)}
	\begin{algorithmic}
		\State {\bf Input:} initial point $x_0= z_0$, nonnegative parameters $\{\alpha_t\}$, $\{\eta_t\}$ and $\{\beta_t\}$, and batch size $\{m_t\}$.
		\State \tli{$G_0 =G_{-1} = \tfrac{1}{m_0} \tsum_{i=1}^{m_0} \cG(x_0, \xi_{0,i})$.}{$x_0 = x_{-1}$.}
		\For {$t =1, 2, \ldots,$}
		%\State
\begin{subequations}
\label{eq:GGE}\aic{
\begin{align}
			\widetilde G_t &= G_{t-1}(x_{t-1}) + \alpha_t (G_{t-1}(x_{t-1}) - G_{t-1}(x_{t-2})), \label{eq:GE}\\
    G_t(\cdot) &= \tfrac{1}{m_{t-1}} \tsum_{i=1}^{m_{t-1}} \cG(\cdot, \xi_{t-1,i}),\\
			z_t &= \argmin_{x \in \feaReg} \{ \langle \widetilde G_t, x \rangle
% 			+ \mu \alpha_t V(\bar x_{t-1}, x)
			%+ \mu \alpha_t [V(\bar x_{t-1}, x) - V(\bar x_{t-2}, x)]  \nn\\
			% & \quad \quad \quad \quad
			+\eta_t V_{x_0}(z_{t-1}, x)
% 			+ \lambda_t V_{x_0}(x_0, x)
			\}, \label{eq:ad_SGE2}\\
			x_t& = {(1-\beta_t)x_{t-1} + \beta_t z_{t}}.\label{eq:ad_SGE3}
			% G_t(\cdot) &= \tfrac{1}{m_t} \tsum_{i=1}^{m_t} \cG(\cdot, \xi_{t,i}). \nn
		\end{align}}{
\begin{align}
			\widetilde G_t &= G_{t-1}(x_{t-1}) + \alpha_t (G_{t-1}(x_{t-1}) - G_{t-1}(x_{t-2})), \label{eq:GE}\\
 			z_t &= \argmin_{x \in \feaReg} \{ \langle \widetilde G_t, x \rangle
% 			+ \mu \alpha_t V(\bar x_{t-1}, x)
			%+ \mu \alpha_t [V(\bar x_{t-1}, x) - V(\bar x_{t-2}, x)]  \nn\\
			% & \quad \quad \quad \quad
			+\eta_t V_{x_0}(z_{t-1}, x)
% 			+ \lambda_t V_{x_0}(x_0, x)
			\}, \label{eq:ad_SGE2}\\
			x_t& = {(1-\beta_t)x_{t-1} + \beta_t z_{t}}.\label{eq:ad_SGE3}
			% G_t(\cdot) &= \tfrac{1}{m_t} \tsum_{i=1}^{m_t} \cG(\cdot, \xi_{t,i}). \nn
		\end{align}}
\end{subequations}
where
\beq
\tli{}{
   G_s(\cdot)= \tfrac{1}{m_{s}} \tsum_{i=1}^{m_{s}} \cG(\cdot, \xi_{s,i}).\label{eq:GE_new}}
\eeq
		\EndFor
	\end{algorithmic} \label{adaptive_GEM}
\end{algorithm}
\paragraph{Remarks.} {The} basic {iterative} scheme {presented in Algorithm~\ref{adaptive_GEM}} is conceptually simple. It involves two sequences of search points $\{z_t\}$ and $\{x_t\}$,  the latter being weighted averages of the former. Note that both $x_t$ and $z_t$ are $\calF_{t-1}$-measurable. Notably, the stochastic gradients are estimated at the search points $\{x_t\}$, which are also the approximate output solutions generated by the algorithm at each iteration.
This property brings benefits for dealing with state-dependent noise {of} the gradient estimation over the stochastic accelerated gradient descent (SAGD) method, where the output and gradient estimation use different sequences.

\tli{}{Notice that in step \eqref{eq:GE} we utilize random noises $\xi_{t-1,1},...,\xi_{t-1,m_{t-1}}$ to compute gradient estimates $G_{t-1}(x_{t-1})$ and $G_{t-1}(x_{t-2})$ at search points $x_{t-2}$ and $x_{t-1}$. In other words, it requires using two series of stochastic gradient computations at each search point $x_{s}$ of the method. This is indeed allowed in classical (gray-box) applications like the GLR problem introduced in Section~\ref{sec_intro}, and it ensures that both $\bbec{t-2}[G_{t-1}(x_{t-1})]=g(x_{t-1})$ and $\bbec{t-2}[G_{t-1}(x_{t-2})]=g(x_{t-2})$ hold. However, in the ``black-box'' setting where only the value of $\cG(x,\xi)$ may be computed at a search point $x$, this will require performing two series of queries to the stochastic oracle at each search point to ensure the above-mentioned conditions.}

The special relationship between SGE and SAGD merits an explanation. It has been discussed in detail in \cite[Section 3]{lan2018random} and  \cite[Section 5.2]{LanBook2020} in the deterministic setting. For the sake of completeness, we summarize the corresponding argument here.

For the sake of simplicity, let us consider the problem of unconstrained minimization
\[
\min_x f(x)
\] where $f:\,E\to \bbr$ is strictly convex and continuously differentiable. For $\varsigma\in E$ let us denote
\[
\varphi(\varsigma)=\max_x\big\{\langle x, \varsigma \rangle - f(x)\big\},
\]
so that $\varphi:\,E\to \bbr$ is strictly convex and continuously differentiable on $E$. Then $f$ has the Fenchel representation
\[f(x)=\max_{\varsigma} \big\{\underbrace{\langle x, \varsigma \rangle - \varphi(\varsigma)}_{=:F(x,\varsigma)}\big\},\;\; x\in E,
\]
and we can reformulate the original minimization problem as a saddle point problem:
\[
% f^* := \min_{x \in X} \left\{ \max_{g \in {\cal G} } \{\langle x, g \rangle - J_f(g)\}+ \mu \, \w(x) \right\}.
f^* := \min_{x} \left\{ \max_\varsigma F(x,\varsigma)\right\}.
\]
Let us define the Bregman divergence associated with $\varphi$ according to
\[
W_f(\chi,\varsigma)= \varphi(\varsigma)- [\varphi(\chi) + \langle \nabla\varphi(\chi), \varsigma -\chi \rangle],\;\;\chi,\varsigma\in E,\]
and the (generalized) prox-mapping
\begin{align}\label{gendprox}
\argmax_{\varsigma} \Big\{
F(z,\varsigma)-\tau W_f(\chi, \varsigma)\Big\},\;\;  z,\chi\in E.
\end{align}
As shown in Lemma 1 of \cite{lan2015optimal} (cf. also Lemma 3.6 of \cite{LanBook2020}),
%the generalized  prox-mapping associated with $W_f$---
the maximizer of \rf{gendprox} is the value of $\nabla f$ at certain $x\in E $, specifically,
\[
\nabla f(x) = \argmax_\varsigma \left\{F(z,\varsigma) - \tau W_f(\chi, \varsigma)\right\}
\]
where $x = [z + \tau \nabla\varphi(\chi)] / (1 + \tau)$.
\par Using the above observation, we can rewrite the corresponding deterministic version of the SGE recursion
in a primal-dual form. It is initialized with $(\varsigma_{-1}, \varsigma_0)$ and $x_0=\nabla \varphi(\varsigma_0)$, with the updates $(x_t,\varsigma_t)$ computed according to
\begin{subequations}
\label{eq:dsge}
\begin{align}
\tilde\varsigma_{t}&= \varsigma_{t-1} + \alpha_t (\varsigma_{t-1}-\varsigma_{t-2}), \label{def_alt_1}\\
z_{t} &=  \argmin_{x}\left\{\langle \tilde\varsigma_t,x\rangle +\eta_t V_{x_0}(z_{t-1},x) \right\}, \label{def_alt_2}\\
%&\mbox{Define $G_t$ as a stochastic observation of}\nn\\
\varsigma_{t} &= \argmax_\varsigma \left\{
F(z_t,\varsigma)- \tau_t W_f(\varsigma_{t-1}, \varsigma)\right\}.
%\Pr_{{\cal G}} (- x^t, g^{t-1}, \tau_t),
\label{def_alt_3}
\end{align}
\end{subequations}
Because $\nabla \varphi(\varsigma_{t-1})=x_{t-1}$, by the above, $\varsigma_t = \nabla f(x_t)$ with $x_t = (z_t + \tau_t x_{t-1}) / (1 + \tau_t)$ which is the definition of $x_t$ in \eqref{eq:ad_SGE3} with $\beta_t = 1/ (1+\tau_t)$.
The corresponding stochastic iteration \rf{eq:GGE} is obtained from \rf{eq:dsge} by replacing $\varsigma_t$ with its estimation $G_t$---the mean of $m_t$ stochastic gradients $\cG(x_t, \xi_{t,i})$. Similarly (see \cite[Section 2.2]{lan2015optimal} and \cite[Section 3.4]{LanBook2020}), one can show that SAGD iteration
can be viewed as a specific stochastic version of
the following (deterministic) primal-dual update:
\begin{subequations}
\label{eq:dsagd}
\begin{align}
\tilde z_{t}&= z_{t-1} + \alpha_t (z_{t-1}-z_{t-2}), \label{def_txt}\\
%&\mbox{Set $G_t$ to be a stochastic estimator of}\nn\\
\varsigma_{t} &= \argmax_{\varsigma} \left\{
F(\tilde z_t,\varsigma)- \tau_t W_f(\varsigma_{t-1}, \varsigma)\right\},
%\Pr_{{\cal G}} (-\tilde x^t, g^{t-1}, \tau_t),
\label{def_yt}\\
z_{t} &=  \argmin_{x}\left\{\langle \varsigma_t,x\rangle+\eta_t V_{x_0}(z_{t-1},x) \right\}.
%\Pr_{X} (g^t, x^{t-1}, \eta_t).
\label{def_xt_GEM}
 \end{align}
 \end{subequations}
Recursion \rf{eq:dsge} %\eqref{def_alt_1}--\eqref{def_alt_3}
can be seen as a dual version of \rf{eq:dsagd}; %\eqref{def_txt}--\eqref{def_xt_GEM},
%in which dual variables replace the primal ones in the equations \eqref{def_txt}--\eqref{def_xt_GEM}.
the principal difference between the two resides in the extrapolation step which is performed in
the dual space in \rf{def_alt_1} and in the primal space in \rf{def_txt}.

We now establish the convergence guarantees of SGE method in expectation, i.e., $\bbe[ f(x_k) - f^*]$.
We should stress here that the convergence analysis of
SGE under Assumption~\eqref{assump:variance}
is much more involved than the ones for its basic scheme in
\cite{lan2018random}, and the SAGD method in Section \ref{sec_SAGD}.
Therefore, the details are deferred to the appendix.

\begin{theorem}\label{thm:main}
Suppose Assumption \eqref{assump:variance} holds. Assume that the parameters $\{\alpha_t\}$, $\{\eta_t\}$ and $\{\beta_t\}$
of Algorithm~\ref{adaptive_GEM} satisfy \revision{}{and a nonnegative sequence $\{\theta_t\}$} satisfy
\begin{subequations}
\label{eq:def_rellb}
\begin{align}
\theta_{t-1} &= \alpha_t \theta_t, ~~\eta_t \le \alpha_t \eta_{t-1} , ~~ t= 2, \ldots,k \label{eq:def_eta_relb}\\
% % \Leftrightarrow  \eta_1\eta_2 \geq 25 L^2\alpha_2,\label{eq:def_alpha}
% \label{eq:def_Lip_1}\\
\tfrac{\eta_t{(1-\beta_t)}}{\alpha_{t}} &\ge {4}L\beta_t,
% \Leftrightarrow \eta_t \tau_{t-1} \ge 5 L \alpha_t,
\quad t= 3, \ldots, k \label{eq:def_Lip_relb}\\
 \tfrac{\eta_1\eta_2}{\alpha_2} &\geq  {16}L^2, ~~\eta_k (1-\beta_k)  \ge L\beta_k, \label{eq:def_Lip_rel1b}
\end{align}
\end{subequations}
Denote
\beq \label{eq:def_epsilon_t_adaptive}
q_t= {\tfrac{ \theta_{t+1} (1 + \alpha_{t+1}^2)}{\eta_{t+1}} +  \tfrac{ \theta_{t+2} \alpha_{t+2}^2  }{\eta_{t+2}}} \quad \text{and}\quad
\epsilon_t= {\tfrac{2q_t}{m_t}}.
\eeq
If $\beta_1=1$ and the parameters satisfy
\beq \label{eq:batchsizeBt_adaptiveb}
\tfrac{\theta_t(1-\beta_t)}{\beta_t}+\cL\epsilon_{t-1}\leq \tfrac{\theta_{t-1}}{\beta_{t-1}}, ~ t\ge 2,
\eeq
then
\begin{align}
\tfrac{\theta_k}{\beta_k}\bbe[ f(x_k) - f^*]
  \le \theta_1  \eta_1 V(x_0, x^*) +   \tfrac{\cL\epsilon_0 L}{2}\|x_0-x^*\|^2
 + \sigma_*^2 \sum_{t=0}^{k-1} \epsilon_t . \label{eq:GEM_theorem_result_general}
\end{align}
\end{theorem}

We now specify a particular stepsize policy in order to establish the convergence guarantee of SGE.
\begin{corollary}\label{coro_non_strongly_convex}
Let
$
% \label{stepsize_non_strong_convex_2}
	\theta_t = t, ~\alpha_t = \tfrac{t-1}{t},   ~m_t= m,
	$
	$\beta_t=\tfrac{3}{t+2}$, and $\eta_t = \tfrac{\eta}{t}$ for $\eta >0$.
Suppose that $V(x_0,x^*)\leq D ^2$ and that
\begin{align*}
\tli{}{
    \eta = \max\left\{\tli{30}{24}L, \tfrac{\tli{30 \revision{\Omega}{\newOmega}}{18} (k+2) \cL}{m}, \tfrac{\sigma_*}{D}\sqrt{\tfrac{2 (k+1)^3 }{m}}\right\}.}
\end{align*}
Then	\begin{align}\label{stepsize_non_strong_convex_2_4}
\tli{}{
	    \bbe[ f(x_k) - f^*] \leq  \tfrac{\tli{91}{73} L D ^2}{k(k+2)} + \tfrac{\tli{90\revision{\Omega}{\newOmega}}{54}  \cL D ^2}{m k} +  \tfrac{6\sigma_* D\sqrt{2}}{\sqrt{m k }}.}
	\end{align}

\end{corollary}
% ({\bf GL:Please add some explanation the difference the two results in the corollary, and why the second one is needed.})
\paragraph{Remark.}
\aic{The bounds of Corollary \ref{coro_non_strongly_convex} merit some comments.} Observe first that SGE achieves the optimal sample complexity
\begin{align}\label{complexity_expectation}
	\mathcal{O}\bigg\{ \sqrt{\tfrac{L D ^2 }{\epsilon}} +\tfrac{\cL  D ^2 }{\epsilon} + \tfrac{ D ^2\sigma_*^2}{\epsilon^2}\bigg\}
\end{align}
in the case of $m=1$.
In the mini-batch setting\tli{\revision{}{(where $\newOmega = \Omega$)}}{}, by setting the batch size of $m\geq \max\left\{1, \tfrac{\tli{\Omega}{} \cL k}{L}, \tfrac{\tli{\Omega}{} k^3 \sigma_*^2}{D ^2 L^2}\right\}$, the iteration complexity of the algorithm
% , i.e., $\bar x \in \feaReg$ such that $\bbe[f(\bar x) - f^*]\leq \epsilon$,
is bounded by $\O\left(\sqrt{L D ^2 /{\epsilon}}\right)$ and the overall sample complexity is
\begin{align}\label{complexity_expectation_3}
\tli{}{
	\mathcal{O}\bigg\{ \sqrt{\tfrac{L D ^2 }{\epsilon}} +\tfrac{\cL \tli{\Omega}{} D ^2 }{\epsilon} + \tfrac{\tli{\Omega}{} D ^2\sigma_*^2}{\epsilon^2}\bigg\}.}
\end{align}
 Notably, {the SGE attains the optimal iteration and sample complexity bounds.}

{In particular,} when $\|x_0 - x^*\|\leq R$ and {$V(x_0, x^*) \leq D^2:=\tfrac{\Omega}{2}\|x_0-x^*\|^2$, setting}
 \begin{align}\label{stepsize_non_strong_convex_2_3}
 \tli{}{
	\eta = \max\left\{\tli{30}{24}L, \tfrac{\tli{30 \revision{\Omega}{\newOmega}}{18} (k+2) \cL}{m}, \aic{\tfrac{2\sigma_*}{D}\sqrt{\tfrac{ (k+1)^3}{m}}}{\tfrac{2\sigma_*}{R}\sqrt{\tfrac{(k+1)^3}{\Omega m}}}\right\},}
	\end{align}
we obtain
 \begin{align}\label{stepsize_non_strong_convex_2_5}
 \tli{}{
	    \bbe[ f(x_k) - f^*] \leq  \tfrac{\tli{91}{73} L \Omega R^2}{2k(k+2)} + \tfrac{\tli{45\revision{\Omega^2}{\Omega\newOmega}}{27\Omega}  \cL R^2}{m k} +  {6\sigma_*R}\sqrt{\tfrac{\Omega}{m k }}.}
	\end{align}
This bound will be used in the proof of the multi-stage SGE method in the next section.
\aic{
Because $V(x_0, x^*) \leq \tfrac{\Omega}{2}\|x_0-x^*\|^2$ (cf. \eqref{bound_bregman}), \eqref{stepsize_non_strong_convex_2_5} states the convergence rate under condition $\|x_0-x^*\|\leq R$. This bound will be used in the proof of the multi-stage SGE method in the next section.
}{}%
\section{SGE for convex problem with quadratic growth condition}\label{sec:quadratic_growth}
In this section, we consider the problem setting in which the smooth objective $f$ in \rf{cp-ASGD} satisfies the quadratic growth condition (cf. \eqref{eq:general_quadractic_growth_1}), i.e.,
 when for some $\mu>0$ and \revision{$x_*\in \feaReg$}{$x^*\in \feaReg$}
\begin{align}\label{eq:general_quadractic_growth}
    f(x) - f^*\geq \tfrac{\mu}{2}\|x-x^*\|^2, ~~\forall x \in \feaReg.
\end{align}

We propose a multi-stage routine with restarts that utilizes Algorithm~\ref{adaptive_GEM} as a working horse.

% Before we formally introducing the algorithm, we introduce some additional notation. Given an initialization $x_0\in \feaReg$, we define the $x_0$-associated Bregman's distance between $x,y\in \feaReg$ as
% \begin{align*}
%     V_{x_0} (x,y) = \omega(y-x_0) - \omega(x-x_0) - \langle \nabla \omega(x-x_0),y-x\rangle.
% \end{align*}
% Clearly, for any $y, x, x_0\in \feaReg$, we have
% \begin{align}\label{bound_bregman}
%      V_{x_0}(x_0,y) \leq \tfrac{\Omega}{2}\|y-x_0\|^2 \quad \text{and} \quad V_{x_0}(x, y) \geq \tfrac{1}{2}\|x-y\|^2.
% \end{align}
% We will use the shorthand notation $V(x,y)$ for $V_{x_0}(x,y)$ when $x_0$ is clear in the content.
% \subsection{Convergence in expectation}\label{sec_strongly_convex_expectation}

\begin{algorithm}[H]
	\caption{Multi-stage stochastic gradient extrapolation method}
	\begin{algorithmic}
		\State {\bf Input:} initial point $\stx^0 \in \feaReg$. Let $R_0$ be a positive real.% which satisfies $\| \stx^0-x^*\|^2\leq R_0^2$.
		\For {$k =1, 2, \ldots, K$}
		\State (a) Set $R_k=R_02^{-k/2}$, $N=\left\lceil 10 \sqrt{\tfrac{2 \Omega L}{\mu}}\,\right\rceil$. Run $N$ iterations of SGE (Algorithm~\ref{adaptive_GEM}) with $x_0 = z_0 = \stx^{k-1}$ and
		\begin{align*}
		    \theta_t &= t, ~~\alpha_t = \tfrac{t-1}{t}, ~~{\beta_t=\tfrac{3}{t+2}},  ~~ \eta_t = \tfrac{\eta}{t},~~ t= 1,...,N\nn\\
      \eta &= \max\left\{{24}L, \tfrac{{18} ({N}+2) \cL}{{m^k}}, \tfrac{\sigma_*}{R_k}\sqrt{\tfrac{{2} ({N}+1)^3}{\Omega m^k}}\right\},\nn\\
		    m^k & =  \max\left\{1, \left\lceil \tfrac{{3} \cL ({N}+2)}{ L} \right\rceil , \left\lceil \tfrac{{8}  {N} ({N}+2)^2 \sigma_*^2}{{9 \Omega}L^2R_k^2} \right\rceil\right\}\aic{,~~ t= 0,...,{N}}.
	\end{align*}
	\State (b) Set $\stx^k = x_{N}$, where $x_{N}$ is the SGE solution obtained in Step (a).
		\EndFor
	\end{algorithmic} \label{adaptive_GEM_epoch}
\end{algorithm}

 Algorithm~\ref{adaptive_GEM_epoch} has a simple structure. {Each stage $k$ of the routine, consists of $N_k=N$ iterations of the SGE with initial condition $x_0=\stx^{k-1}$ being the approximate solution at the end of the stage  $k-1$. Method parameters are selected in such a way  that the upper bound $R_k^2$ for the expected squared distance $\bbe\|\stx^k - x^*\|^2$ between the approximate solution $\stx^k$ at the end of the $k$-th stage and the optimal solutions reduces by factor 2.}
\begin{corollary}\label{coro_strongly_convex}
	Let $\{\stx^K\}$ be {approximate solution} by Algorithm~\ref{adaptive_GEM_epoch} {after $K\geq 1$ stages}. Assume that $\| \stx^0-x^*\|^2\leq R_0^2$. Then 	\begin{align*}
		\bbe[f(\stx^K) - f^*] \leq  2^{-K - 1}\mu R_0^2 \quad \text{and} \quad \bbe[\|\stx^K - x^*\|^2]\leq  2^{-K}R_0^2.
	\end{align*}
\end{corollary}

\paragraph{Remarks.}
By Corollary~\ref{coro_strongly_convex}, the number of stages of Algorithm \ref{adaptive_GEM_epoch} required to attain the expected inaccuracy $\epsilon$ is bounded with $\mathcal{O}\left(\log(\mu R_0^2/{\epsilon})\right)$. When recalling what is the total number $N$ of iterations at each stage, we conclude that the ``total'' iteration complexity of the method is $\mathcal{O}\left( \sqrt{L\Omega/{\mu}}\,\log(\mu  R_0^2/{\epsilon}) \right)$. Consequently, the corresponding sample complexity  $\sum_{k=1}^{K} \sum_{t=1}^{N} m^k$ is of the order of
\begin{align*}
\tli{}{
\mathcal{O}\left\{\sqrt{\tfrac{L\Omega}{\mu}}\log\left(\tfrac{\mu  R_0^2}{\epsilon}\right)  + \tfrac{ \cL{\Omega} }{\mu}\log\left(\tfrac{\mu R_0^2}{\epsilon}\right) + \tfrac{\tli{\Omega^2}{\Omega}\sigma_*^2}{\mu \epsilon}\right\}.}
\end{align*}
Similarly, the iteration complexity of solution $\stx^K$ satisfying $\mathbb{E}[\|\stx^K-x^*\|^2] \leq \epsilon^2$ does not exceed \[\mathcal{O}\left(\sqrt{L\Omega/{\mu}}\,\log(  R_0/{\epsilon})\right);
\]
the corresponding sample complexity is then bounded by
\begin{align*}
\tli{}{
	\mathcal{O}\bigg\{ \sqrt{\tfrac{L \Omega}{\mu}}\log\left(\tfrac{ R_0}{\epsilon}\right)  + \tfrac{\cL \tli{\Omega^2}{\Omega} }{\mu}\log\left(\tfrac{  R_0}{\epsilon}\right) + \tfrac{\tli{\Omega^2}{\Omega}\sigma_*^2}{\mu^2 \epsilon^2}\bigg\}.}
\end{align*}
Similar to Corollary~\ref{coro_non_strongly_convex}, the three terms in the sample complexity bounds represent the deterministic error, the state-dependent stochastic error, and the state-independent stochastic error, respectively. Regardless of the dependence on $\Omega$, the multi-stage SGE achieves the optimal iteration complexity and sample complexity simultaneously, \textcolor{black}{supported by the lower bound in Theorem 5 of \cite{woodworth2021even}.}

% How do we know this? Both the first and the last term are fine. But, what about the middle term?

\section{SGE for sparse recovery}\label{sec_sparse}
An interesting application of stochastic optimization with state-dependent noise arises in the relation to the problem of sparse recovery when it is assumed that \eqref{cp-ASGD} has a sparse or low-rank solution $x^*$. This problem is motivated by applications in high-dimensional statistics, where stochastic estimation under sparsity or low-rank constraints has garnered significant attention.
To solve this problem, one usually builds a sample average approximation (SAA) of the expected risk and solves the resulting minimization problem. To enhance sparsity, a classical approach consists of incorporating a regularization term, e.g., $\ell_1$- or trace-norm penalty (as in Lasso), and minimizing the norm of the solution under constraint (as in Dantzig Selector), see, e.g.,  \cite{bickel2009simultaneous, candes2007dantzig, candes2011probabilistic,candes2011tight, candes2012exact,candes2006stable,candes2009near,fazel2008compressed,juditsky2011accuracy,tsybakov2011nuclear} among many others.

Stochastic approximation also serves as a standard approach to deal with the sparse recovery problem. However, utilizing traditional Euclidean stochastic approximation usually leads to sub-optimal complexity bounds:  in this setting, the expected squared $\ell_2$-error of the stochastic operator $\cG(x, \xi)$ is usually proportional to the problem dimension $n$. Therefore, non-Euclidean stochastic mirror descent (SMD) methods have been applied to address this issue. In particular, the SMD algorithm
% \revision{in \cite{shalev2009stochastic, srebro2010smoothness} attains}{described in \cite{shalev2009stochastic, srebro2010smoothness}
% computes an estimate $\widehat{x}$ after $N$ iterations, aimed at finding an $s$-sparse vector $x^*$, the proposed method achieves}
in \cite{shalev2009stochastic, srebro2010smoothness} attains
the high probability complexity bound $f(\widehat x) - f^* =\O(\sigma\sqrt{s / N})$ (up to some ``logarithmic factors'') \revision{under the sub-Gaussian noise assumption}{in recovering an $s$-sparse signal $x^*$ under the sub-Gaussian noise assumption with {subgaussianity parameter} $\sigma^2$}, referred to as ``slow rates'' for sparse recovery. To improve the convergence rate to $\O(\sigma^2 s/N)$, multi-stage routines exploiting the properties similar to strong/uniform convexity could be used, cf. \cite{ghadimi2012optimal,juditsky2011first,juditsky2014deterministic}. In \cite{agarwal2012stochastic, gaillard2017sparse}, the authors utilize the ``restricted'' strong convexity condition to establish \revision{$\O(\tfrac{\sigma^2 s}{N})$}{$\O(\sigma^2 s/N)$} complexity bounds when assuming that $N\gg s^2$. The latter assumption means that the optimal rates are only valid in the range $s\ll \sqrt{N}$ of sparsity parameter.  Recently, in \cite{juditsky2023sparse,ilandarideva2023stochastic}, new multi-stage stochastic mirror descent algorithms were proposed which rely on the idea of variance reduction. The proposed method improved the required number of iterations in each stage to $\O(s)$, thus attaining the best-known state-dependent stochastic error. However, the iteration complexity of those routines is still sub-optimal in the mini-batch setting, leaving room for further acceleration.

In this section, we suppose that stochastic optimization problem \rf{cp-ASGD} admits a sparse solution $x^*\in \feaReg$.
{Standard examples of the sparsity assumptions are as follows:}
% \footnote{For the sake of simplicity, we consider here two standard sparsity structures on $\bbr^n$. The case of another ``classical'' matrix sparsity---low rank structure may be addressed in the same way, see \cite{juditsky2014unified,juditsky2023sparse}.}
\begin{itemize}
    \item ``Vanilla'' sparsity: we assume that an optimal solution $x^*\in \feaReg$ has at most $s\ll n$ nonvanishing entries. We put $\|\cdot\|=\|\cdot\|_1$ and $\|\cdot\|_*=\|\cdot\|_\infty$.
    \item Group sparsity: let us partition the set $[n]$ into $b$ subsets $\{I_1, ..., I_B\}$, and let $x_b$ the $b$th block of $x$, meaning that $[x_b]_i=0$ for all $i\notin I_b$. We assume that the optimal $x^*\in \feaReg $ is a block vector with at most $s\leq B$ nonvanishing blocks $x_b$. We define
        $\|x\| = \tsum_{b=1}^B \|x_b\|_2$ (block $\ell_1/\ell_2$-norm) and $\|x\|_* = \max_{b\leq B} \|x_b\|_2$ (block $\ell_\infty/\ell_2$-norm).
    \item Low rank sparsity: consider the matrix space $\bbr^{p\times q}$ where $p\geq q$ \revision{eqipped}{equipped} with Frobenius inner product. We assume that the optimal $x^*\in \feaReg$ satisfies $\text{rank}(x^*)\leq s$. We consider the nuclear norm  $\|x\| = \tsum_{i=1}^q\sigma_i(x)$ where $\sigma_i(\cdot)$ are the singular values of $x$, so that $\|x\|_* = \max_{i\in [q]}\sigma_i(x)$ is the spectral norm.
\end{itemize}
In what follows we assume that problem \rf{cp-ASGD} and norms $\|\cdot\|$ and $\|\cdot\|_*$ satisfy conditions \rf{smoothness_ASGD}, \rf{assump:bias_SO}, and
state-dependent noise conditions in Assumptions  \eqref{assump:variance} and \eqref{assump:variance_1} of Section \ref{sec_assump}.
Also, instead of  $\mu$-quadratic growth (with $\mu>0$) condition with respect to $\|\cdot\|$ as in the previous section, we suppose $f$ and $x^*$ verify the quadratic growth condition with respect to the Euclidean norm, i.e.,
\begin{align}\label{eq:quadractic_growth}
    f(x) - f^* \geq \tfrac{1}{2}\underline \kappa\|x-x^*\|_2^2, ~~\forall x \in \feaReg,
\end{align}
for some $\underline \kappa>0$.
\subsection{Application: sparse generalized linear regression}\label{sec:glr}
Let us check that problem assumptions of this section hold in the case of generalized linear regression problem (GLR), as described in the introduction. %We are now looking for estimating the \emph{sparse} signal $x^*\in\bbr^n$ with sparse level $s$ from i.i.d. observations $(\phi_t,\eta_t)$ as in \rf{prob_GLM}.
Let us assume that
\begin{itemize}
    \item regressors $\phi_t$ satisfy $\mathbb{E}[\phi_1 \phi_1^\top] = \Sigma \succeq \kappa I$ with $\kappa > 0$ and $\|\Sigma\|_\infty \revision{}{:=\max_{i\in[n], j\in [n]}\Sigma_{i,j}} \leq \nu$;
    \item noises $\zeta_t$ are zero mean with bounded variance, i.e., $\mathbb{E}[\zeta_1^2] \leq 1$ without loss of generality;
    \item activation function $u$ is strongly monotone and Lipschitz continuous, i.e.,
%     \begin{align}\label{mono_lipstichz_old}
%     \underline r (t - t') \leq u(t) - u(t') \leq \bar r (t-t').
% \end{align}
\revision{}{for some $\bar r \geq  \underline r \geq 0$,}
\begin{align}\label{mono_lipstichz}
\revision{}{\left(u(t) - u(t')\right)(t - t') \geq \underline r (t - t')^2, \quad \text{and} \quad |u(t) - u(t')| \leq \bar r |t- t'|, \quad \forall~ t, t' \in \bbr.}
\end{align}
\end{itemize}
As already explained in the introduction, estimation of $x^*\in\inter\feaReg$ may be addressed through solving the stochastic optimization problem
\begin{align}\label{opt_GLM}
    \min_{x \in \feaReg} \left\{f(x) := \bbe[v(\phi^\top x) - \phi^\top x \eta]\right\}
\end{align}
where $v'(t) = u(t)$. The gradient of the problem objective and its stochastic estimate are given by
\begin{align*}
    g(x) = \bbe[\phi \big(u(\phi^\top x) - \eta\big) ] \quad \text{and} \quad \cG(x, (\phi, \zeta)) := \phi\big(u(\phi^\top x) - \eta\big)=
    \phi\big(u(\phi^\top x) - u(\phi^\top x^*)\big)-\phi\zeta.
\end{align*}
It is easy to see that condition \eqref{assump:bias_SO} is verified in this case, and invoking $\bbe[\eta] = \bbe[u(\phi^\top x)]$,  we conclude that $g(x^*)= 0$. To check the quadratic growth condition~\eqref{eq:quadractic_growth}, we write
\begin{align}\label{bound_f_f_star}
    f(x) - f^* &= \int_{0}^{1} g(x^* + t(x-x^*))^\top (x-x^*) dt\nn\\
    & = \int_{0}^{1} \bbe\big\{ \phi \big[\revision{u(\phi^\top x^* + t(x-x^*))}{u\big(\phi^\top (x^* + t(x-x^*))\big)} - u\big(\phi^\top x^*\big)\big]\big\}^\top(x-x^*) dt\nn\\
    \text{[by \eqref{mono_lipstichz}]}\quad&\geq \int_0^1\underline r \bbe\{[\phi^\top (x-x^*)]^2 \} t dt = \tfrac{\underline r}{2} \|x-x^*\|_\Sigma^2 \geq \tfrac{\underline r \kappa}{2} \|x-x^*\|_2^2.
\end{align}
Therefore, condition~\eqref{eq:quadractic_growth} holds with $\underline \kappa = \underline r \, \kappa$ which is independent of problem dimension $n$.

Since we are interested in the high-dimensional setting, the desired recovery error should have, at most, logarithmic dependence in the problem dimension $n$. However, the $\ell_2$ variance of the stochastic first-order information $\bbe\|\mathcal{G} - g\|_2^2$ is proportional to the problem dimension $n$, making the standard Euclidean SA methods not applicable. To address this issue, we work in the non-Euclidean setting with $\|\cdot\| = \|\cdot\|_1$ and $\|\cdot\|_* = \|\cdot\|_\infty$. Next, let us examine the smoothness condition \eqref{smoothness_ASGD} and state-dependent variance condition \eqref{assump:variance}. Note that for all $x, x' \in \bbr^n$,
\begin{align}\label{derive_smooth}
    \|g(x) - g(x')\|_\infty &= \sup_{\|z\|_1\leq 1} \langle g(x) - g(x'), z \rangle = \sup_{\|z\|_1\leq 1} \bbe\{\phi^\top z [u(\phi^\top x)- u(\phi^\top x')]\} \nn\\
    & \overset{(i)}\leq  \sup_{\|z\|_1\leq 1}\bar r \bbe\{|\phi^\top z| |\phi^\top (x-x')| \} \overset{(ii)}\leq \bar r  \sup_{\|z\|_1\leq 1}\sqrt{\bbe\{(\phi^\top z)^2 \}}  \|x-x'\|_\Sigma\nn\\ &\leq \bar r\sqrt{\nu} \|x-x'\|_\Sigma,
\end{align}
where (i) is a consequence of~\eqref{mono_lipstichz} and (ii) follows from the Cauchy inequality. Consequently, we have
\begin{align*}
    \bbe\|\cG(x, (\phi, \eta)) - g (x)\|_\infty^2 & =  \bbe\|\phi\big(u(\phi^\top x) - u(\phi^\top x^*)\big) + \phi \zeta + [g(x^*) -  g (x)] \|_\infty^2 \\
   \text{[by \eqref{derive_smooth}]}  &\leq 3\bbe\|\phi\big(u(\phi^\top x) - u(\phi^\top x^*)\big)\|_\infty^2 + 3\bbe\{\|\phi\|_\infty^2\}\sigma^2+ 3 \bar r^2 \nu \|x - x^*\|_\Sigma^2\\
   \text{[by \eqref{mono_lipstichz}]}\quad & \leq 3\bar r^2\bbe\{\|\phi\|_\infty^2\big(\phi^\top x- \phi^\top x^*\big)^2\} + 3\bbe\{\|\phi\|_\infty^2\}\sigma^2+ 3 \bar r^2 \nu \|x - x^*\|_\Sigma^2.
\end{align*}
By \eqref{bound_f_f_star}, we conclude that the condition \eqref{assump:variance} holds whenever
\begin{align*}
    \bbe[\|\phi\|_\infty^2\big(\phi^\top x- \phi^\top x^*\big)^2] \lesssim \bbe[\big(\phi^T(x- x^*)\big)^2] \lesssim\|x-x^*\|_\Sigma^2;
\end{align*}
e.g., when the regressor $\phi$ is bounded or sub-Gaussian.  Finally, under similar assumptions on the regressors and sub-Gaussian assumption on the additive noise $\zeta$, condition \eqref{assump:variance_1} naturally follows.
%To summarize, the stochastic optimization problem arising in the generalized linear regression satisfies our assumptions.
%
%
%

\subsection{SGE-SR: stochastic gradient extrapolation for sparse recovery}
We extend SGE to solve the sparse recovery problem. Similarly to Algorithm~\ref{adaptive_GEM_epoch}, sparse recovery routine is organized in stages; each stage represents a run of SGE (Algorithm~\ref{adaptive_GEM}). The principal difference with Algorithm~\ref{adaptive_GEM_epoch}, apart from the different choice of algorithm parameters, is the sparsity enforcing step (see, e.g., \cite{blumensath2009iterative, jain2014iterative, liu2020between}) implemented at the end of each stage.

Observe that for $x\in \feaReg$ one can efficiently compute a sparse approximation of $x$, specifically, $x_s=\mathsf{sparse}(x)$, an optimal solution to
\beq\label{eq:sparsex}
    \min \|x-z\|_2 ~\text{over $s$-sparse $z\in \feaReg$}.
\eeq
For instance, in the ``vanilla sparsity'' case, when the set $\feaReg$ is positive monotone,\footnote{\revision{}{For $x\in \bbr^n$, let $|x|$ denote a vector in $\bbr^n_+$ whose entries are absolute values of the corresponding entries of $x$.} We say that $\feaReg$ is positive monotone if whenever $x\in \feaReg$ and $|y|\leq |x|$ (the inequality is understood coordinate-wise), one also has $y\in \feaReg$. A typical example of a monotone convex set $\feaReg$  is a ball of an absolute norm in $\bbr^n$.} $x_s$ is obtained by zeroing all but $s$ largest in amplitude entries of $x$.\footnote{In the block sparsity case, when $\feaReg$ is positive block-monotone, the corresponding ``sparsification'' amounts to zeroing out all but $s$ largest (in $\ell_2$-norm) blocks of $x$; when $\feaReg$ is a ball of a  Schatten norm in the space of $p\times q$ real matrices, low rank $x_s$ may be obtained from $x$ by trimming the singular values of $x$.}

\begin{algorithm}[H]
	\caption{Stochastic gradient extrapolation method for sparse recovery (SGE-SR)}
	\begin{algorithmic}
		\State {\bf Input:} initial point $\sstx^0 \in X$.
		\For {$k =1, 2, \ldots,K$}
		\State (a) Set $N=40\sqrt{  sL\Omega/{\underline \kappa}}$ and $R_k=2^{-k/2}R_0$. Run $N$ iterations of SGE (Algorithm~\ref{adaptive_GEM}) with $x_0 = z_0 = \sstx^{k-1}$ and
		\begin{align*}
		    \theta_t &= t, ~~\alpha_t = \tfrac{t-1}{t}, ~~{\beta_t=\tfrac{3}{t+2}},  ~~ \eta_t = \tfrac{\eta}{t},~~ t= 1,...,N\\
      \eta &= \max\left\{{24}L, \tfrac{{18}(N+2) \cL}{m^k}, \tfrac{\sigma_*}{ R_k}\sqrt{\tfrac{2 (N+1)^3}{{\Omega}m^k}}\right\},\\
		    m^k &= \max\left\{1,  \left\lceil \tfrac{{3}\cL (N+2)}{ L} \right\rceil ,
\left\lceil \tfrac{{8} N (N+2)^2 \sigma_*^2}{{9\Omega}L^2  R_k^2 }\right\rceil \right\},~~ t= 0,...,N
		\end{align*}
	\State (b) Set $\stx^k = x_{N}$, where $x_{N}$ is the solution obtained in Step (a).  Calculate $$\sstx^k = \mathsf{sparse}(\stx^k).$$
		\EndFor
	\end{algorithmic} \label{adaptive_GEM_epoch_SR}
\end{algorithm}

The following corollary characterizes the convergence rate of SGE-SR for solving the sparse recovery problem.
\begin{corollary}\label{coro_SR}
	Let $\{\stx^k, \sstx^k\}$ be computed by Algorithm~\ref{adaptive_GEM_epoch_SR}. Assume $\|\sstx^0 - x^*\|^2\leq R_0^2$. Then we have for $k\geq 1$
	\begin{align*}
		\bbe[f(\stx^k) - f^*] \leq \underline \kappa s^{-1}  2^{-k+4}R_0^2  \quad \text{and} \quad \bbe[\|\sstx^k-x^*\|^2]\leq2^{-k} R_0^2.
	\end{align*}
\end{corollary}
\paragraph{Remarks.}
From the result of Corollary \ref{coro_SR} we conclude that the SGE-SR algorithm finds an $s$-sparse $\sstx^k \in \feaReg$ such that $\bbe[\|\revision{y}{\sstx^k} - x^*\|^2]\leq \epsilon^2$ for any $\epsilon \in (0, R_0)$ in at most $k=\mathcal{O}\left(\log(R_0/\epsilon)\right)$ stages. The corresponding iteration complexity of the SGE-SR is  $\mathcal{O}\left( \sqrt{\tfrac{s L\Omega}{\underline \kappa}}\log\left(\tfrac{R_0}{\epsilon}\right)\right)$, and  the overall sample complexity is
\begin{align*}
\tli{}{
    \mathcal{O}\bigg\{ \sqrt{\tfrac{s L\Omega}{\underline \kappa}}\log\left(\tfrac{R_0}{\epsilon}\right)  + \tfrac{s\cL \tli{\Omega^2}{\Omega} }{\underline \kappa}\log\left(\tfrac{ R_0}{\epsilon}\right) + \tfrac{ \tli{\Omega^2}{\Omega} s^2   \sigma_*^2}{\underline \kappa^2 \epsilon^2}\bigg\}.}
\end{align*}
Similarly, the iteration complexity of the solution $\stx^k\in \feaReg$ (which is not $s$-sparse in general) such that $\bbe[f(\stx^k) - f^*]\leq \epsilon$ is $\mathcal{O}\left( \sqrt{\tfrac{s L\Omega}{\underline \kappa}}\log\left(\tfrac{\underline \kappa R_0^2}{s \epsilon}\right)\right)$, while the total sample complexity is
\begin{align*}
\tli{}{
    \mathcal{O}\bigg\{ \sqrt{\tfrac{s L\Omega}{\underline \kappa}}\log\left(\tfrac{\underline \kappa R_0^2}{s \epsilon}\right)  + \tfrac{s\cL \tli{\Omega^2}{\Omega} }{\underline \kappa}\log\left(\tfrac{ \underline \kappa R_0^2}{s\epsilon}\right) + \tfrac{\tli{\Omega^2}{\Omega} s   \sigma_*^2}{\underline \kappa \epsilon}\bigg\}.}
\end{align*}
The above results may be compared to the convergence guarantees obtained in \cite{juditsky2023sparse} in the similar setting of the sparse recovery problem. The iteration complexity of the SGE-SR algorithm attains the optimal dependence on the problem's condition number, $\sqrt{sL /\underline{\kappa}}$, which improves over the corresponding result in \cite{juditsky2023sparse} by a factor of $\mathcal{O}(\sqrt{sL/ \underline{\kappa}})$.
Furthermore, the proposed solution matches the best-known sample complexity bounds for the stochastic error in \cite{juditsky2023sparse}.

\section{Numerical experiments} \label{sec:numeric}
\begin{figure}[ht]
    \centering
    \includegraphics[width=0.4\textwidth]{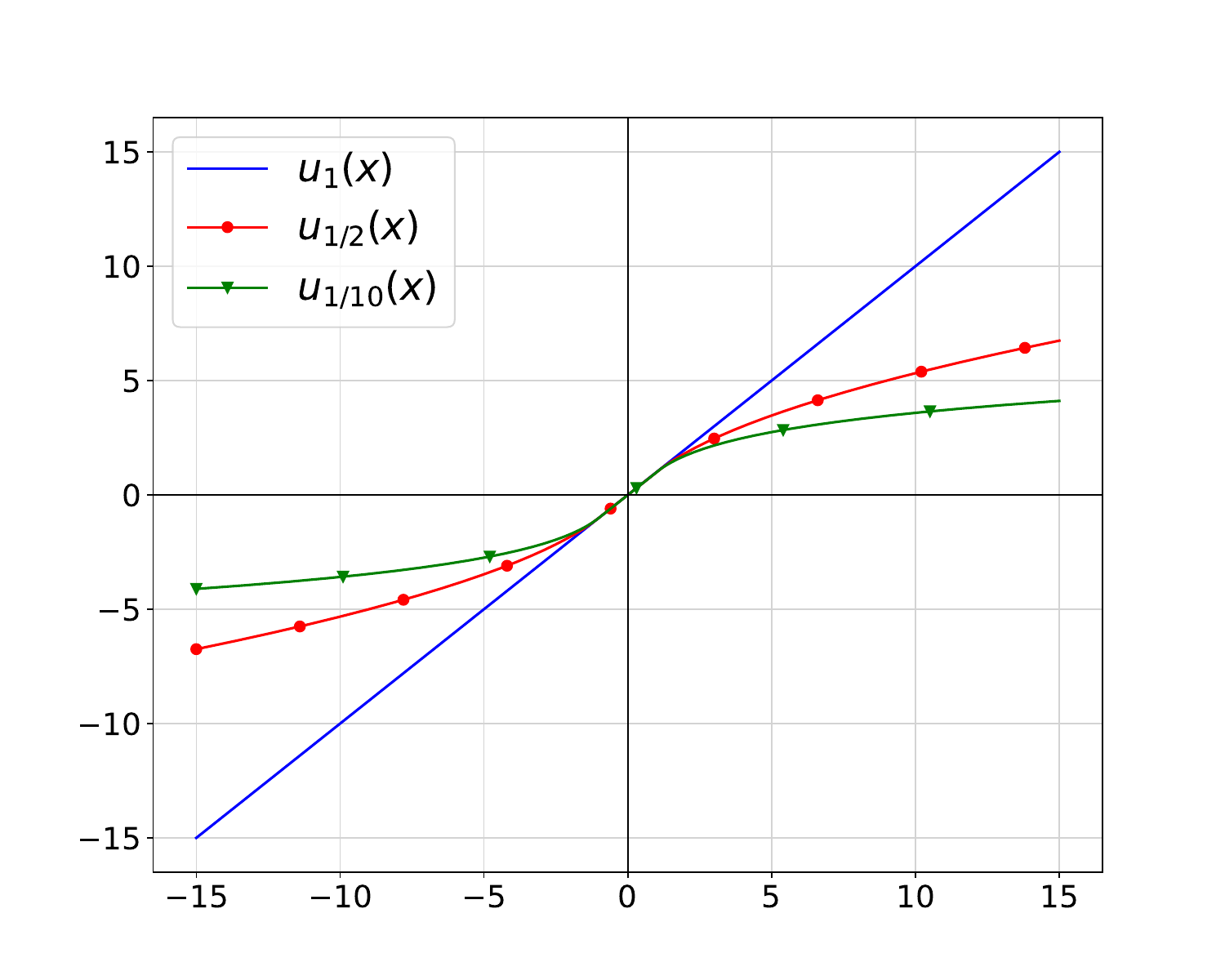}
    \caption{Activation functions}
    \label{fig:link functions}
\end{figure}

In this section we present a simulation study illustrating numerical performance of the proposed routines. We consider the sparse recovery problem in generalized linear regression (GLR) model with  random design as discussed in the previous section. Recall that we are looking to recover the $s$-sparse vector $x^* \in \bbr^n$ from i.i.d observations $$\eta_i = u(\phi_i^\top x^*) + \sigma \zeta_i, ~ i=1,2,\dots, N.$$
In  experiments we report on below, the activation function $u(\cdot)$ is of the form $$
u_{\alpha}(x) = x \mathbf{1}\{|x| \leq \alpha\} + \mathrm{sign}(x) [\alpha^{-1}(|x|^{\alpha} -1) + 1]\mathbf{1}\{|x| > \alpha\},
\quad \alpha >0, \;x \in \bbr.
$$
We consider three different activations, namely, the linear link function $u_{1}(\cdot)$, $u_{1/2}(\cdot)$, and $u_{1/10}(\cdot)$ (cf. Figure \ref{fig:link functions}).
%In our simulations, the $s$ nonvanishing components of the signal $x^*$ are sampled from the $s$-dimensional standard Gaussian distribution; regressors $\phi_i$ are independently drawn from a multivariate Gaussian distribution $\phi_i \sim \mathcal{N}(0, \Sigma)$, where $\Sigma$ is a diagonal covariance matrix with diagonal entries $0<\Sigma_{1,1} \leq \dots \leq \Sigma_{n,n}$. The condition number $\kappa$  of the problem is defined as the ratio between the largest and the smallest eigenvalues of $\Sigma$. The additive noise of the model is the zero-mean Gaussian noise with variance $\sigma^2$. Because of the memory limitations, observations $(\eta_i, \phi_i)$ are generated on the fly at each oracle calls.
In our simulations, $s$ nonvanishing components of the signal $x^*$ are sampled from the $s$-dimensional standard Gaussian distribution. We explore two setups---light tailed and heavy tailed---for generating regressors and additive noises. In the light-tail setup, regressors $\phi_i$ are independently drawn from a multivariate Gaussian distribution $\phi_i \sim \mathcal{N}(0, \Sigma)$, where $\Sigma$ is a diagonal covariance matrix with diagonal entries $0<\Sigma_{1,1} \leq \dots \leq \Sigma_{n,n}$. In the heavy-tail setup, regressors are independently drawn from a multivariate Student distribution $\phi_i \sim t_{n}(\nu, 0, \Sigma)$, $\nu$ being the corresponding degree of freedom \cite{kotz2004multivariate}. The condition number $\kappa$  of the problem is defined as the ratio of the largest and the smallest eigenvalues of $\Sigma$. The additive noise of the model in the light-tail setup is the zero-mean Gaussian noise with variance $\sigma^2$; in the heavy-tail setup, the additive noises have (scaled) univariate Student distribution  $\eta_i\sim \lambda t(\nu)$, $\nu\geq 3$, with scale parameter $ \lambda=\sqrt{(\nu - 2)/ \nu}$ with unit variance. Because of the memory limitations, observations $(\eta_i, \phi_i)$ are generated on the fly at each oracle call.

In all our experiments, we run 50 simulation trials (with randomly generated regressors and noises); then we trace in the plots the median and the first and the last deciles of the error $\|x_t-x^*\|_2$.

% The aim of the first series of experiments is to compare the procedure described in Section \ref{sec_sparse} to the  SMD-SR algorithm of \cite{juditsky2023sparse};\footnote{SMD-SR is a stochastic approximation algorithm for sparse recovery utilizing hard thresholding which relies upon ``vanilla'' non-Euclidean mirror descent; both algorithms use the same distance generating function $\omega(x) = c(n) \|x\|^2_p$.} the corresponding results are presented in Figure \ref{fig:SGE-SR vs SMD-SR}.
% In these experiments $n = 500\,000$, the maximal number of calls to the stochastic oracle (estimation sample size) $N =250\,000$, and sparsity level $s =250$; unless stated otherwise, problem condition number is set to $\kappa=1$.

\begin{figure}[ht]
  \centering
  \begin{subfigure}[b]{0.32\linewidth}
    \includegraphics[width=\linewidth]{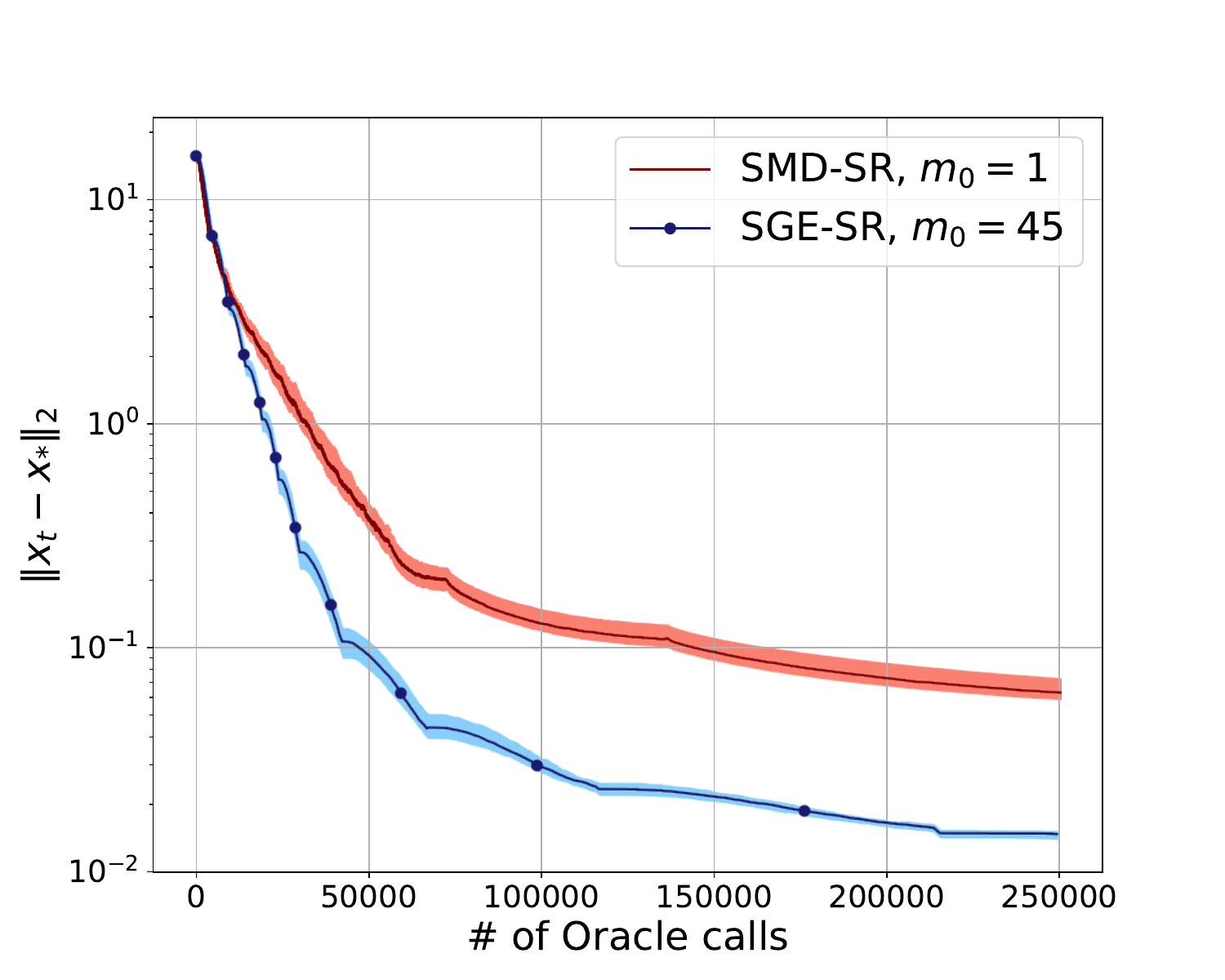}

  \end{subfigure}
  \begin{subfigure}[b]{0.32\linewidth}
    \includegraphics[width=\linewidth]{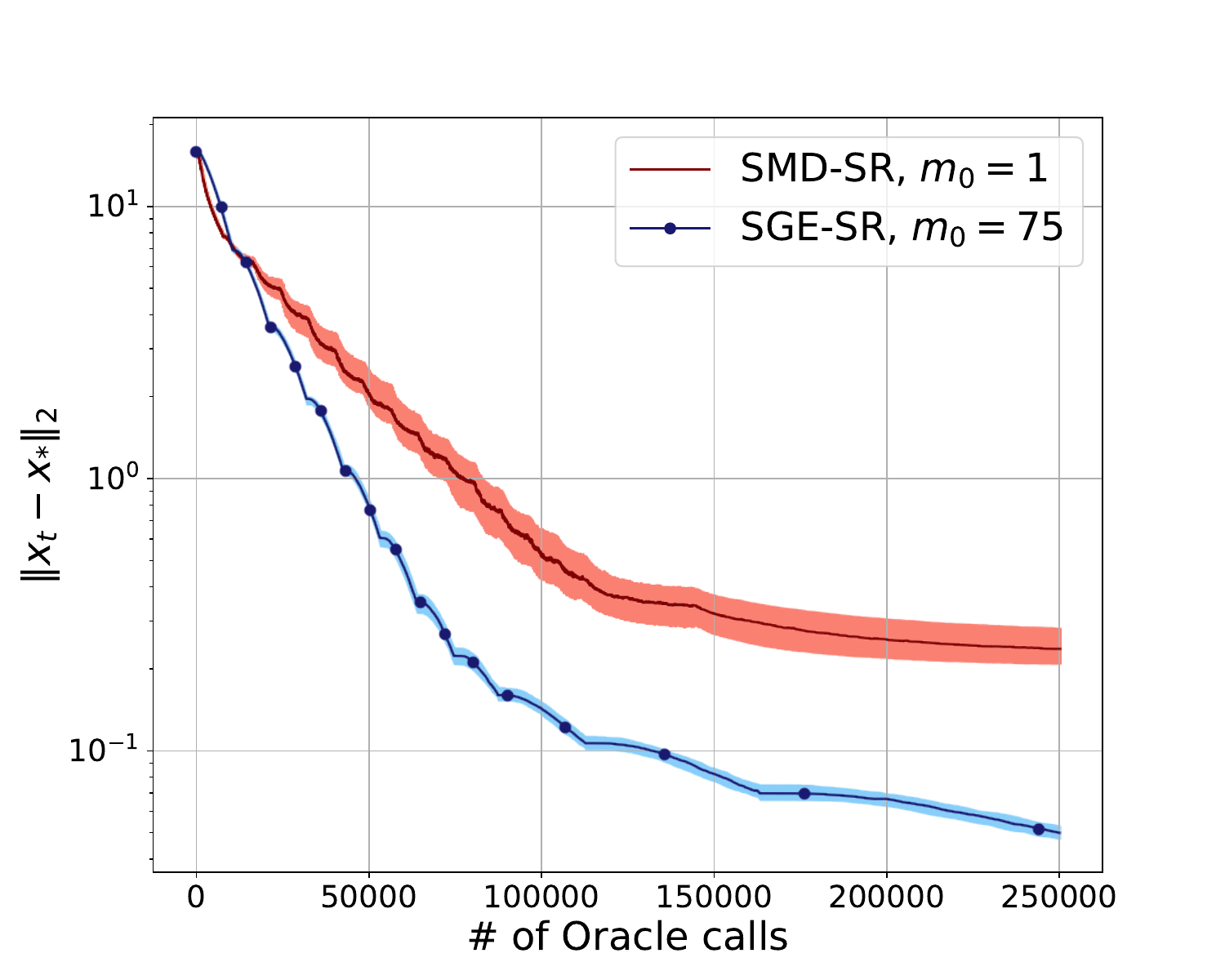}

  \end{subfigure}
  \begin{subfigure}[b]{0.32\linewidth}
    \includegraphics[width=\linewidth]{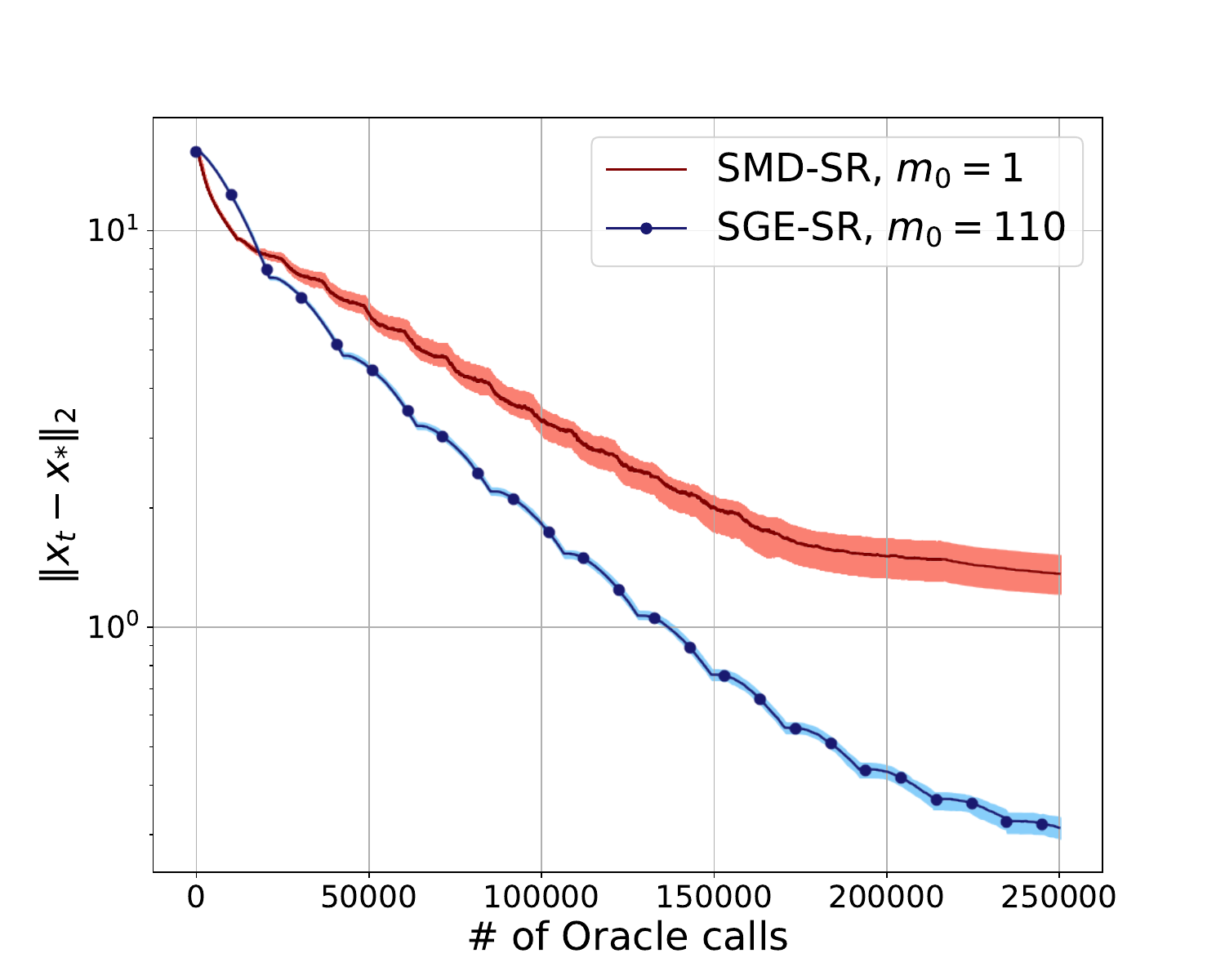}

  \end{subfigure}
  \medskip
  \centering
  \begin{subfigure}[b]{0.32\linewidth}
    \includegraphics[width=\linewidth]{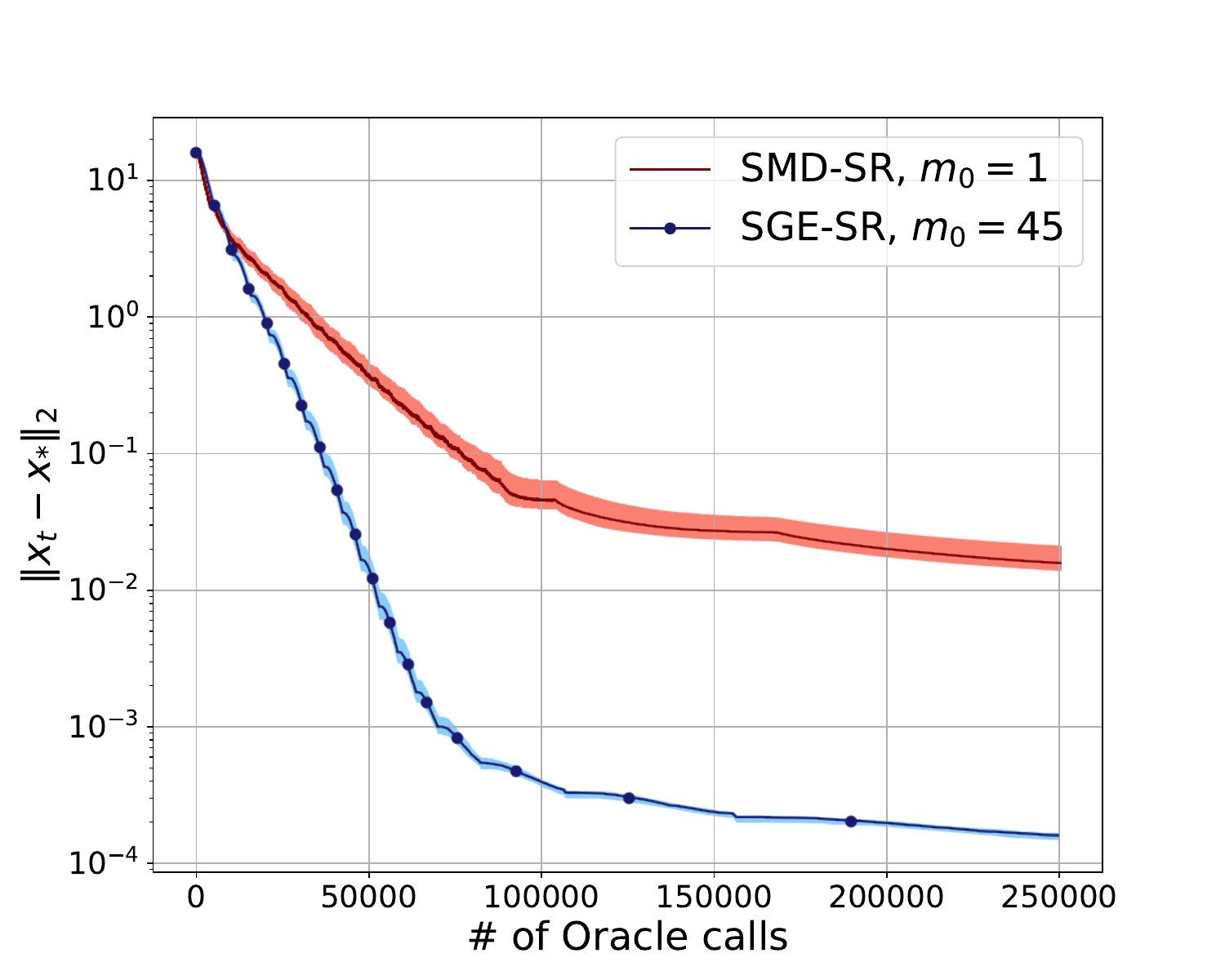}

  \end{subfigure}
  \begin{subfigure}[b]{0.32\linewidth}
    \includegraphics[width=\linewidth]{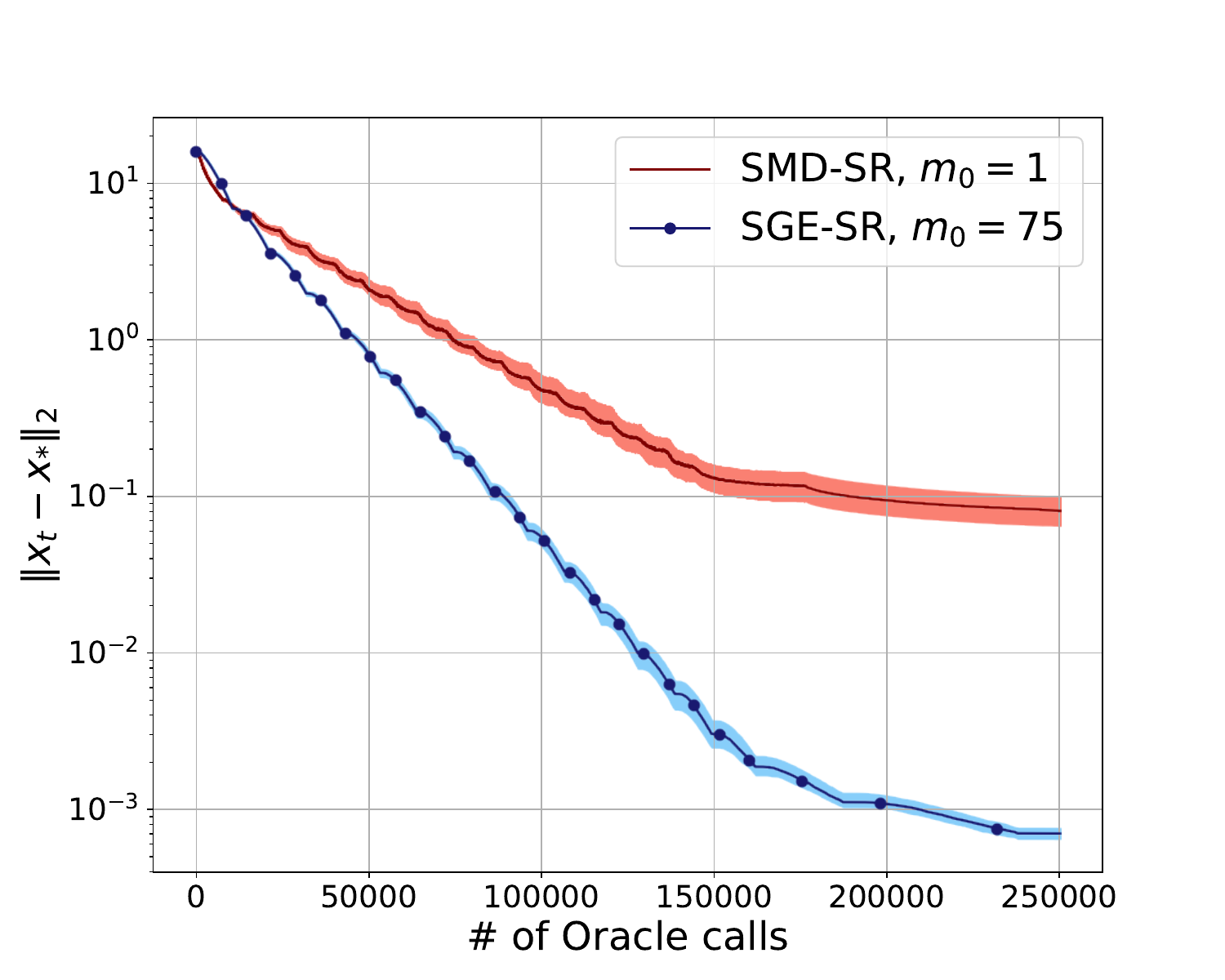}

  \end{subfigure}
  \begin{subfigure}[b]{0.32\linewidth}
    \includegraphics[width=\linewidth]{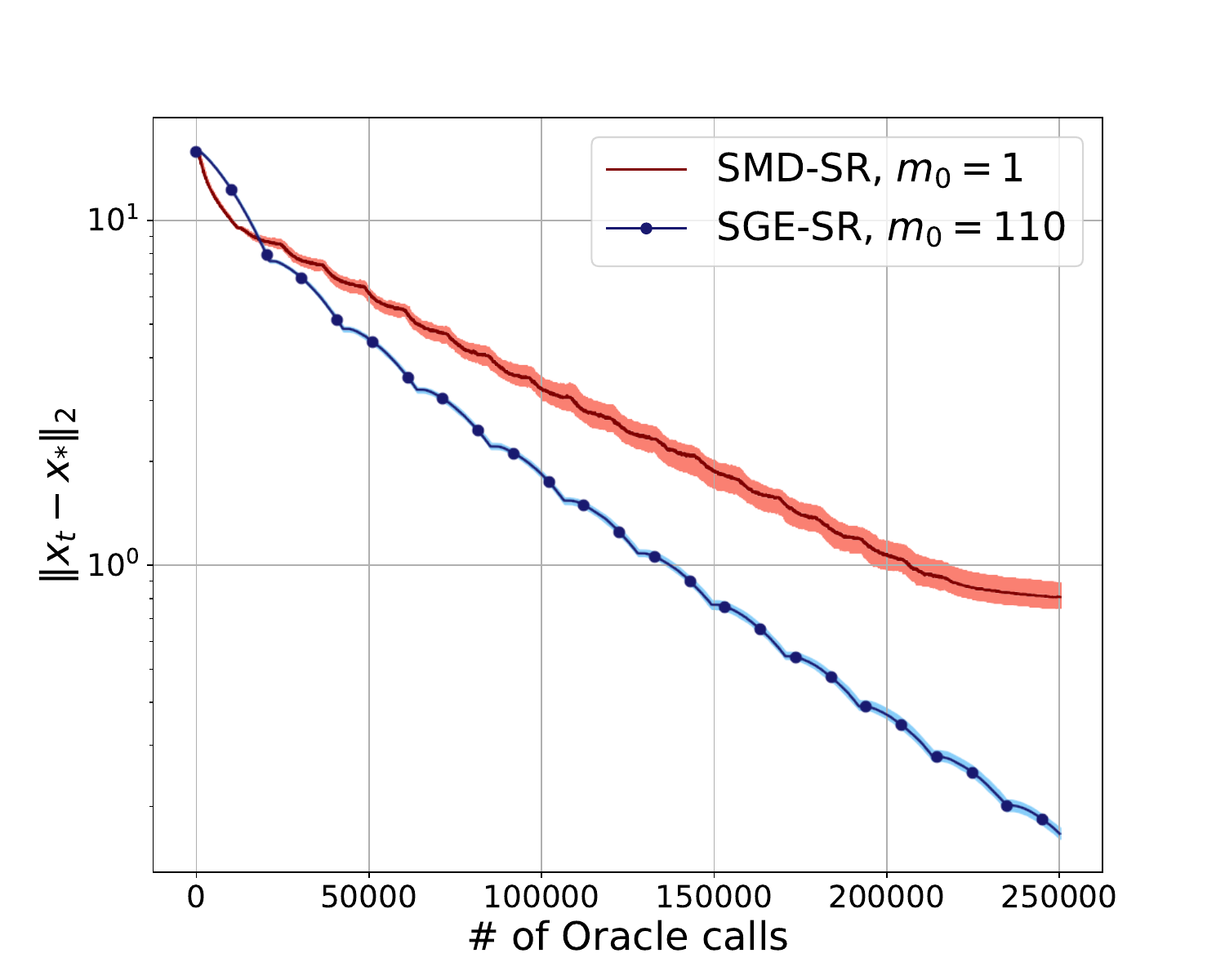}

  \end{subfigure}
  \caption{Estimation error $\|x_t -x^*\|_2$ against the number of stochastic oracle calls for SGE-SR and SMD-SR algorithms. In the left, middle, and right columns of the plot we show results for the linear activation $u_{1}$, and nonlinear $u_{1/2}$ and $u_{1/10}$, respectively. Two figure rows correspond to two different noise levels, $\sigma=0.1$ (the upper row) and $\sigma=0.001$ (the bottom row). The legend specifies the value $m_0$ of the batch size of the preliminary phase of the algorithm for both routines.}
  \label{fig:SGE-SR vs SMD-SR}

\end{figure}

\begin{figure}[h!]
  \centering
  \begin{subfigure}[b]{0.32\linewidth}
    \includegraphics[width=\linewidth]{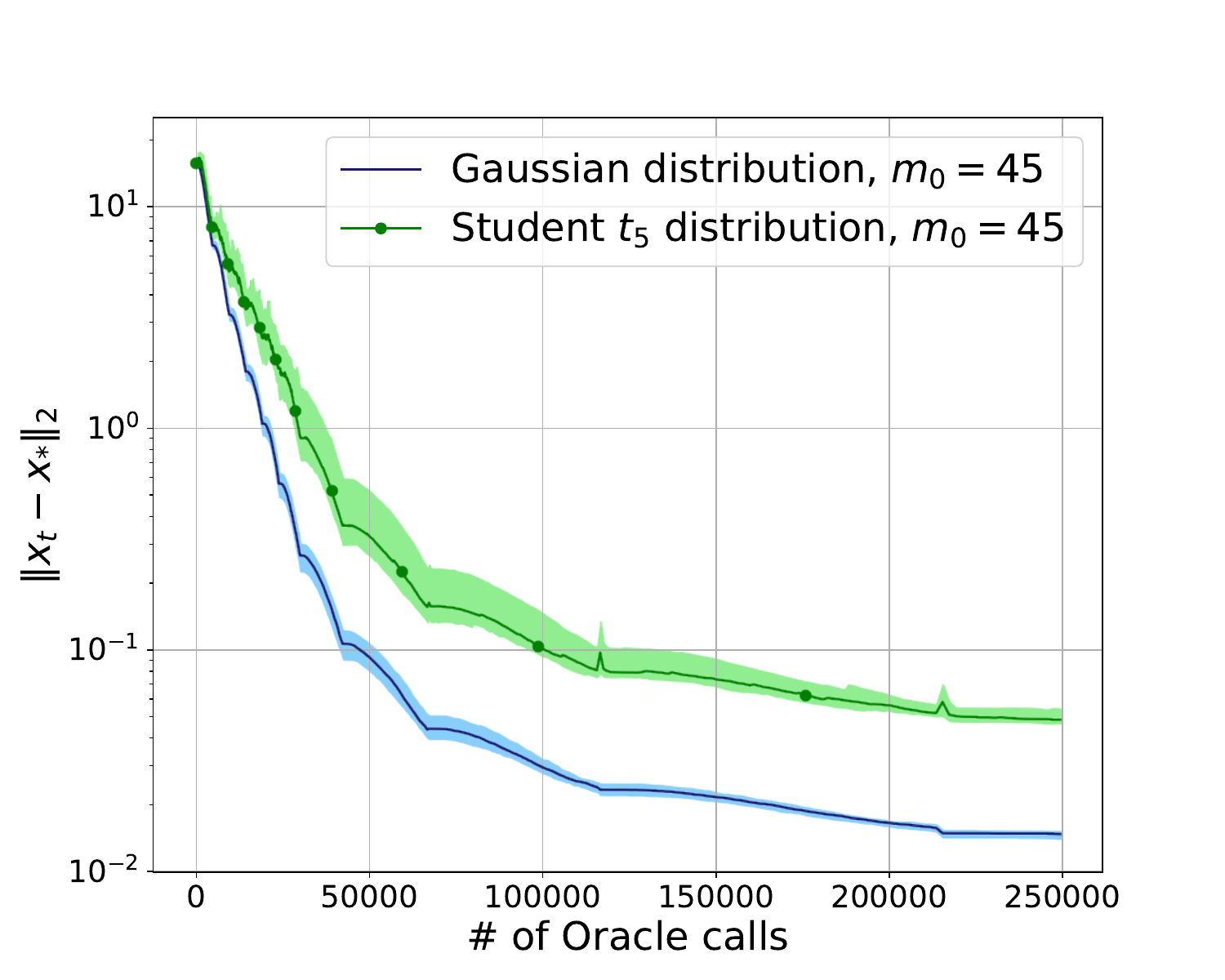}

  \end{subfigure}
  \begin{subfigure}[b]{0.32\linewidth}
    \includegraphics[width=\linewidth]{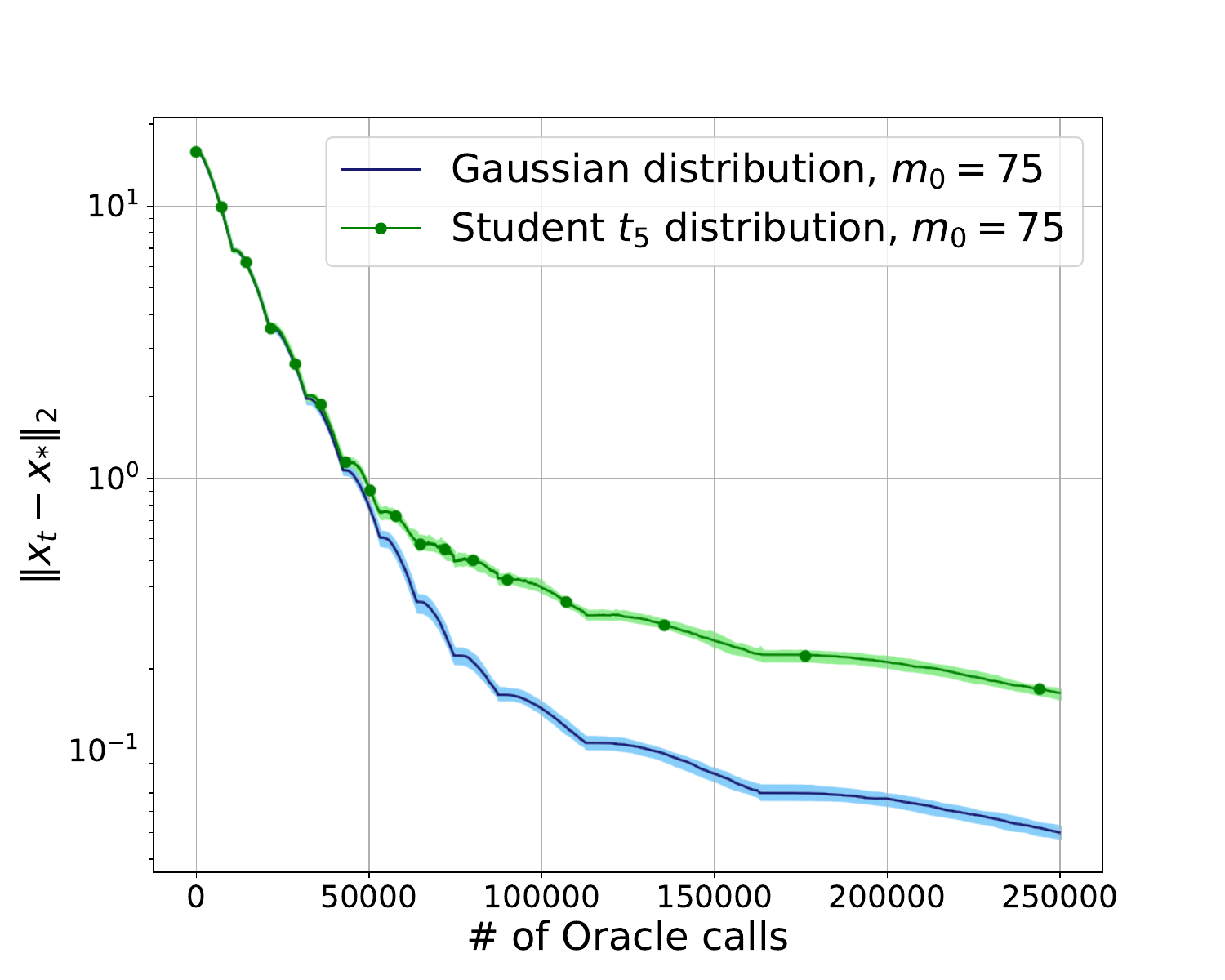}

  \end{subfigure}
  \begin{subfigure}[b]{0.32\linewidth}
    \includegraphics[width=\linewidth]{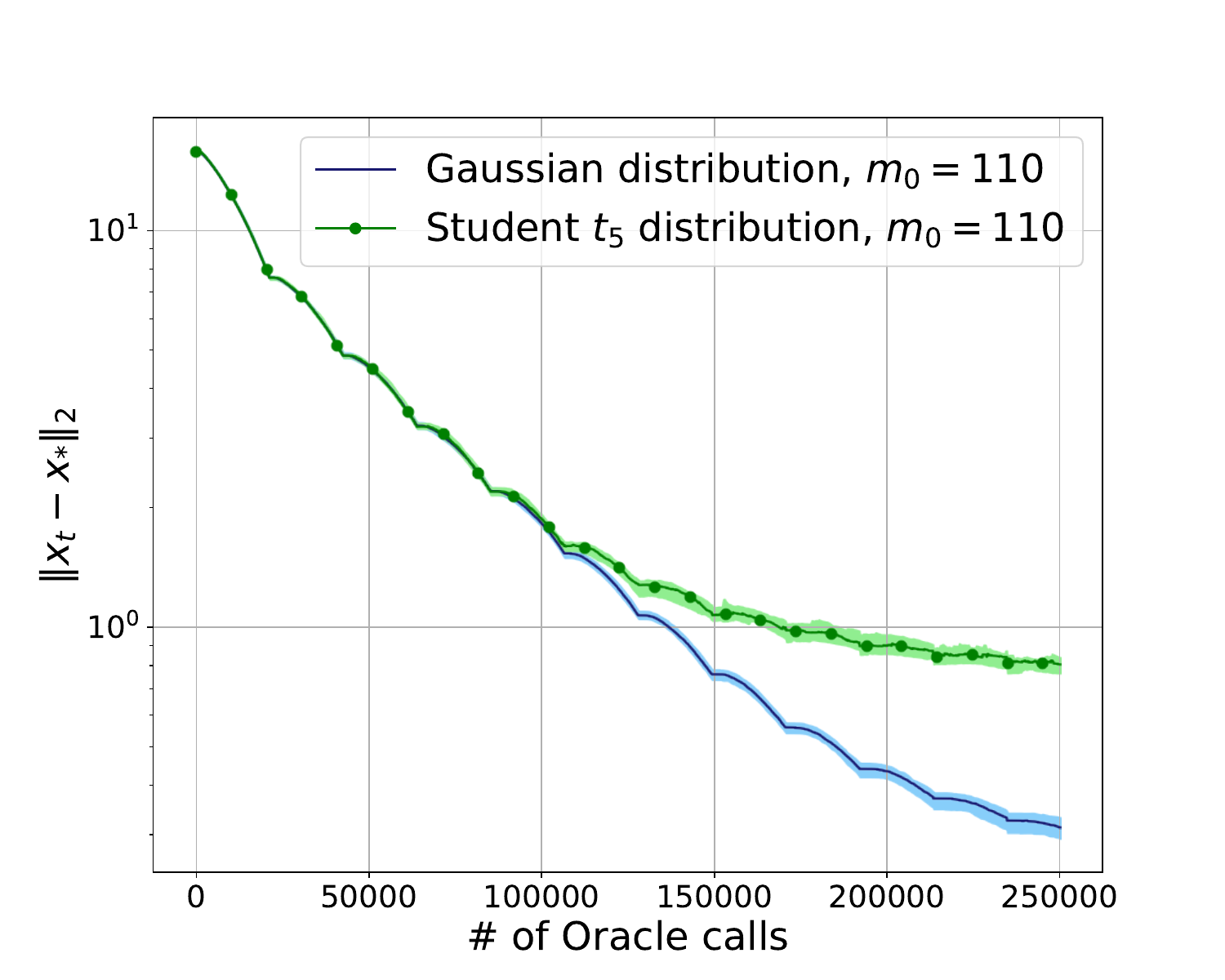}

  \end{subfigure}
  \medskip
  \centering
  \begin{subfigure}[b]{0.32\linewidth}
    \includegraphics[width=\linewidth]{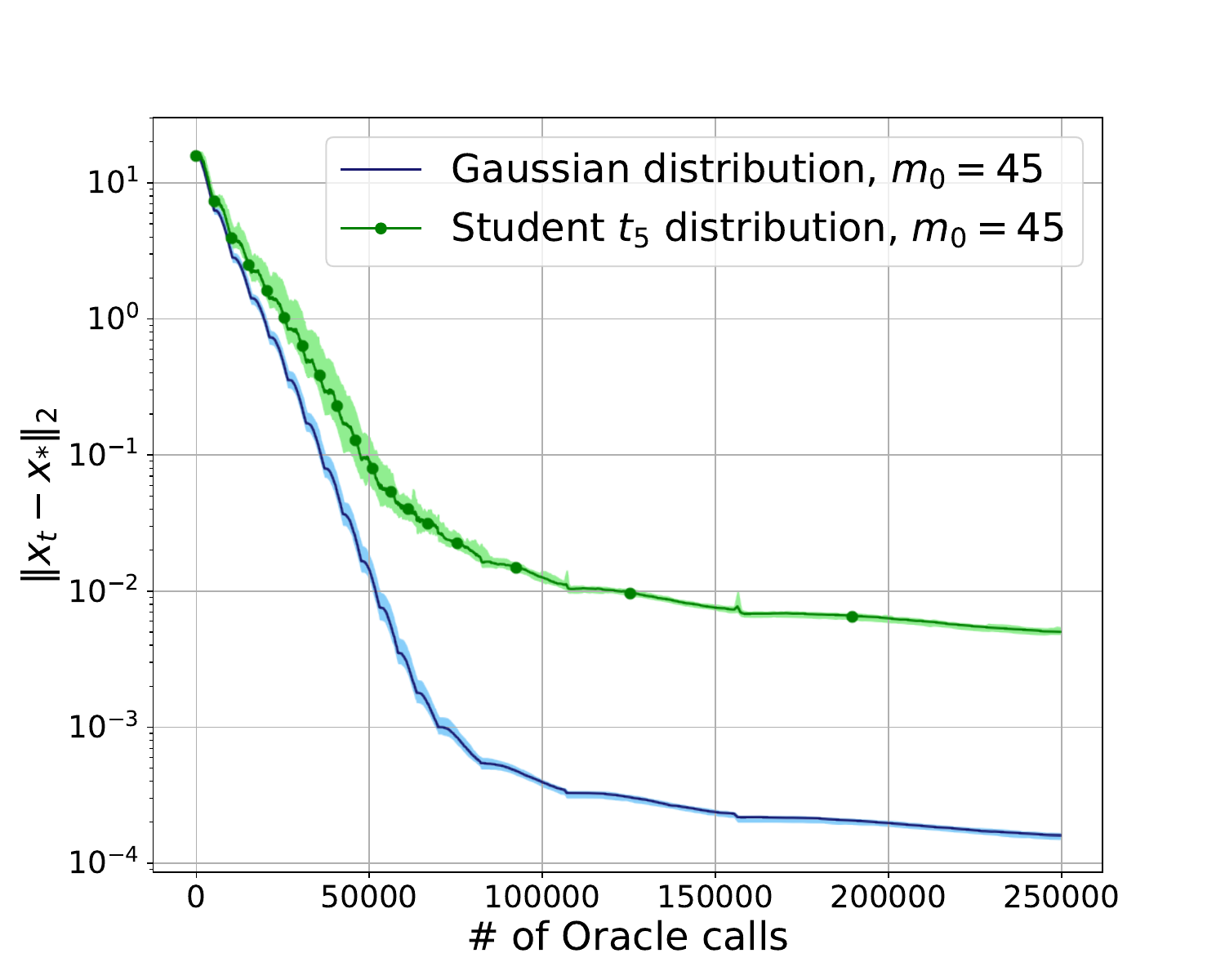}

  \end{subfigure}
  \begin{subfigure}[b]{0.32\linewidth}
    \includegraphics[width=\linewidth]{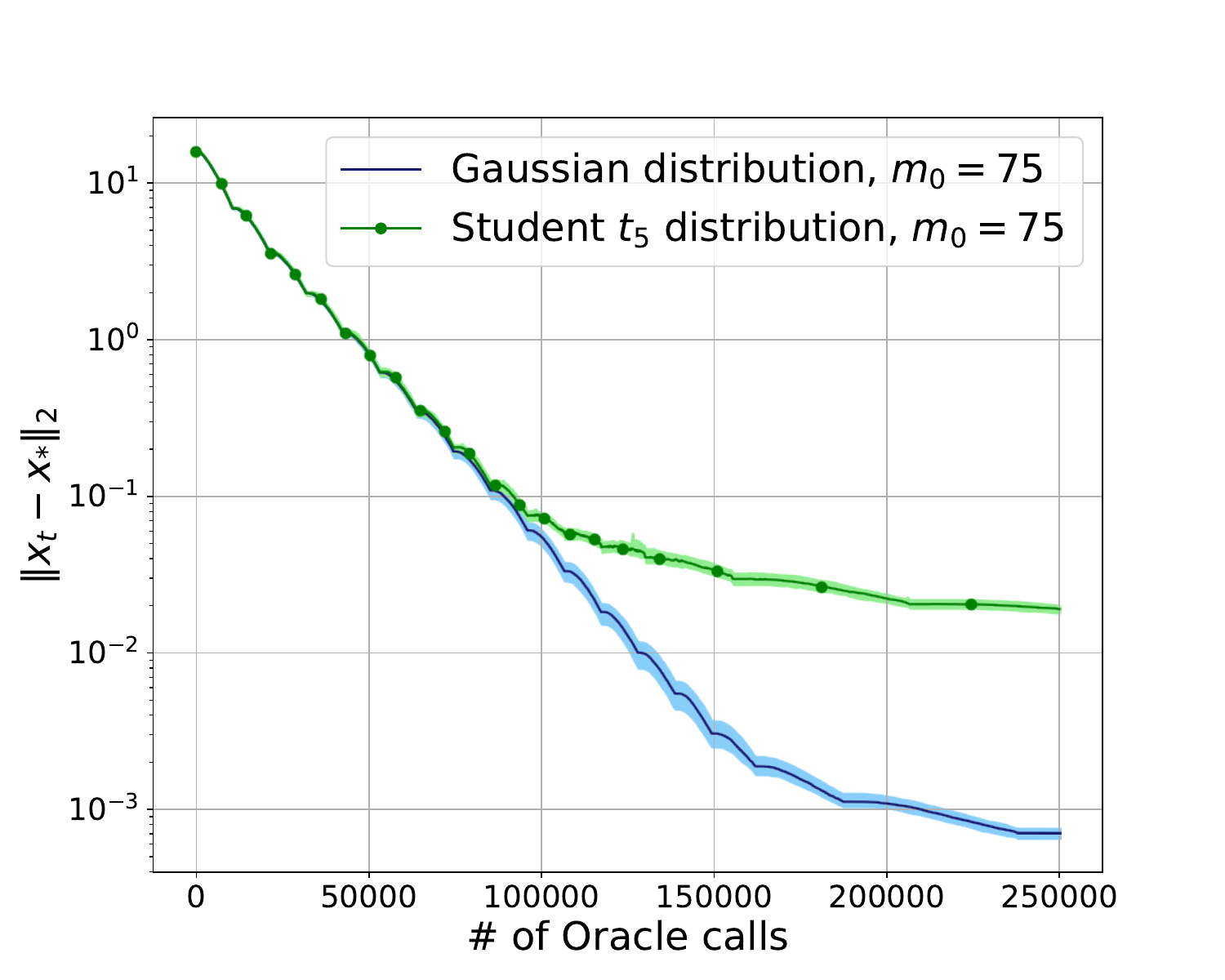}
  \end{subfigure}
  \begin{subfigure}[b]{0.32\linewidth}
    \includegraphics[width=\linewidth]{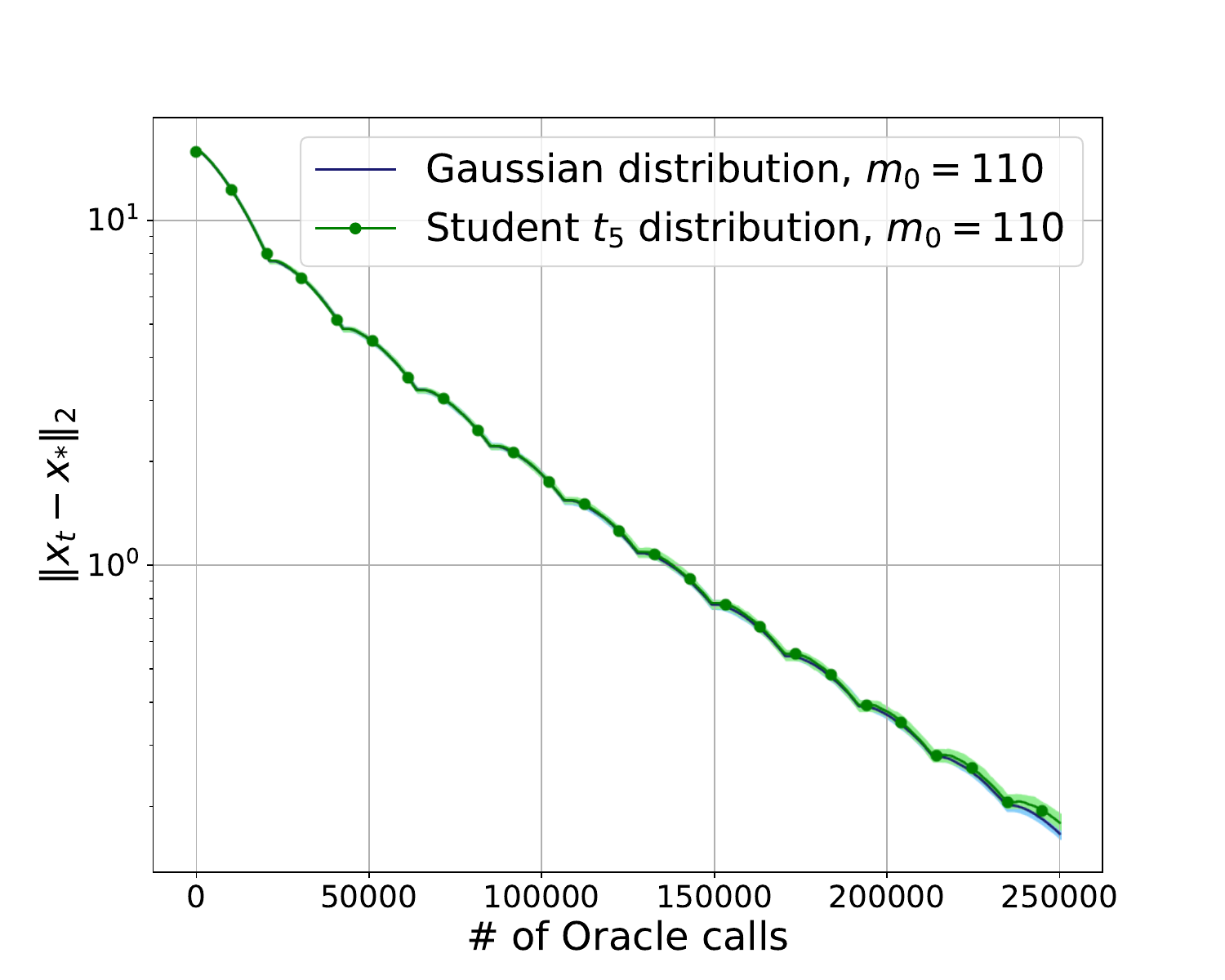}

  \end{subfigure}
  \caption{Estimation error $\|x_t -x^*\|_2$ against the number of stochastic oracle calls for SGE-SR in Gaussian (light-tail) and Student $t_5$ (heavy-tail) regressor and noise generation setups. In the left, middle, and right columns of the plot we show results for the linear activation $u_{1}$, and nonlinear $u_{1/2}$ and $u_{1/10}$, respectively. Two figure rows correspond to two different noise levels, $\sigma=0.1$ (the upper row) and $\sigma=0.001$ (the bottom row). The legend specifies the value $m_0$ of the batch size of the preliminary phase of the algorithm for both routines.}
  \label{fig:light-tail vs heavy-tail}

\end{figure}

The aim of the first series of experiments is to compare the procedure described in Section \ref{sec_sparse} to the  SMD-SR algorithm of \cite{juditsky2023sparse};\footnote{SMD-SR is a stochastic approximation algorithm for sparse recovery utilizing hard thresholding which relies upon ``vanilla'' non-Euclidean mirror descent; both algorithms use the same distance generating function $\omega(x) = c(n) \|x\|^2_p$.} in the light-tail noise setting, the corresponding results are presented in Figure~\ref{fig:SGE-SR vs SMD-SR}. In Figure \ref{fig:light-tail vs heavy-tail}, we present results of simulations of the accelerated algorithm in the light-tail and heavy-tail setup. We used the same algorithmic parameters in both simulation setups. In the above experiments, $n = 500\,000$, the maximal number of calls to the stochastic oracle (estimation sample size) $N =250\,000$, and sparsity level $s =250$; unless stated otherwise, the problem condition number is set to $\kappa=1$.

\begin{figure}[ht]
  \centering
  \begin{subfigure}[b]{0.32\linewidth}
    \includegraphics[width=\linewidth]{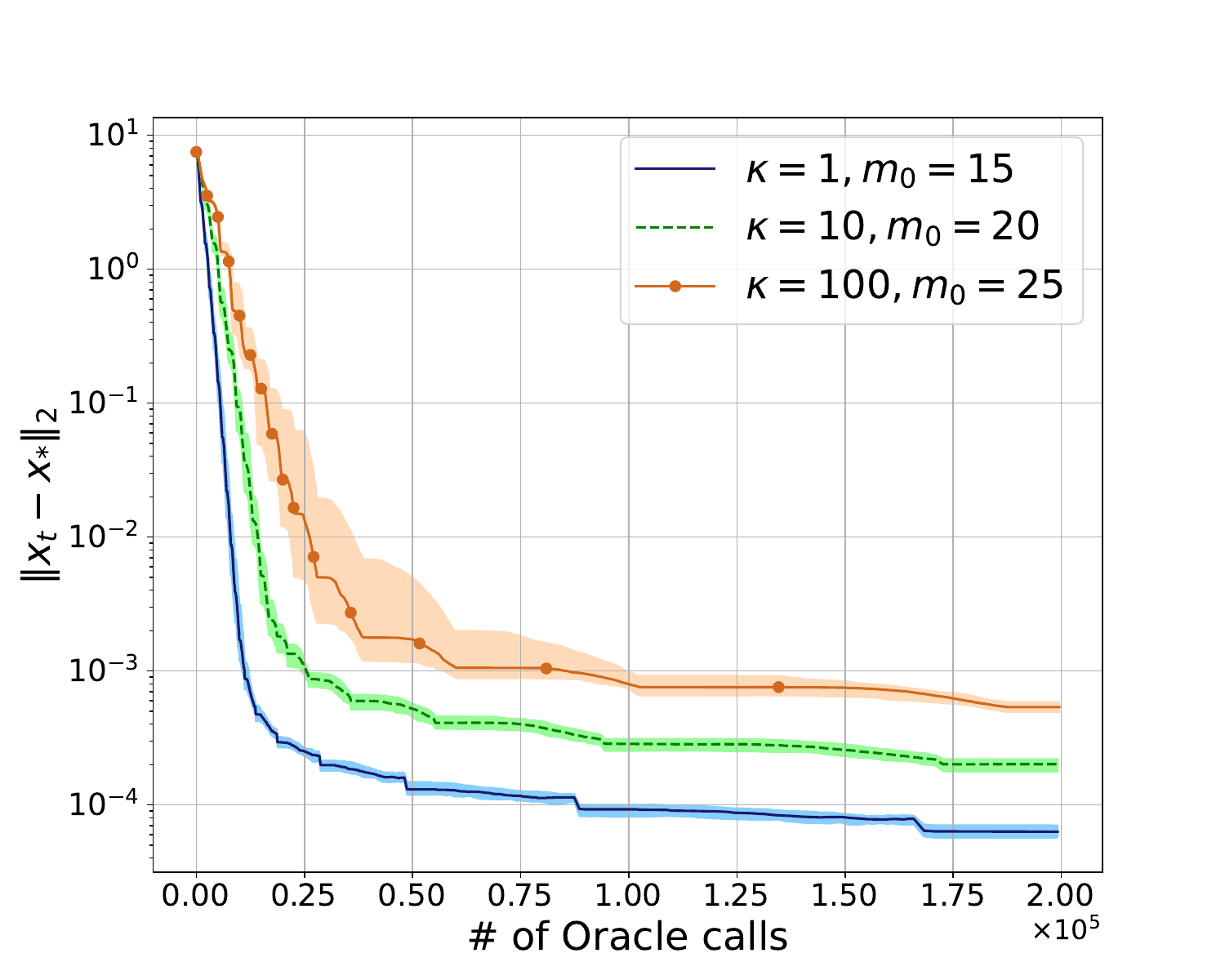}

  \end{subfigure}
  \begin{subfigure}[b]{0.32\linewidth}
    \includegraphics[width=\linewidth]{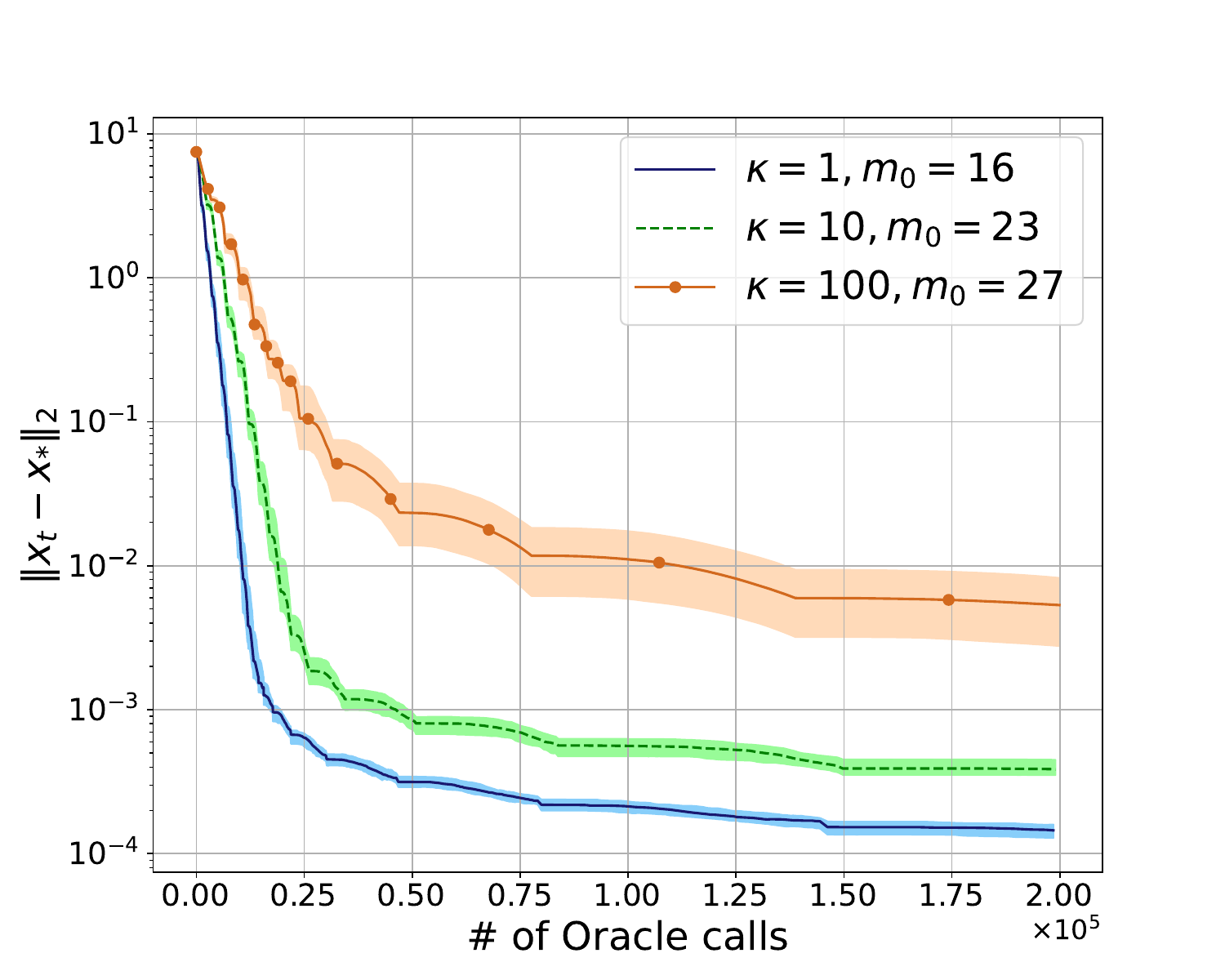}

  \end{subfigure}
  \begin{subfigure}[b]{0.32\linewidth}
    \includegraphics[width=\linewidth]{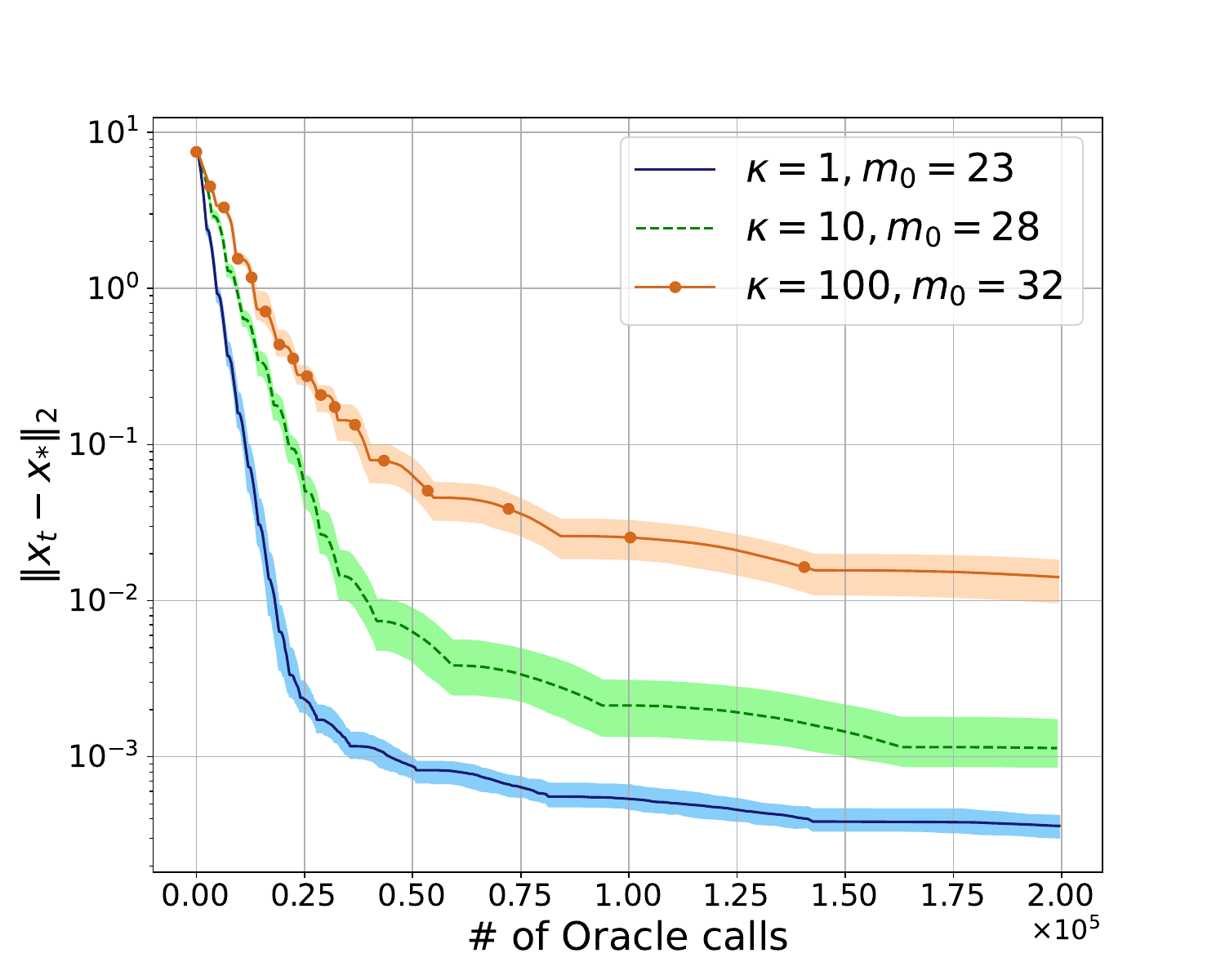}

  \end{subfigure}
  \medskip
  \centering
  \begin{subfigure}[b]{0.32\linewidth}
    \includegraphics[width=\linewidth]{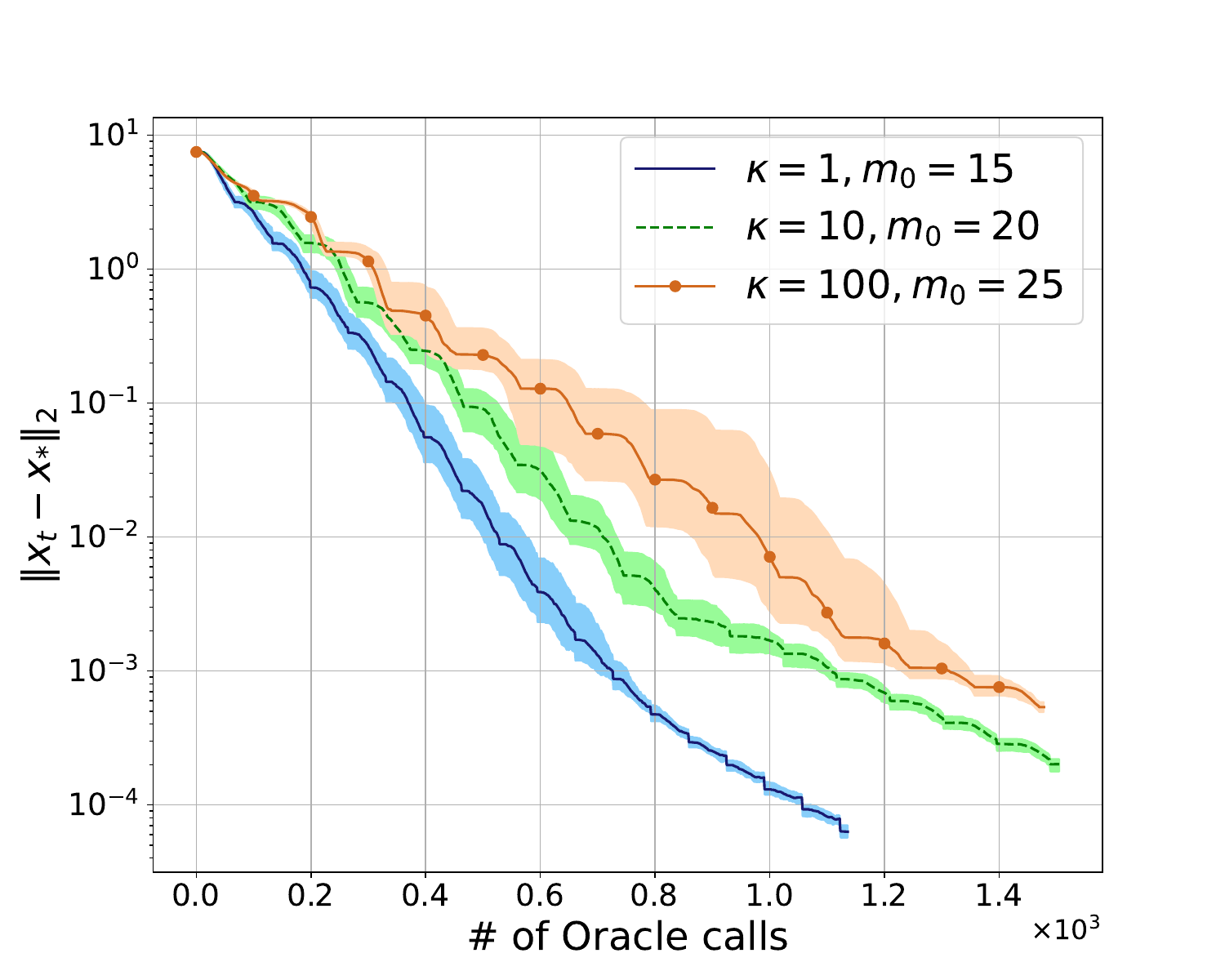}

  \end{subfigure}
  \begin{subfigure}[b]{0.32\linewidth}
    \includegraphics[width=\linewidth]{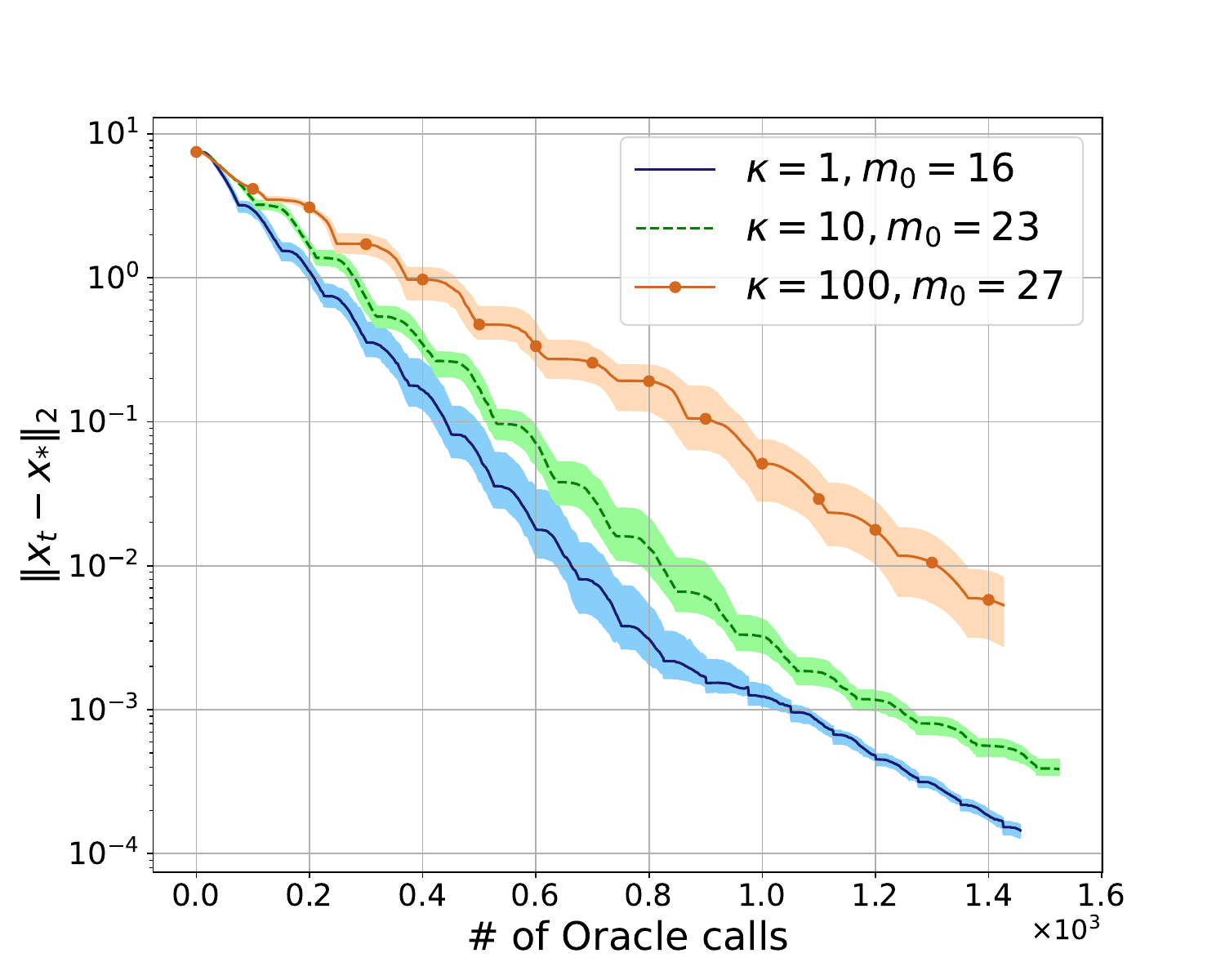}

  \end{subfigure}
  \begin{subfigure}[b]{0.32\linewidth}
    \includegraphics[width=\linewidth]{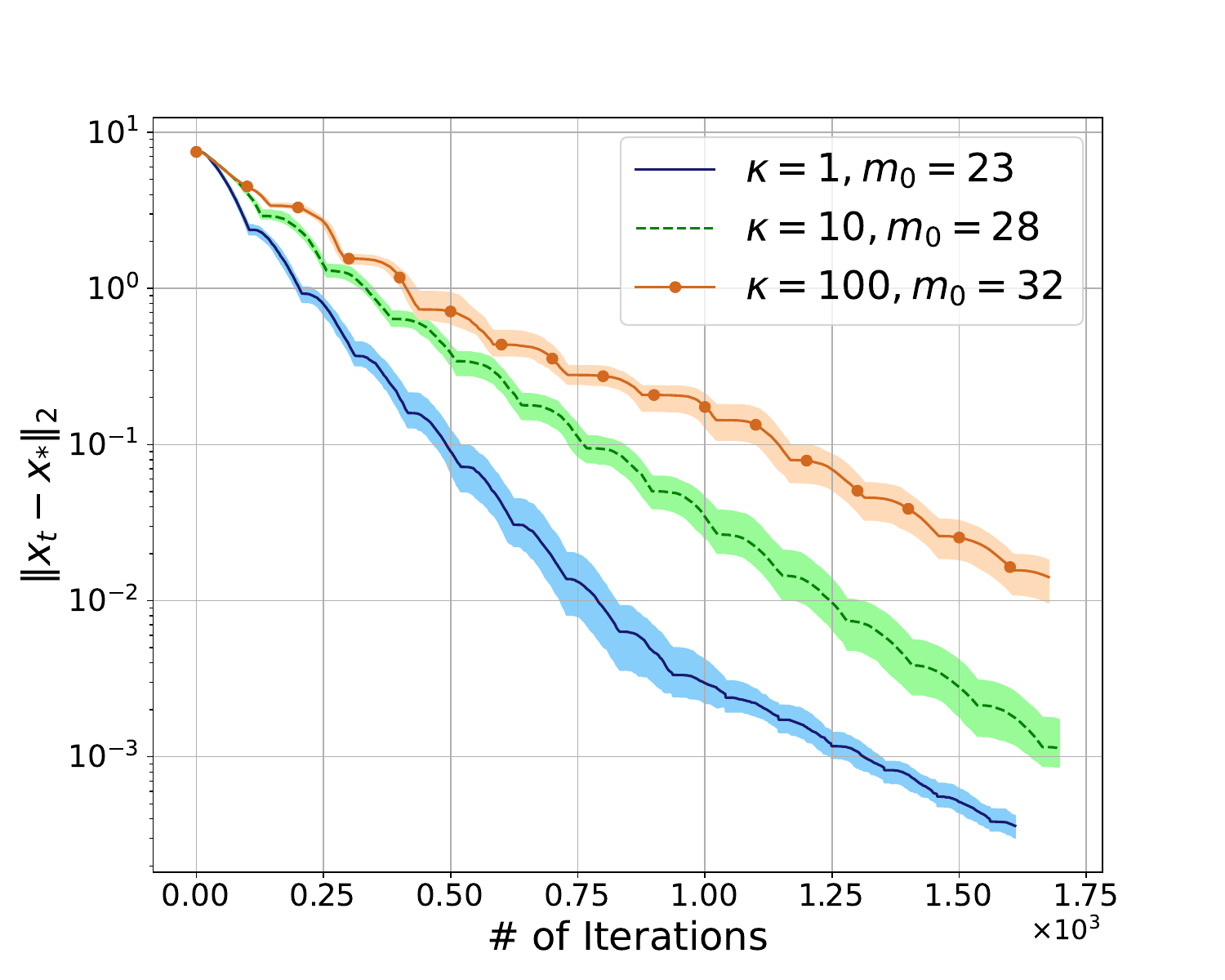}

  \end{subfigure}
  \caption{Error $\|x_t -x^*\|_2$ of the SGE-SR algorithm for three values of the problem condition number. First row:  algorithm error against the number of oracle calls; second row:  error against the number of algorithm iterations.   Figure columns correspond to the results for  $u_1$, $u_{1/2}$, and $u_{1/10}$ activation functions and $\sigma=0.001$.}
  \label{fig:SGE-SR different condition number}

\end{figure}

In the second series of experiments we put $n = 100\,000$, $N =200\,000$, and $s = 50$. Experiments reported in Figure \ref{fig:SGE-SR different condition number} illustrate the impact of the condition number on the convergence of the SGE-SR algorithm.

\begin{figure}[h!]
  \centering
  \begin{subfigure}[b]{0.32\linewidth}
    \includegraphics[width=\linewidth]{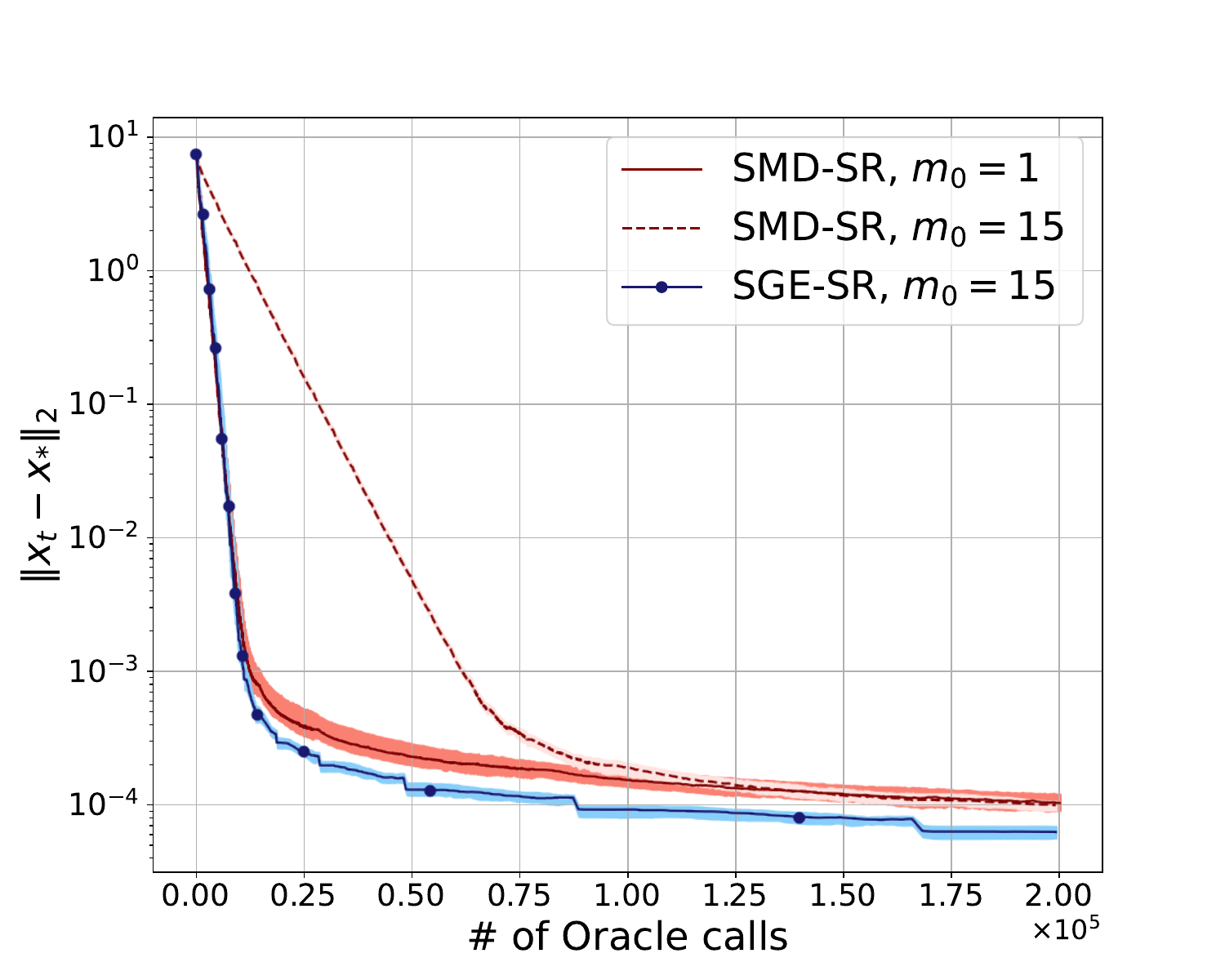}

  \end{subfigure}
  \begin{subfigure}[b]{0.32\linewidth}
    \includegraphics[width=\linewidth]{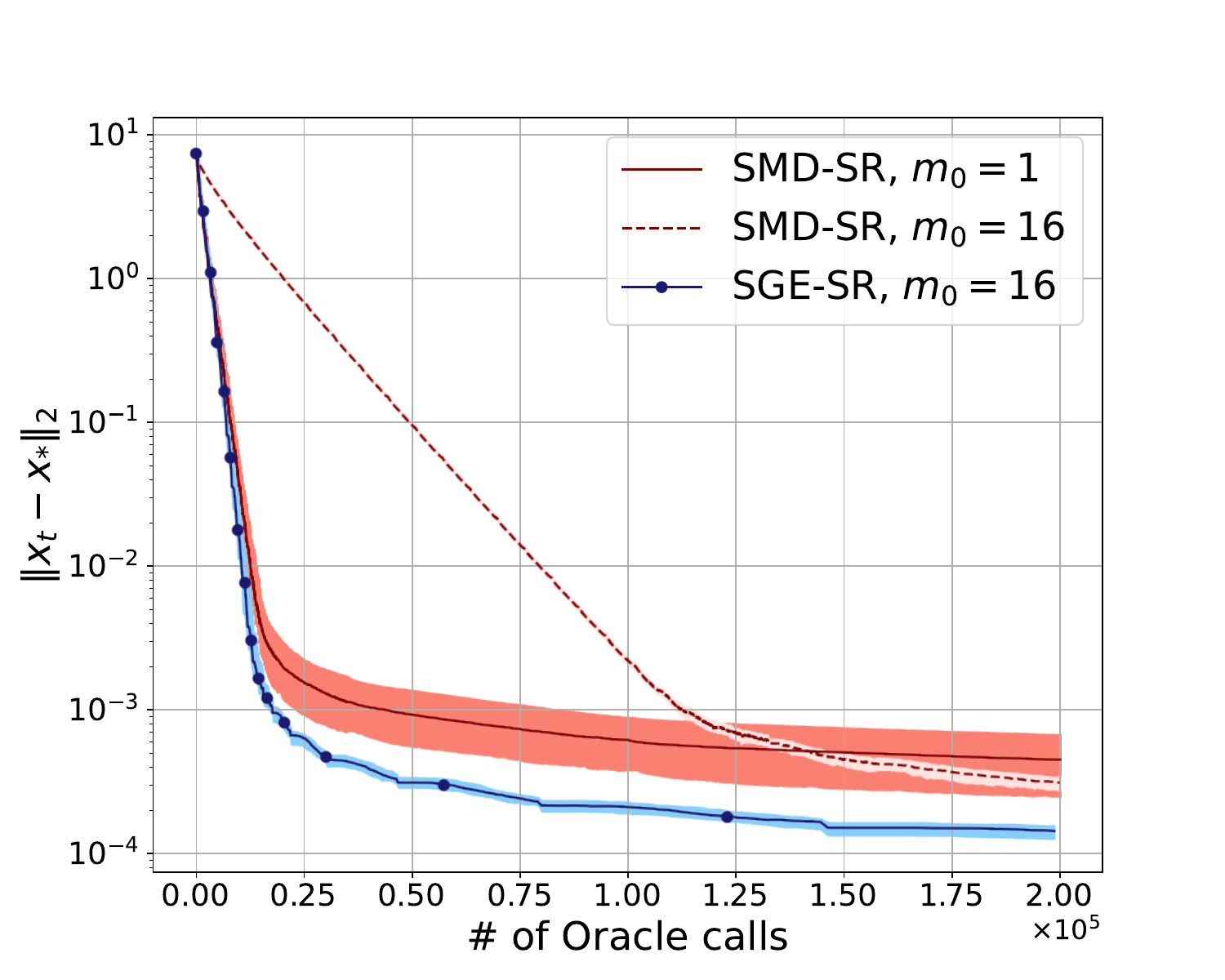}

  \end{subfigure}
  \begin{subfigure}[b]{0.32\linewidth}
    \includegraphics[width=\linewidth]{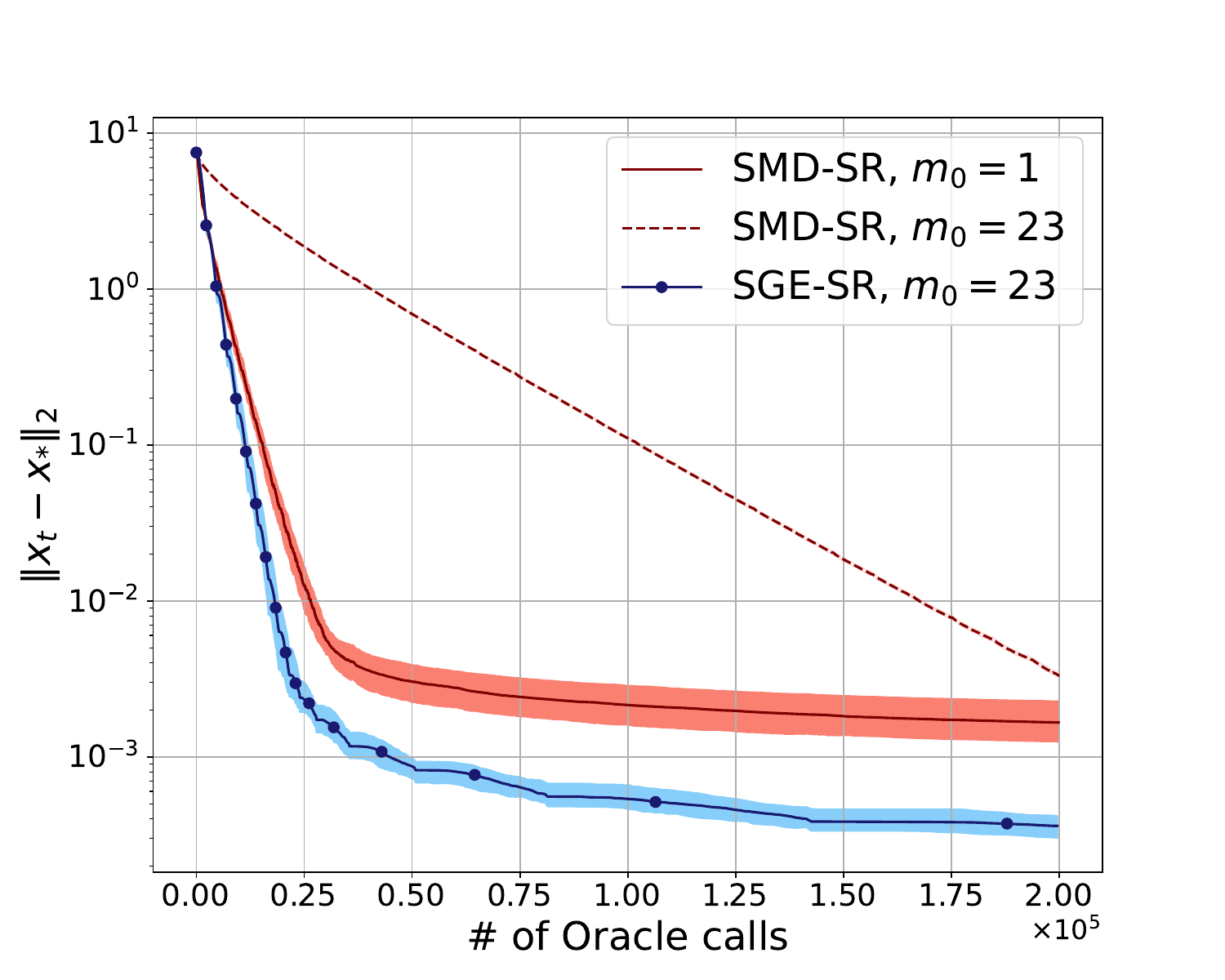}

  \end{subfigure}
  \medskip
  \centering
  \begin{subfigure}[b]{0.32\linewidth}
    \includegraphics[width=\linewidth]{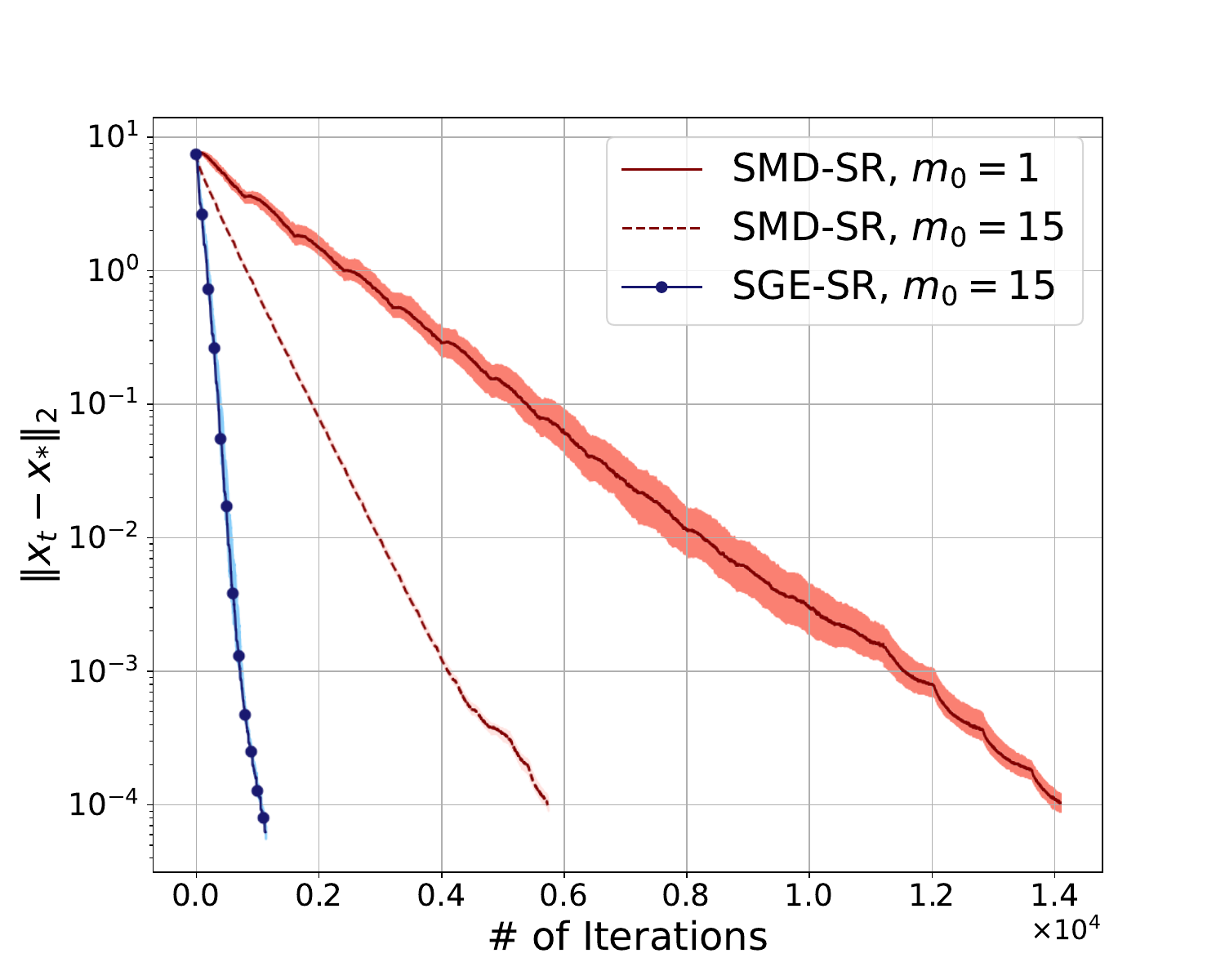}

  \end{subfigure}
  \begin{subfigure}[b]{0.32\linewidth}
    \includegraphics[width=\linewidth]{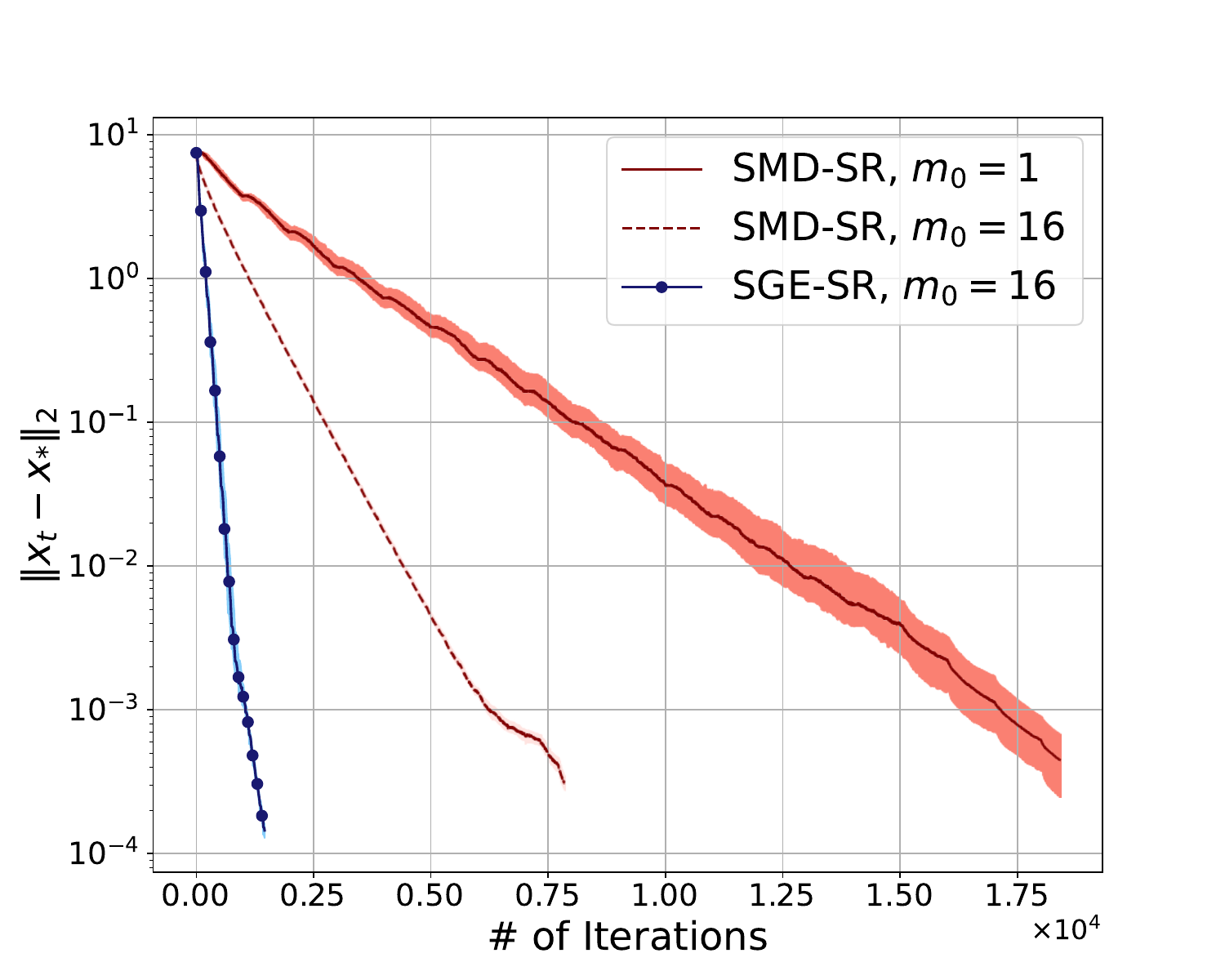}

  \end{subfigure}
  \begin{subfigure}[b]{0.32\linewidth}
    \includegraphics[width=\linewidth]{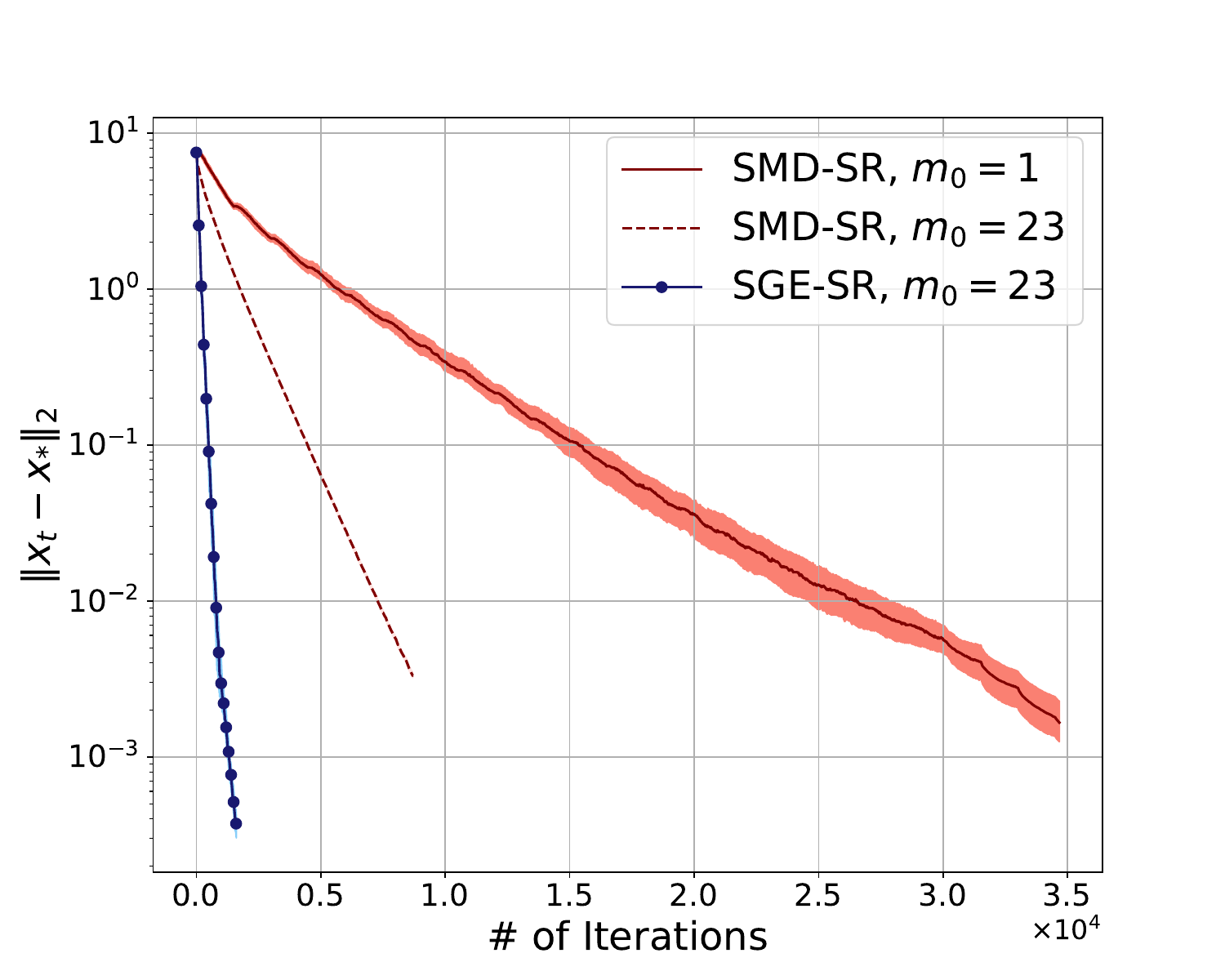}
  \end{subfigure}
  \caption{SGE-SR compared to ``vanilla'' SMD-SR  and its mini-batch variant. First row:  error $\|x_t -x^*\|_2$ against the number of oracle calls; second row: the error against the number of algorithm iterations.   Figure columns correspond to the results for  $u_1$, $u_{1/2}$, and $u_{1/10}$ activation functions and $\sigma=0.001$.}
  \label{fig:SGE-SR vs SMD-SR batch-size experiment}

\end{figure}

Finally, we illustrate the performance of SGE-SR and SMD-SR algorithms which share the same size of mini-batch. Recall that both  algorithms if ``normally set''---SGE-SR with mini-batch of optimal size and SMD-SR with ``trivial'' mini-batch (of size $m_0=1$)---converge linearly during the preliminary phase. In Figure \ref{fig:SGE-SR vs SMD-SR batch-size experiment} we report on the simulation of algorithms with the same  (optimal for the accelerated algorithm) size of the mini-batch \cite{juditsky2023sparse}.

% \textcolor{purple}{\textbf{Sasila : }
The series of experiments conducted indicates that the SGE-SR algorithm outperforms its non-accelerated counterpart. Despite both algorithms exhibiting linear rates of convergence, SGE-SR reaches a better precision for a fixed number of samples in every setting. This advantage is also observed when the two algorithms are compared in terms of the number of iterations, where the accelerated algorithm clearly outperforms its non-accelerated counterpart by a significant margin. The improved iteration complexities make SGE-SR a viable solution for distributed sparse recovery problems, where it is crucial to reduce the number of communication rounds between the servers and the clients while keeping high precision.
% }

% We believe that these numerical results speak for themselves ({\bf GL: I think it is worth to summarize the findings in words, rather than ask the readers to do the job.}).

%Finally, Figure 5 compares the accelerated method with the Euclidean distance generating function $\omega(x) = c\|x\|_2^2$ against the non-Euclidean variant. Figure \ref{fig:SGE-SR vs SMD-SR}, \ref{fig:SGE-SR vs SMD-SR batch-size experiment} and \ref{fig:SGE-SR different condition number} share the same setting i.e., $n = 100000$, $N \leq 200000$ and $s = 50$. In Figure 5, we vary the dimension of the problem but fix at the same values the maximum sample size $N$ and the sparsity level $s$.

%
%
% ({\bf GL: I recommend to add a short concluding remarks section. Here we can also mention some futher directions that will possibly relate to the second part of this series.})

\section{Concluding remarks}
In this paper, we investigate the problem of stochastic smooth convex optimization with ``state-dependent'' variance of stochastic gradients. We study two non-Euclidean accelerated stochastic approximation algorithms, stochastic accelerated gradient descent (SAGD) and stochastic gradient extrapolation (SGE), and provide optimal iteration and sample complexities for both algorithms under appropriate conditions. However, the optimal convergence guarantees for SGE require less restrictive assumptions, thus leading to wider applications such as statistical estimation problems with heavy tail noises. In addition, we propose a multi-stage routine of SGE to solve problems that satisfy the quadratic growth condition and further extend it to the sparse recovery problem. Our theoretical guarantees are corroborated by numerical experiments in high-dimensional settings. Further research will be directed to proving large deviation bounds to ensure the reliability and robustness of the solutions.

\appendix
\section{Appendix: proofs}
\subsection{Proof of Lemmas {\ref{assump:expectation_revised}} and \ref{assump:expectation_2}}
\paragraph{Proof of Lemma \ref{assump:expectation_revised}}
\tli{}{
To save notation, we put
    \begin{align*}
         \delta_{t}= G_t({u}) - g({u}), \quad \delta_{t,i}= \cG_t({u}, \xi_{t, i}) - g({u}).
    \end{align*}
    Recall that, by definition,
    \begin{align*}
    V_{x, x_0}^*(y) = \max_{z \in E} \left\{\langle y, z - x \rangle - V_{x_0}(x,z) \right\} = \max_{z \in E} \left\{\langle y + \nabla \omega (x-x_0), z - x \rangle - \omega(z-x_0) + \omega(x-x_0) \right\}.
    \end{align*}
    Then, by the $1$-strong convexity of $\omega$ w.r.t. the $\|\cdot\|$-norm, we have that $V_{x, x_0}^*$ is smooth with $1$-Lipschitz continuous gradient w.r.t. the dual norm  $\|\cdot\|_*$. Thus,
    \begin{align*}
        V_{x, x_0}^*\big(\gamma\delta_t\big) \leq V_{x, x_0}^*\left(\tfrac{\gamma}{m_t} \sum_{i=1}^{m_t - 1} \delta_{t, i}\right)  + \left\langle \nabla V_{x, x_0}^*\left(\tfrac{\gamma}{m_t} \sum_{i=1}^{m_t - 1} \delta_{t, i}\right), \tfrac{\gamma\delta_{t, m_t}}{m_t}  \right\rangle + \tfrac{1}{2}\left\|\tfrac{\gamma\delta_{t, m_t}}{m_t} \right\|_*^2.
    \end{align*}
When using the above relationship recursively, by the independence of $\delta_{t, i},~i\in [m_t]$ and the fact that $V_{x, x_0}^*(0) = 0$, we get
    $$
    \aic{\bbec{t-1}}{\bbe}[V_{x, x_0}^*\big( \gamma \delta_t\big)] \leq \tfrac{\gamma^2}{2m_t}  \aic{\bbec{t-1}}{\bbe}\|\delta_{t,1}\|_*^2,
    $$
which completes the proof.\qed}

\paragraph{Proof of Lemma \ref{assump:expectation_2}.}
% Let $\bar G_t := \tfrac{1}{m_t} \tsum_{i=1}^{m_t} \cG(\xavg_{t-1}, \xi_{t,i})$\textcolor{red}{(This quantity is introduced but is not used in the proof)}.
\tli{}{
Note that $\xavg_{t-1}$ and $\xmid_t$ are $\mathcal{F}_{t-1}$-measurable.
By Lemma \ref{assump:expectation_revised},
\begin{align*}
&\bbec{t-1}\big[V^*_{z_{t-1}}(\gamma[G_t - g(\xmid_t)])\big]\leq \tfrac{\gamma^2}{m_t}\bbec{t-1}\big[\|\cG(\xmid_t, \xi_{t,1}) - g(\xmid_t)\|_*^2\big]\nn\\
&\leq\tfrac{3{\gamma^2}}{m_t}\Big\{\bbec{t-1}\big[\|\cG(\xmid_t, \xi_{t,1}) - \cG(\xavg_{t-1}, \xi_{t,1})\|_*^2\big] + \bbec{t-1}\big[\|\cG(\xavg_{t-1}, \xi_{t,1}) - g(\xavg_{t-1})\|_*^2\big]\nn\\&\qquad\qquad+ \bbec{t-1}\big[\| g(\xavg_{t-1}) - g(\xmid_t)\|_*^2\big]\Big\}\nn\\
&\leq \tfrac{3{\gamma^2}}{m_t}\left\{\left(\bbe_{\xi_{t,1}}[\cK(\xi_{t,1})^2] + L^2\right) \|\xmid_t -\xavg_{t-1} \|^2 +\big[\cL[f(\xavg_{t-1}) - f^* - \langle g(x^*), \xavg_{t-1} - x^*\rangle] + \sigma_*^2\big]\right\}%,\label{error_decompo}
\end{align*}
(the last inequality being the consequence of Assumptions~\eqref{assump:variance_1} and \eqref{assump:variance}).
When utilizing $\xavg_{t-1}- \xmid_t = \xavg_{t-1}-(1-\beta_t)\xavg_{t-1} - \beta_t z_{t-1} = \beta_t (\xavg_{t-1}-z_{t-1})$,
we obtain the desired result.\qed}

\subsection{Proof of Theorem \ref{thm:SAGD_2}}
For \revision{sack}{sake} of completeness, we start with the following statement similar to the previous results for ASGD (e.g., Proposition 3.1 in \cite{LanBook2020}).

% We start with the following technical statement similar to .
%The following proposition characterizes some convergence properties of the SAGD.

% \subsection{Analysis}
\begin{proposition}\label{pro:sagd1}
Let $\{z_t\}, \{\xmid_t\}$ and $\{\xavg_t\}$ be generated by Algorithm~\ref{alg:SAGD}. Suppose that $\{\beta_t\}$ and $\{\eta_t\}$ satisfy \rf{SAGD_cond0} for some $\theta_t \geq 0$.
%\begin{align}
%    \theta_t \beta_t \eta_t &\leq \theta_{t-1} \beta_{t-1} \eta_{t-1}, ~t=2,..., k \label{SAGD_cond_2}\\
%    \eta_t &> L\beta_t, ~t=1,..., k \label{SAGD_cond_3}
%\end{align}
Then for any $x \in X$ one has
\begin{align}
    & \sum_{t=1}^{k} \theta_t[f(\xavg_t) - f(x)] + \theta_k\beta_k \eta_k V(z_k,x)\leq \sum_{t=1}^k \theta_t (1-\beta_t)[f(\xavg_{t-1}) - f(x)]+ \theta_1 \beta_1\eta_1 V(z_0,x) \nn\\&\qquad+ \sum_{t=1}^k \theta_t \beta_t \langle\delta_t, x - z_{t-1} \rangle +
    \sum_{t=1}^k {\theta_t\beta_t (\eta_t - L\beta_t)V^*_{z_{t-1}}\left(-\tfrac{\delta_t}{\eta_t - L\beta_t}\right)} \label{eq:prop_SAGD}
\end{align}
where $\delta_t := G_t - g(\xmid_t)$ and $V^*_z(\cdot)$ is defined in \rf{vstar}.
\end{proposition}
%\begin{proof}
\paragraph{Proof of the proposition.}{}%
First, by  convexity of $f$ and due to the definition of $\xavg_t$, we have for all $z \in X$,
\begin{align*}
    f(z) + \langle g(z), \xavg_t - z \rangle &= (1-\beta_t)[f(z) + \langle g(z), \xavg_{t-1}-z\rangle ] + \beta_t [f(z) + \langle g(z), z_t -z\rangle]\\
    & \leq (1-\beta_t) f(\xavg_{t-1})+ \beta_t [f(z) + \langle g(z), z_t -z\rangle].
\end{align*}
Now, by the smoothness of $f$,
\begin{align*}
    f(\xavg_t) &\leq f(z) + \langle g(z), \xavg_t - z\rangle  + \tfrac{L}{2}\|\xavg_t - z\|^2\\
    & \leq (1-\beta_t) f(\xavg_{t-1}) + \beta_t[f(z) + \langle g(z), z_t - z\rangle] + \tfrac{L}{2}\|\xavg_t - z\|^2,
\end{align*}
so that for $z = \xmid_t$ we have
\begin{align}\label{SAGD_eq1}
f(\xavg_t) \leq (1-\beta_t) f(\xavg_{t-1}) + \beta_t[f(\xmid_t) + \langle g(\xmid_t), z_t - \xmid_t\rangle] + \tfrac{L}{2}\|\xavg_t - \xmid_t\|^2.
\end{align}
On the other hand, the optimality condition for \eqref{eq:ad_SAGD2} in Algorithm~\ref{alg:SAGD} implies the following relationship (see, e.g., Lemma 3.5 of \cite{LanBook2020}):
% cf.  \cite{chen1993convergence} ({\bf GL: what is this reference? Three-point lemma?})
\begin{align}\label{SAGD_eq2}
\langle G_t, z_t - \xmid_t\rangle + \eta_t V(z_{t-1}, z_t) + \eta_t V(z_t,x) \leq \langle G_t, x - \xmid_t \rangle + \eta_t V(z_{t-1}, x), \quad \forall x\in X.
\end{align}
By combining \eqref{SAGD_eq1} and \eqref{SAGD_eq2} we obtain for all $x\in X$,
\begin{align}\label{SAGD_eq3}
&f(\xavg_t) + \beta_t \eta_t V(z_{t-1},z_t) + \beta_t \eta_t V(z_t, x)\nn\\
&\qquad\leq (1-\beta_t) f(\xavg_{t-1}) + \beta_t [f(\xmid_t) + \langle g(\xmid_t), x - \xmid_t\rangle]\nn\\
&\qquad\quad + \beta_t \langle \delta_t, x - z_t\rangle + \beta_t \eta_t V(z_{t-1},x) + \tfrac{L}{2}\|\xavg_t - \xmid_t\|^2\nn\\
&\qquad\leq (1-\beta_t) f(\xavg_{t-1}) + \beta_t [f(\xmid_t) + \langle g(\xmid_t), x - \xmid_t\rangle]\nn\\
&\qquad\quad + \beta_t \langle \delta_t, z_{t-1} - z_t\rangle + \beta_t \langle \delta_t, x - z_{t-1}\rangle + \beta_t \eta_t V(z_{t-1},x) + \tfrac{L}{2}\|\xavg_t - \xmid_t\|^2.
\end{align}
From \eqref{eq:ad_SAGD1} and \eqref{eq:ad_SAGD3} we have $\xavg_t - \xmid_t = \beta_t(z_t - z_{t-1})$; after substituting into \eqref{SAGD_eq3} and taking into account \eqref{bound_bregman}, we obtain
\begin{align*}
    &f(\xavg_t) + \beta_t \eta_t V(z_t, x)\nn\\
&\quad\leq (1-\beta_t) f(\xavg_{t-1}) + \beta_t [f(\xmid_t) + \langle g(\xmid_t), x - \xmid_t\rangle]\nn\\
&\quad\quad + \beta_t \langle \delta_t, x - z_{t-1}\rangle + \beta_t \eta_t V(z_{t-1},x)+ {\beta_t \langle \delta_t, z_{t-1}-z_t \rangle} +
{(L\beta_t^2 - \eta_t \beta_t) V(z_{t-1}, z_t)}\nn\\
&\quad\leq (1-\beta_t) f(\xavg_{t-1}) + \beta_t [f(\xmid_t) + \langle g(\xmid_t), x - \xmid_t\rangle] + \beta_t \langle \delta_t, x - z_{t-1}\rangle + \beta_t \eta_t V(z_{t-1},x)+ {\beta_t (\eta_t - L\beta_t)V^*_{z_{t-1}}\left(-\tfrac{\delta_t}{\eta_t - L\beta_t}\right)},
\end{align*}
the last inequality being a consequence of \eqref{eq:ad_SAGD3} and the definition of $V^*_x(\cdot)$.  When subtracting $f(x)$ on both sides we get
\begin{align*}
&[f(\xavg_t) -f(x)] + \beta_t \eta_t V(z_t, x)\\
&\leq (1-\beta_t)[f(\xavg_{t-1})-f(x)] + \beta_t [\overbrace{f(\xmid_t) + \langle g(\xmid_t), x - \xmid_t\rangle - f(x)}^{\leq 0}]+ \beta_t \eta_t V(z_{t-1}, x)\\
&\qquad+ \beta_t \langle \delta_t, x - z_{t-1}\rangle+ {\beta_t (\eta_t - L\beta_t)V^*_{z_{t-1}}\left(-\tfrac{\delta_t}{\eta_t - L\beta_t}\right)}\\
& \leq (1-\beta_t)[f(\xavg_{t-1})-f(x)]+ \beta_t \eta_t V(z_{t-1}, x)+ \beta_t \langle \delta_t, x - z_{t-1}\rangle + {\beta_t (\eta_t - L\beta_t)V^*_{z_{t-1}}\left(-\tfrac{\delta_t}{\eta_t - L\beta_t}\right)}.
\end{align*}
After multiplying by $\theta_t$ and summing up from $t=1$ to $k$ we arrive at % and invoking condition \eqref{SAGD_cond_2}, we obtain
\begin{align*}
    &\sum_{t=1}^{k} \theta_t[f(\xavg_t) - f(x)] + \theta_k\beta_k \eta_k V(z_k,x)\leq \sum_{t=1}^k \theta_t (1-\beta_t)[f(\xavg_{t-1}) - f(x)] + \theta_1 \beta_1\eta_1 V(z_0,x) \\&\qquad+ \sum_{t=1}^k \theta_t \beta_t \langle\delta_t, x - z_{t-1} \rangle + \sum_{t=1}^k {\theta_t\beta_t (\eta_t - L\beta_t)V^*_{z_{t-1}}\left(-\tfrac{\delta_t}{\eta_t - L\beta_t}\right)},
\end{align*}
which is \rf{eq:prop_SAGD}.\qed
%\end{proof}

\paragraph{Proof of the theorem.}
When setting $x = x^*$ and taking expectation on both sides of \eqref{eq:prop_SAGD}, using {\eqref{eq:assum_minibatch_revised} in Lemma~\ref{assump:expectation_revised}},  we get
\begin{align*}
&\sum_{t=1}^k \theta_t\bbe[f(\xavg_t) -f^*] + \theta_k\beta_k \eta_k \bbe[V(z_k, x^*)]\\
% &\leq (1-\beta_t)\bbe[f(\xavg_{t-1})-f^*] + \beta_t \bbe[\underbrace{f(\xmid_t) + \langle g(\xmid_t), x^* - \xmid_t\rangle - f^*}_{\leq 0}]+ \beta_t \eta_t \bbe[V(z_{t-1}, x^*)]+ \tfrac{\beta_t \bbe\|\delta_t\|_*^2}{2(\eta_t-L\beta_t)}\\
&\quad  \leq \sum_{t=1}^k \theta_t(1-\beta_t)\bbe[f(\xavg_{t-1})-f^*]+ \theta_1\beta_1 \eta_1 V(z_0, x^*)+ \sum_{t=1}^k\theta_t r_t \bbe\big[\cL\big(f(\xmid_t) - f^* - \langle g(x^*), \xmid_t - x^* \rangle \big) + \sigma_*^2 \big],
\end{align*}
where $r_t$ is defined in Theorem~\ref{thm:SAGD_2}.
Recall that $\xmid_t = (1-\beta_t) \xavg_{t-1} + \beta_t z_{t-1}$. On the other hand, by smoothness of $f$, when recalling that $\langle g(x^*), \xavg_{t-1} - x^* \rangle\geq 0$,
\begin{align*}
    &f(\xmid_t) - f^* - \langle g(x^*), \xmid_t - x^* \rangle\\
    &\quad\leq (1-\beta_t)\big[ f(\xavg_{t-1}) - f^* - \langle g(x^*), \xavg_{t-1} - x^* \rangle\big] + \beta_t \big[ f(z_{t-1}) - f^* - \langle g(x^*), z_{t-1} - x^* \rangle\big]\\
    &\quad\leq (1-\beta_t)\big[ f(\xavg_{t-1}) - f^* - \langle g(x^*), \xavg_{t-1} - x^* \rangle\big] + \tfrac{L\beta_t}{2}\|z_{t-1}-x^*\|^2\\
    &\quad\leq (1-\beta_t)\big[ f(\xavg_{t-1}) - f^*\big] + L\beta_t V(z_{t-1}, x^*).
\end{align*}
When combining the above inequalities we obtain
\begin{align*}
    &\sum_{t=1}^k \theta_t\bbe[f(\xavg_t) -f^*]+ \theta_k\beta_k \eta_k \bbe[V(z_k, x^*)]\\
    & \leq \sum_{t=1}^k\theta_t(1-\beta_t)(1+r_t\cL) \bbe[f(\xmid_{t-1}) - f^*] + \theta_1\beta_1\eta_1 V(z_0,x^*)+\sum_{t=1}^k \left\{\theta_t r_t \beta_t L\cL\bbe[V(z_{t-1}, x^*)] + \theta_t r_t \sigma_*^2\right\}.
\end{align*}
Due to \eqref{SAGD_cond_1}, and taking into account that $\beta_1 = 1$, we conclude that
\begin{align*}
    \theta_k \bbe[f(\xavg_k) &- f^*] + \theta_k \beta_k \eta_k \bbe[V(z_k, x^*)]\leq \theta_1\beta_1\eta_1 V(z_0,x^*) + \sum_{t=1}^k \theta_t r_t \beta_t L \cL \bbe[V(z_{t-1}, x^*)] + \sum_{t=1}^k \theta_t r_t \sigma_*^2,
\end{align*}
what is  \eqref{eq:thm_SAGD_1}.\qed
\subsection{Proof of Corollary \ref{cor:sagd2}}
\paragraph{1$^o$.} Observe first that for $\eta\geq 4L$ and $m_t = m$ one has
\begin{align}\label{bound_r_t_2}
   r_t = \tfrac{\beta_t \tli{\revision{\Omega}{\newOmega}}{}}{2[\eta/(t+1) - 3L/(t+2)]m} %= \tfrac{\beta_t \revision{\Omega}{\newOmega}}{2[\eta/(4t+4)+ 3\eta/(4t+4) - 3L/(t+2)]m}
   \leq \tfrac{2\beta_t \tli{\revision{\Omega}{\newOmega}}{} (t+1)}{\eta m}=\tfrac{2\beta_t \tli{\revision{\Omega}{\newOmega}}{}}{\eta_t m}.
\end{align}
As a result, when $\eta \geq \tfrac{6\tli{\revision{\Omega}{\newOmega}}{} (k-1 )\cL}{m}$ and given the definition of $\beta_t$ and $\theta_t$, we have
\begin{align*}
    \theta_t (1-\beta_t) (1+r_t \cL) \leq (t+1)(t-1)\left(1 + \tfrac{6 \tli{\revision{\Omega}{\newOmega}}{} (t+1) \cL}{(t+2)\eta m}\right) \leq (t+1)(t-1)\left(1 + \tfrac{(t+1)}{(t+2)(k-1)}\right) \leq t(t+1) = \theta_{t-1},
\end{align*}
so that \eqref{SAGD_cond_1} holds. On the other hand,
\begin{align}\label{theta_beta_eta}
    \theta_t \beta_t\eta_t = (t+1)(t+2) \, \tfrac{3 \eta}{(t+2)(t+1)} = 3\eta = \theta_{t-1} \beta_{t-1}\eta_{t-1},
\end{align}
so condition~\eqref{SAGD_cond_2} is satisfied, so Theorem \ref{thm:SAGD_2} applies.  By \eqref{bound_r_t_2}, we also have
\begin{align}\label{bound_r_t_3}
    \theta_t r_t \beta_t L \cL \leq \tfrac{2\theta_t\beta_t^2\tli{\revision{\Omega}{\newOmega}}{}L \cL}{\eta_t m} \leq \tfrac{18\tli{\revision{\Omega}{\newOmega}}{} L \cL}{\eta_t m}.
\end{align}

\paragraph{2$^o$.} Next, let us check that $\bbe[V(z_t, x^*)] \leq 3D ^2$ for all $t \geq 1$. Indeed, for $t=1$ we have by \eqref{eq:thm_SAGD_1}, the definition of $\theta_t, \beta_t, \eta_t$ and due to  \eqref{bound_r_t_2}--\eqref{bound_r_t_3}
\begin{align*}
\bbe[V(z_1, x^*)] &\leq V(z_0,x^*) + \tfrac{\theta_1 \beta_1  r_1L \cL}{\theta_1 \beta_1\eta_1} V(z_0,x^*) + \tfrac{\theta_1 r_1}{\theta_1 \beta_1\eta_1} \sigma_*^2\leq V(z_0,x^*) + \tfrac{{8} \tli{\revision{\Omega}{\newOmega}}{} L \cL}{m\eta^2} V(z_0,x^*) + \tfrac{8\tli{\revision{\Omega}{\newOmega}}{}}{\eta^2 m } \sigma_*^2\\
& \leq V(z_0, x^*) + V(z_0,x^*) + D ^2 \leq 3D ^2.
\end{align*}
where the last inequality follows from  $\eta^2 \geq \max\left\{{\tfrac{9(k+{1})^2 \tli{\revision{\Omega}{\newOmega}}{}L \cL}{m}}, {\tfrac{2(k+2)^3 \tli{\revision{\Omega}{\newOmega}}{} \sigma_*^2}{3 D ^2 m}}\right\}$.

Now we assume that for $s=1,..., t-1$, $\bbe[V(z_{s}, x^*)] \leq 3D ^2$. Then
\begin{align*}
    \bbe[V(z_t, x^*)]&\leq \tfrac{\theta_1\beta_1\eta_1}{\theta_t \beta_t \eta_t } V(z_0,x^*) + \sum_{s=1}^t \tfrac{\theta_s r_s \beta_s L \cL}{\theta_t \beta_t \eta_t } \bbe[V(z_{s-1}, x^*)] + \sum_{s=1}^t \tfrac{\theta_s r_s}{\theta_t \beta_t \eta_t} \sigma_*^2\\
    &\leq V(z_0,x^*) +\sum_{s=1}^t \tfrac{r_s}{\eta_s }  L \cL\bbe[V(z_{s-1}, x^*)] + \sum_{s=1}^t \tfrac{r_s}{\beta_t\beta_s \eta_s} \sigma_*^2\\
    &\leq V(z_0,x^*) + \tfrac{3 (t+{1})^2 \tli{\revision{\Omega}{\newOmega}}{} L \cL }{\eta^2 m} \cdot 3 D ^2 + \tfrac{2 (t+2)^3 }{3\eta^2 m } \sigma_*^2\\
    &\leq V(z_0,x^*) + D ^2 + D ^2 \leq 3D ^2
\end{align*}
where the last inequality is, again, a consequence of $\eta^2 \geq \max\left\{{\tfrac{9(k+{1})^2 \tli{\revision{\Omega}{\newOmega}}{} L \cL}{m}}, {\tfrac{2(k+2)^3 \tli{\revision{\Omega}{\newOmega}}{}\sigma_*^2}{3 D ^2 m}}\right\}$.
\paragraph{3$^o$.} Finally, when substituting the above estimates into~\eqref{eq:thm_SAGD_1}, we obtain
\begin{align*}
    \bbe[f(\xavg_k) - f^*] &\leq \tfrac{1}{\theta_k}\left(\theta_1\beta_1 \eta_1 V(z_0,x^*) + \sum_{t=1}^k 3\theta_t r_t \beta_t L \cL D ^2  + \sum_{t=1}^k \theta_t r_t\sigma_*^2\right)\\
    & \overset{(i)}\leq \tfrac{1}{(k+1)(k+2)} \left( 3\eta D ^2 + \tfrac{27 (k+{1})^2\tli{\revision{\Omega}{\newOmega}}{} L \cL D ^2}{\eta m}  + \tfrac{2(k+2)^3 \tli{\revision{\Omega}{\newOmega}}{} \sigma_*^2}{\eta m} \right)\\
    & \overset{(ii)}\leq \tfrac{12LD ^2}{(k+1)(k+2)}+ \tfrac{18 \tli{\revision{\Omega}{\newOmega}}{} \cL D ^2}{(k+2) m} + \tfrac{18\sqrt{\tli{\revision{\Omega}{\newOmega}}{} L \cL}D ^2}{(k+1)\sqrt{m}} + \tfrac{2\sqrt{6(k+2) \tli{\revision{\Omega}{\newOmega}}{}\sigma_*^2D ^2}}{(k+1) \sqrt{m}}
\end{align*}
where (i) follows from \eqref{bound_r_t_2}--\eqref{bound_r_t_3}, and (ii) is a consequence of the definition of $\eta$. This implies \rf{eq:thm_SAGD_2} due to $\tfrac{k+2}{k+1}\leq \tfrac{4}{3}$ for $k\geq 2$. \qed

\subsection{Proof of Theorem \ref{thm:SAGD_1} and Corollary \ref{cor:sagd22}}
The subsequent proofs follow those of  Theorem \ref{thm:SAGD_2} and Corollary \ref{cor:sagd2}. We present them here for reader's convenience.
\paragraph{Proof of Theorem \ref{thm:SAGD_1}.}
Note that assumptions of Proposition \ref{pro:sagd1} hold. When taking expectation on both sides of \eqref{eq:prop_SAGD} and using the bound
Lemma~\ref{assump:expectation_2}, we obtain for $x=x_*$,
\begin{align*}
&\sum_{t=1}^k \theta_t\bbe[f(\xavg_t) -f^*] + \theta_k\beta_k \eta_k \bbe[V(z_k, x^*)]\\
% &\leq (1-\beta_t)\bbe[f(\xavg_{t-1})-f^*] + \beta_t \bbe[\underbrace{f(\xmid_t) + \langle g(\xmid_t), x^* - \xmid_t\rangle - f^*}_{\leq 0}]+ \beta_t \eta_t \bbe[V(z_{t-1}, x^*)]+ \tfrac{\beta_t \bbe\|\delta_t\|_*^2}{2(\eta_t-L\beta_t)}\\
% &\quad  \leq \tsum_{t=1}^k \theta_t(1-\beta_t)\bbe[f(\xavg_{t-1})-f^*]+ \theta_1\beta_1 \eta_1 V(z_0, x^*)\\
% &\quad  + \tsum_{t=1}^k\theta_t r_t \bbe\big[\cK^2\|\xavg_{t-1} - \xmid_t\|^2+ \cK\big(f(\xavg_{t-1}) - f^* - \langle g(x^*), \xavg_{t-1} - x^* \rangle \big) + \sigma_*^2 \big]\\
&\leq \sum_{t=1}^k \theta_t(1-\beta_t)\bbe[f(\xavg_{t-1})-f^*]+ \theta_1\beta_1 \eta_1 V(z_0, x^*)\\
&\qquad  + \sum_{t=1}^k\theta_t r_t \bbe\big[\revision{\cK}{\bar \cK}^2\beta_t^2\|\xavg_{t-1} - z_{t-1}\|^2+ 3\cL\big(f(\xavg_{t-1}) - f^* - \langle g(x^*), \xavg_{t-1} - x^* \rangle \big) + 3\sigma_*^2 \big].
\end{align*}
% where step (i) follows from $\xavg_{t-1}- \xmid_t = \xavg_{t-1}-(1-\beta_t)\xavg_{t-1} - \beta_t z_{t-1} = \beta_t (\xavg_{t-1}-z_{t-1})$.
By rearranging the terms and utilizing \eqref{eq:SAGD_5} and $\beta_1=1$, we get
\begin{align*}
    &\theta_k \bbe[f(\xavg_k) - f^*] + \theta_k \beta_k \eta_k \bbe[V(z_k, x^*)]\\
    &\leq \theta_1\beta_1\eta_1 V(z_0,x^*) + 3\theta_1 r_1 \cL \big(f(z_{0}) - f^* - \langle g(x^*), z_{0} - x^* \rangle \big) + \sum_{t=1}^k \theta_t r_t \bbe\big[\revision{\cK}{\bar \cK}^2\beta_t^2 \bbe\|\xavg_{t-1} - z_{t-1}\|^2 + 3\sigma_*^2\big]\\
    &\leq \theta_1\beta_1\eta_1 V(z_0,x^*) + \tfrac{3\theta_1 r_1 \cL L }{2} \|z_0-x^*\|^2 + \sum_{t=1}^k \theta_t r_t \bbe\big[\revision{\cK}{\bar \cK}^2\beta_t^2 \bbe\|\xavg_{t-1} - z_{t-1}\|^2 + 3\sigma_*^2\big],
\end{align*}
what is \eqref{eq:thm_SAGD_3}.\qed
\paragraph{Proof of Corollary \ref{cor:sagd22}.}
Observe that for $\eta\geq 4L$ and $m_t = m$ (cf. \rf{bound_r_t_2})
\begin{align}\label{bound_r_t_2_1}
   r_t \leq \tfrac{2\beta_t \tli{\revision{\Omega}{\newOmega}}{} (t+1)}{\eta m}.
\end{align}
Then, due to $\eta \geq \tfrac{18\tli{\revision{\Omega}{\newOmega}}{}(k+1 )\cL}{m}$,
\begin{align*}
    \theta_t (1-\beta_t+3r_t \cL) \leq (t+1)(t+2)\left(\tfrac{t-1}{t+2} + \tfrac{18 \tli{\revision{\Omega}{\newOmega}}{} (t+1) \cL}{(t+2)\eta m}\right) \leq (t+1)\left(t-1 + \tfrac{t+1}{k+1}\right) \leq t(t+1) = \theta_{t-1},
\end{align*}
thus \eqref{eq:SAGD_5} holds. We have (cf. \rf{theta_beta_eta})
\begin{align*}%\label{bound_r_t_2_4}
    \theta_t \beta_t\eta_t = (t+1)(t+2) \cdot\tfrac{3 \eta}{(t+2)(t+1)} = 3\eta = \theta_{t-1} \beta_{t-1}\eta_{t-1},
\end{align*}
so \eqref{SAGD_cond_2} is verified and the conclusion of Theorem \ref{thm:SAGD_1} applies. We also have
\begin{align}\label{bound_r_t_2_3}
    \theta_t r_t \beta_t^2 \leq \tfrac{2\theta_t\beta_t^3\tli{\revision{\Omega}{\newOmega}}{}(t+1)}{\eta m} \leq \tfrac{54\tli{\revision{\Omega}{\newOmega}}{}}{\eta m}.
\end{align}
Let us now check that $\bbe[V(z_t, x^*)] \leq 3 D ^2$ for all $t \geq 1$. First,  for $t=1$, applying \eqref{eq:thm_SAGD_3} and taking into account \eqref{bound_r_t_2_1}--\eqref{bound_r_t_2_3},  we get
\begin{align*}
\bbe[V(z_1, x^*)] &\leq V(z_0,x^*) + \tfrac{3\theta_1 r_1 \cL L }{2\theta_1 \beta_1\eta_1} \|z_0-x^*\|^2 + \tfrac{3\theta_1 r_1}{\theta_1 \beta_1\eta_1} \sigma_*^2\\
& \leq V(z_0,x^*) + \tfrac{12 \tli{\revision{\Omega}{\newOmega}}{}\cL L }{m\eta^2} \|z_0-x^*\|^2 + \tfrac{24\tli{\revision{\Omega}{\newOmega}}{}}{\eta^2 m } \sigma_*^2\leq V(z_0, x^*) + \tfrac{1}{24}V(z_0,x^*) + D ^2 \leq 3D ^2,
\end{align*}
 the second inequality being due to \eqref{bound_bregman} and $\eta \geq
\max\left\{ 4L, \tfrac{18\tli{\revision{\Omega}{\newOmega}}{}(k+1 )\cL}{m}, \sqrt{\tfrac{2(k+2)^3 \tli{\revision{\Omega}{\newOmega}}{} \sigma_*^2}{ D ^2 m}}\right\}$.

Now, assume that for $s=1,..., t-1$, $\bbe[V(z_{s-1}, x^*)] \leq 3D ^2$. Because $\xavg_{s}$ is a weighted average of the previous iterates, by the convexity of $\|\cdot\|^2$, we have for $s <t$,
\begin{align*}
    \bbe\|\xavg_{s} - z_{s}\|^2 &\leq 2\bbe\|\xavg_{s} - x^*\|^2 + 2
    \bbe\|z_{s} - x^*\|^2\nn\\
    &\leq 2 \max_{0\leq i \leq s} \bbe\|z_i-x^*\|^2 +  2\bbe\|z_{s} - x^*\|^2\nn\\
    &\leq 4\max_{0\leq i \leq s} \bbe[V(z_i,x^*)] + 4 \bbe[V(z_{s},x^*)] \leq 24 D ^2.
 \end{align*}
As a consequence, when substituting into \eqref{eq:thm_SAGD_3} the bounds of \eqref{bound_r_t_2_1}--\eqref{bound_r_t_2_3}
 \begin{align*}
  \bbe[V(z_k, x^*)]
    & \leq  V(z_0,x^*) +  \tfrac{3\theta_1 r_1 \cL L}{2 \theta_k \beta_k \eta_k} \|z_0-x^*\|^2+ \sum_{t=2}^k \tfrac{24\theta_t r_t \beta_t^2 \cK^2 D ^2}{\theta_k \beta_k \eta_k}  + \sum_{t=1}^k \tfrac{3\theta_t r_t \sigma_*^2}{\theta_k \beta_k \eta_k}\\
    & \overset{(ii)}\leq V(z_0,x^*) +  \tfrac{12 \tli{\revision{\Omega}{\newOmega}}{} \cL L }{m\eta^2} \|z_0-x^*\|^2 + \tfrac{144 \tli{\revision{\Omega}{\newOmega}}{}(k-1) \cK^2 D ^2}{\eta^2 m} + \tfrac{2(k+2)^3 \tli{\revision{\Omega}{\newOmega}}{} \sigma_*^2}{\eta^2 m} \overset{(i)}\leq 3D ^2,
\end{align*}
due to $$\eta \geq \max\left\{ 4L, \tfrac{18\tli{\revision{\Omega}{\newOmega}}{}(k+1 )\cL}{m}, 12\sqrt{\tfrac{k \tli{\revision{\Omega}{\newOmega}}{} \cK^2}{m}}, \sqrt{\tfrac{2(k+2)^3 \tli{\revision{\Omega}{\newOmega}}{}\sigma_*^2}{ D ^2 m}}\right\}.$$

Finally,  from \eqref{eq:thm_SAGD_3} and the definition of $\eta$ we conclude that
\begin{align*}
    \bbe[f(\xavg_k) - f^*] &\leq \tfrac{1}{\theta_k} \left(\theta_1\beta_1\eta_1 V(z_0,x^*) +  \tfrac{3\theta_1 r_1 \cL L}{2} \|z_0-x^*\|^2+ \sum_{t=2}^k \theta_t r_t \beta_t^2 \cK^2 \bbe\|\bar z_{t-1} -z_{t-1}\|^2 + \sum_{t=1}^k 3\theta_t r_t \sigma_*^2\right)\\
    & \leq \tfrac{1}{(k+1)(k+2)} \left( 3\eta D ^2 + \tfrac{36 \tli{\revision{\Omega}{\newOmega}}{} \cL L}{m \eta} + \tfrac{432 (k-1) \tli{\revision{\Omega}{\newOmega}}{} \cK^2 D ^2}{\eta m}  + \tfrac{6(k+2)^3 \tli{\revision{\Omega}{\newOmega}}{} \sigma_*^2}{\eta m} \right)\\
    & \leq \tfrac{13LD ^2}{(k+1)(k+2)}+ \tfrac{54  \cL D ^2}{(k+2) m} + \tfrac{72 \cK D ^2}{(k+2)\sqrt{m (k+1)}} + \tfrac{6}{k+1 }  \sqrt{\tfrac{2(k+2) \tli{\revision{\Omega}{\newOmega}}{} \sigma_*^2D ^2}{m}}.
\end{align*}
This completes the proof due to $\tfrac{k+2}{k+1}\leq \tfrac{4}{3}$ for $k\geq 2$.\qed

\subsection{Proof of Theorem \ref{thm:main}}
For $0\leq \beta_t\leq 1$, denote $\tau_t=\tfrac{1-\beta_t}{ \beta_t}$. Relationships \eqref{eq:def_rellb} in variables $\theta_t,\alpha_t,\eta_t$, and $\tau_t$ become
\begin{subequations}
\label{eq:def_rell}
\begin{align}
\theta_{t-1} &= \alpha_t \theta_t, ~~\eta_t \le \alpha_t \eta_{t-1} , ~~ t= 2, \ldots,k \label{eq:def_eta_rel}\\
% % \Leftrightarrow  \eta_1\eta_2 \geq 25 L^2\alpha_2,\label{eq:def_alpha}
% \label{eq:def_Lip_1}\\
\tfrac{\eta_t \tau_{t-1}}{\alpha_{t}} &\ge{4}L,
% \Leftrightarrow \eta_t \tau_{t-1} \ge 5 L \alpha_t,
\quad t= 3, \ldots, k \label{eq:def_Lip_rel}\\
 \tfrac{\eta_1\eta_2}{\alpha_2} &\geq  {16}L^2, ~~\eta_k \tau_k  \ge L, \label{eq:def_Lip_rel1}
\end{align}
\end{subequations}
with \rf{eq:ad_SGE3} of Algorithm \ref{adaptive_GEM}
replaced with
\begin{equation}\label{eq:ad_SGE31}
x_t = (z_t + \tau_t x_{t-1})/(1+\tau_t).
\end{equation}
We need the following technical statement.
\begin{proposition}\label{prop:general_convergence}
Let the algorithm parameters $\{\alpha_t\}$, $\{\eta_t\}$ and $\{\tau_t\}$ satisfy $\tau_1 \geq 0,$ $\tau_t > 0$ for all $t\geq 2$, along with relations
\eqref{eq:def_eta_rel}--\eqref{eq:def_Lip_rel1}
 for some $\theta_t \ge 0$.
Then for all $x \in \feaReg$
\begin{align}
&\sum_{t=1}^k \theta_t \left\{  \tau_t [f(x_t) -f(x_{t-1})]- \langle g(x_t), x - x_t \rangle\right\}\nn\\
&\leq \theta_1  \eta_1  V(x_0, x) + \sum_{t=0}^{k-1}\left[  \tfrac{ \theta_{t+1} (1 + \alpha_{t+1}^2)}{4\eta_{t+1}} V^*_{z_{t}}(4\delta_{t}(x_{t})) + \tfrac{\theta_{t+2}\alpha_{t+2}^2}{4\eta_{t+2}}  V^*_{z_{t+1}}(4\delta_{t+1}(x_{t})) - \theta_{t+1} \langle \delta_t(x_t), z_t - x\rangle \right]\nn\\
&\qquad+ \sum_{t=1}^{k-1}\theta_{t+1}\alpha_{t+1} \langle \delta_t(x_{t-1}) - \delta_t(x_t), z_t - x \rangle.
\label{eq:prop_GEM_main_general}
\end{align}
where
% $p_0 = \theta_1+ \theta_2\alpha_2$,
% $~p_t:= \begin{cases}  ~\theta_t,&t\leq k-2,\\
% ~\theta_{t+1} + \theta_{t}, & t = k-1 \end{cases}$ and
$\delta_{t}(x_{s})= G_t(x_{s}) - g(x_{s})$ with $G_t(\cdot)$ as defined in Algorithm~\ref{adaptive_GEM}.

 \end{proposition}
\paragraph{Proof of the proposition.}
\paragraph{1$^o$.}
% Denote $g_t := g(x_t)$.
By the smoothness of $f$, we have
\[
\tfrac{1}{2L}\|g(x_t) - g(x_{t-1})\|_*^2 \le f(x_{t-1}) - [f(x_t) + \langle g(x_t), x_{t-1} - x_t\rangle].
\]
Then for any $x\in \feaReg$, using the above inequality and the fact that
$x -  x_t  = x -z_t -  \tau_t( x_{t-1} -  x_t)$ due to \eqref{eq:ad_SGE31},
we get
\begin{align*}
\tau_t f( x_t) - \langle g(x_t), x -  x_t \rangle
&= \tau_t [f( x_t) + \langle g(x_t),  x_{t-1} -  x_t \rangle] + \langle g(x_t), z_t - x \rangle  \\
&\le \tau_t [f( x_{t-1}) - \tfrac{1}{2L} \|g(x_t) - g(x_{t-1})\|_*^2] + \langle g(x_t), z_t - x \rangle .
\end{align*}
By the optimality condition of \eqref{eq:ad_SGE2}, we have
\begin{align*}
 &\langle \widetilde G_t, z_t - x \rangle \le  \eta_t V(z_{t-1}, x) - \eta_tV(z_t, x) - \eta_t V(z_{t-1}, z_t).
\end{align*}
%In addition, by the convexity of $h$, we have
%\[
%(1 + \tau_t) h( x_t) \le h(x_t) + \tau_t h( x_{t-1}).
%\]
Combining two previous inequalities, we obtain
\begin{align*}
& \tau_t f( x_t) - \langle {g(x_t)}, x -  x_t \rangle  - \tau_t f( x_{t-1})\\
%&(1+\tau_t) [\Psi( x_t) - \Psi(x)] - \tau_t [\Psi( x_{t-1}) -\Psi(x)]
&\le \eta_t V(z_{t-1}, x) - \eta_tV(z_t, x) + \langle {g(x_t)} - \widetilde G_t, z_t - x \rangle   -\tfrac{\tau_t}{2L} \|{g(x_t)-g(x_{t-1})}\|_*^2 - \eta_t V(z_{t-1}, z_t).
\end{align*}
Note that
\begin{align*}
\langle {g(x_t)} - \widetilde G_t, z_t - x \rangle
&= \langle {g(x_t) - g(x_{t-1}) - \alpha_t (g(x_{t-1}) - g(x_{t-2}))}, z_t - x \rangle  - \langle {\delta_{t-1}(x_{t-1}) + \alpha_t (\delta_{t-1}(x_{t-1}) - \delta_{t-1}(x_{t-2}))}, z_t - x\rangle\\
&= \langle {g(x_t) - g(x_{t-1})}, z_t - x\rangle - \alpha_t \langle  {g(x_{t-1}) - g(x_{t-2})}, z_{t-1} - x \rangle\\
&\quad + \alpha_t \langle {g(x_{t-1}) - g(x_{t-2})}, z_t - z_{t-1} \rangle- \langle {\delta_{t-1}(x_{t-1}) + \alpha_t (\delta_{t-1}(x_{t-1}) - \delta_{t-1}(x_{t-2}))}, z_t - x\rangle.
\end{align*}
When taking the $\theta_t$-weighted sum of the above inequalities for $t= 1, \ldots, k$,
recalling that $x_0 = z_0$, and using \eqref{eq:def_eta_rel}, we obtain
\begin{align}
&\sum_{t=1}^k \theta_t \left\{  \tau_t [f( x_t)- f( x_{t-1})] - \langle {g(x_t)}, x -  x_t \rangle  \right\} \nn\\
&\le \theta_1 \eta_1 V(x_0, x) - \theta_k \eta_k V(z_k, x) +\theta_k \langle {g(x_k) - g(x_{k-1})}, z_k - x \rangle   + \Delta_k \label{eq:GEM_converge_rel1}
\end{align}
where
\begin{align*}
\Delta_k &:= \sum_{t=1}^k \theta_t \big[\alpha_t \langle {g(x_{t-1}) - g(x_{t-2})}, z_t - z_{t-1} \rangle -\tfrac{\tau_t}{2L} \|{g(x_t) - g(x_{t-1})}\|_*^2\\
& ~\qquad\quad\qquad - \langle {\delta_{t-1}(x_{t-1}) + \alpha_t (\delta_{t-1}(x_{t-1}) - \delta_{t-1}(x_{t-2}))}, z_t - x\rangle  - \eta_t V(z_{t-1}, z_t)\big].
\end{align*}
Observe that
\[\theta_t\alpha_t\langle {g(x_{t-1})-g(x_{t-2})},z_t-z_{t-1}\rangle\leq \tfrac{\theta_t^2\alpha_t^2L}{ 2\theta_{t-1}\tau_{t-1}}\|z_t-z_{t-1}\|^2+\tfrac{\theta_{t-1}\tau_{t-1}}{ 2L}\|{g(x_{t-1})-g(x_{t-2})}\|_*^2,
\]
so that,
after rearranging the terms,
\begin{align*}
\Delta_k &\le \underbrace{\theta_2\alpha_2 \langle {g(x_1)-g(x_0)}, z_2-z_1\rangle +  \sum_{t=3}^k \tfrac{L \theta_{t}^2 \alpha_{t}^2}{2 \theta_{t-1} \tau_{t-1}} \|z_t - z_{t-1}\|^2 - \tfrac{\theta_k \tau_k}{2 L}\|{g(x_k) - g(x_{k-1})}\|_*^2 - {\tfrac{1}{4}}\sum_{t=1}^k \theta_t \eta_t V(z_{t-1},z_t)}_{=:\Delta_{k,1}} \\
& \quad \underbrace{-\sum_{t=1}^k \theta_t \langle {\delta_{t-1}(x_{t-1}) + \alpha_t (\delta_{t-1}(x_{t-1}) - \delta_{t-1}(x_{t-2}))}, z_t - x\rangle - {\tfrac{3}{4}}\sum_{t=1}^k \theta_t \eta_t V(z_{t-1},z_t)}_{=:\Delta_{k,2}}.
\end{align*}
\paragraph{2$^o$.} Let us bound $\Delta_{k,1}$. By the Young's inequality,
\[\alpha_2\langle {g(x_1)-g(x_0)},z_2-z_1\rangle\leq \tfrac{2\alpha_2^2 }{\eta_2}\|{g(x_1)-g(x_0)}\|_*^2+\tfrac{\eta_2}{8}\|z_2-z_1\|^2\leq \tfrac{{2}\alpha_2^2 }{\eta_2}\|{g(x_1)-g(x_0)}\|_*^2+\tfrac{\eta_2}{4}V(z_1,z_2).
\]
Thus,
\begin{align}\label{eq:bound_delta_1}
    \Delta_{k,1} &\leq \tfrac{{2}\theta_2\alpha_2^2 \|{g(x_1)-g(x_0)}\|_*^2}{\eta_2} - \tfrac{1}{{4}}\theta_1\eta_1 V(z_1,z_0)+\sum_{t=3}^k \left(\tfrac{L \theta_{t}^2 \alpha_{t}^2}{2 \theta_{t-1} \tau_{t-1}} - \tfrac{\theta_t \eta_t}{{8}}\right)\|z_t - z_{t-1}\|^2- \tfrac{\theta_k \tau_k}{2 L}\|{g(x_k) - g(x_{k-1})}\|_*^2\nn\\
    &\leq - \tfrac{\theta_k \tau_k}{2 L}\|{g(x_k) - g(x_{k-1})}\|_*^2
\end{align}
where the last inequality follows from \eqref{eq:def_Lip_rel} and the bound
\begin{align*}
	\tfrac{{2}\theta_2\alpha_2^2 \|{g(x_1)-g(x_0)}\|_*^2}{\tli{2}{}\eta_2}  \leq \tfrac{2\theta_2\alpha_2^2L^2\| x_1 -  x_0\|^2}{\eta_2} = \tfrac{2\theta_2\alpha_2^2L^2\|z_1 -  z_0\|^2}{\eta_2} \leq \tfrac{{4}\theta_2\alpha_2^2L^2V(z_1,z_0) }{\eta_2} \leq \tfrac{\theta_1 \eta_1 V(z_1,z_0)}{{4}}.
\end{align*}
Note that
\begin{align*}
&  -  \eta_k V(z_k, x)
+ \langle \tli{g_k - g_{k-1}}{g(x_k) - g(x_{k-1})}, z_k - x \rangle
- \tfrac{ \tau_k}{2 L}\|\tli{g_k - g_{k-1}}{g(x_k) - g(x_{k-1})}\|_*^2 \nn \\
&\qquad \le -\tfrac{\eta_k }{2} \|z_k-x\|^2 + \langle \tli{g_k - g_{k-1}}{g(x_k) - g(x_{k-1})}, z_k - x \rangle
- \tfrac{ \tau_k}{2 L}\|\tli{g_k - g_{k-1}}{g(x_k) - g(x_{k-1})}\|_*^2 \nn \\
&\qquad\le - (\tfrac{\eta_k }{2} - \tfrac{L}{2\tau_k}) \|z_k-x\|^2 \le 0
\end{align*}
due to the second relationship in \eqref{eq:def_Lip_rel1}.
Now, substituting the bound \rf{eq:bound_delta_1} into \eqref{eq:GEM_converge_rel1} results in
\begin{align}
&\sum_{t=1}^k \theta_t \left\{  \tau_t [f( x_t)-f( x_{t-1})] - \langle \tli{g_t}{g(x_t)}, x -  x_t \rangle\right\} \nn\\
&\qquad\le \theta_1 \eta_1V(x_0, x) - \theta_k \eta_k V(z_k, x) +\theta_k \langle \tli{g_k - g_{k-1}}{g(x_k) - g(x_{k-1})}, z_k - x \rangle  - \tfrac{\theta_k \tau_k}{2 L}\|\tli{g_k - g_{k-1}}{g(x_k) - g(x_{k-1})}\|_*^2 + \Delta_{k,2}\nn\\
&\qquad \le \theta_1 \eta_1 V(x_0, x)    + \Delta_{k,2}, \label{eq:GEM_converge_rel1_1}
\end{align}
\paragraph{3$^o$.}
To bound $\Delta_{k,2}$ we act as follows. Observe that
\begin{align*}
    \Delta_{k,2} &\le - \sum_{t=1}^k\theta_t [ \langle \tli{\delta_{t-1}}{\delta_{t-1}(x_{t-1})}, z_t - z_{t-1}\rangle +  \langle \tli{\delta_{t-1}}{\delta_{t-1}(x_{t-1})}, z_{t-1} - x\rangle]\\
&\quad - \sum_{t=2}^k \theta_t \alpha_t [\langle \tli{\delta_{t-1}}{\delta_{t-1}(x_{t-1})}, z_t - z_{t-1}\rangle +  \langle \tli{\delta_{t-1}}{\delta_{t-1}(x_{t-1})}, z_{t-1} - x\rangle]\\
&\quad +  \sum_{t=2}^k\theta_t \alpha_t [ \langle \tli{\delta_{t-2}}{\delta_{t-1}(x_{t-2})}, z_t - z_{t-1}\rangle + \tli{\langle \delta_{t-2}, z_{t-1} - z_{t-2}\rangle + \langle \delta_{t-2}, z_{t-2} - x\rangle}{\langle \delta_{t-1}(x_{t-2}), z_{t-1} - x\rangle}]\\
&\quad - {  \tfrac{3}{4}\sum_{t=1}^k \theta_t \eta_t V(z_{t-1},z_t)}.
\end{align*}
Recall that $z_t$ is ${\cal F}_{t-1}$-measurable. Recalling the definition of $V^*_x(y)$ in Lemma~\ref{assump:expectation_revised},  for $a>0$,
\[\langle \delta, z_r - z_{r-1}\rangle - a V(z_{r-1}, z_r) \le a V_{z_{r-1}}^{*}(\tfrac{\delta}{a}).
\]
As a result,
\begin{align}
  \Delta_{k,2} &\le  \sum_{t=1}^k [ \tli{\tfrac{ 5\theta_t}{2 \eta_t}\|\delta_{t-1}\|_*^2}{\tfrac{\theta_t}{4 \eta_t} V^*_{z_{t-1}}(4 \delta_{t-1}(x_{t-1}))} -  \theta_t\langle \tli{\delta_{t-1}}{\delta_{t-1}(x_{t-1})}, z_{t-1} - x\rangle]\nn\\
&\quad + \sum_{t=2}^k  [\tli{\tfrac{5 \theta_t \alpha_t^2}{2 \eta_t}  \|\delta_{t-1}\|_*^2}{\tfrac{\theta_t \alpha_t^2}{4 \eta_t}  V^*_{z_{t-1}}(4 \delta_{t-1}(x_{t-1}))} -  \theta_t \alpha_t \langle \tli{\delta_{t-1}}{\delta_{t-1}(x_{t-1})}, z_{t-1} - x\rangle] \nn\\
&\quad + \sum_{t=2}^k[ \tli{\tfrac{5 \theta_t \alpha_t^2  }{2\eta_t} \|\delta_{t-2}\|_*^2}{\tfrac{\theta_t \alpha_t^2}{4 \eta_t}  V^*_{z_{t-1}}(4 \delta_{t-1}(x_{t-2}))}
+ \tli{\tfrac{ 5\theta_t^2 \alpha_t^2}{2 \theta_{t-1} \eta_{t-1}}  \|\delta_{t-2}\|_*^2}{}
+ \theta_t \alpha_t\tli{\langle \delta_{t-2}, z_{t-2} - x\rangle}{\langle \delta_{t-1}(x_{t-2}), z_{t-1} - x\rangle}] \nn \\
&\le
  \tli{\sum_{t=1}^k\tfrac{5}{2} [ \tfrac{ \theta_t (1 + \alpha_t^2)}{\eta_t} +  \tfrac{ \theta_{t+1} \alpha_{t+1}^2  }{\eta_{t+1}} + \tfrac{ \theta_{t+1}^2 \alpha_{t+1}^2}{\theta_{t} \eta_{t}}]\|\delta_{t-1}\|_*^2 }
  {\sum_{t=1}^k \big[\tfrac{ \theta_t (1 + \alpha_t^2)}{4\eta_t}  V^*_{z_{t-1}}(4\delta_{t-1}(x_{t-1})) + \tfrac{\theta_{t+1}\alpha_{t+1}^2}{4\eta_{t+1}}  V^*_{z_{t}}(4\delta_{t}(x_{t-1})) \big]}\nn\\
&\quad -\sum_{t=1}^k\theta_t\tli{(1+\alpha_t)}{}\langle \delta_{t-1}, z_{t-1} - x\rangle + \tli{\sum_{t=1}^{k-1}\theta_{t} \langle \delta_{t-1}, z_{t-1} - x\rangle}{\sum_{t=2}^{k}\theta_{t}\alpha_t \langle \delta_{t-1}(x_{t-2})-\delta_{t-1}(x_{t-1}), z_{t-1} - x\rangle}.
\label{eq:bound_delta_2}
\end{align}
When substituting the bound \rf{eq:bound_delta_2} for $\Delta_{k,2} $ into \eqref{eq:GEM_converge_rel1_1} we obtain \rf{eq:prop_GEM_main_general}.
\qed
\paragraph{Proof of the theorem.}
By taking expectation on both sides of \eqref{eq:prop_GEM_main_general},
and using \eqref{assump:variance}, \tli{\eqref{eq:assum_minibatch}}{\eqref{eq:assum_minibatch_revised}} along with
the fact that
\begin{align*}
 \langle g(x_t), x - x_t \rangle   \le f(x) - f(x_t)
\end{align*}
we have for all $x \in \feaReg$,
\begin{align*}
&\sum_{t=1}^k \theta_t \bbe \left\{  \tau_t [f(x_t)-f( x_{t-1})] - [ f(x) - f(x_t)]\right\}\nn\\
&\le  \theta_1 \eta_1  V(x_0, x) + \sum_{t=1}^k\epsilon_{t-1} \left\{\cL  \bbe [f(x_{t-1}) - f^* - \langle g(x^*), x_{t-1} - x^* \rangle] + \sigma_*^2\right\},
\end{align*}
where $\epsilon_t$ is defined in \eqref{eq:def_epsilon_t_adaptive}.
Note that \rf{eq:batchsizeBt_adaptiveb} implies that
\beq \label{eq:batchsizeBt_adaptive}
\theta_t \tau_t + \cL  \epsilon_{t-1} \le \theta_{t-1}  (1+ \tau_{t-1}), ~ t\ge 2.
\eeq
Thus, by setting $x = x^*$ and using $\langle g(x^*), x_t - x^* \rangle\geq 0$, $f(x_0) - f^* - \langle g(x^*), x_0 - x^* \rangle \leq \tfrac{L}{2}\|x_0-x^*\|^2$, and $\tau_1 = 0$, we obtain
\begin{align*}
\theta_k  (1+ \tau_k)\bbe[ f(x_k) - f^*]&\leq \sum_{t=1}^k \theta_t  (1+ \tau_t)\bbe[ f(x_t) - f^*]
- \sum_{t=1}^{k-1} \theta_t  (1+ \tau_t)\bbe[ f(x_t) - f^*]  \\
&\leq \sum_{t=1}^k \theta_t  (1+ \tau_t)\bbe[ f(x_t) - f^*]
- \sum_{t=2}^k (\theta_t \tau_t + \epsilon_{t-1}\cL) \bbe [f(x_{t-1}) - f^*]  \\
&\qquad\le  \theta_1 \eta_1 V(x_0, x^*) +\tfrac{\epsilon_0\cL L}{2}\|x_0-x^*\|^2
% + \theta_1\tau_1[f(x_0) - f^*]
+\sum_{t=1}^k\epsilon_{t-1} \sigma_*^2
\end{align*}
which is \eqref{eq:GEM_theorem_result_general}.\qed

\subsection{Proof of Corollary \ref{coro_non_strongly_convex}}
\paragraph{1$^o$.}
Note that in the premise of the corollary one has $\tau_t={t-1\over 3}$. Let us check that with the present choice of stepsize parameters conditions \eqref{eq:def_eta_rel}--\eqref{eq:def_Lip_rel1} and \eqref{eq:batchsizeBt_adaptive} (which are equivalents \rf{eq:def_rellb} and \rf{eq:batchsizeBt_adaptiveb}) are satisfied. It is easy to see that \eqref{eq:def_eta_rel} holds. Because $\eta \geq {24}L$, we have
	\begin{align*}
		\eta_1\eta_2 & {=} \tfrac{1}{2} ({24}L)^2  \geq {16} \alpha_2 L^2,\\
		\eta_t \tau_{t-1} &{=}  \tfrac{\eta}{t} \tfrac{t-2}{3}
		\ge   \tfrac{{8}L(t-2)}{t} \geq {4}L \alpha_t,~t=3,...,k\\
		\eta_k \tau_k &\geq  \tfrac{\eta}{k}  \tfrac{k-1}{3} \ge \tfrac{{8}L(k-1) }{k} \ge L.
	\end{align*}
	To check \eqref{eq:batchsizeBt_adaptive}, notice that for $t \geq 1$,
    \begin{align*}
        q_t = \tli{\tfrac{\theta_{t+1}(1+\alpha_{t+1}^2)}{\eta_{t+1}} + \tfrac{\theta_{t+2} \alpha_{t+2}^2}{\eta_{t+2}} + \tfrac{\theta_{t+2}^2\alpha_{t+2}^2}{\theta_{t+1}\eta_{t+1}}}{\tfrac{\theta_{t+1}(1+\alpha_{t+1}^2)}{\eta_{t+1}} + \tfrac{\theta_{t+2} \alpha_{t+2}^2}{\eta_{t+2}}} = \tfrac{\tli{3}{2} (t+1)^2 + t^2}{\eta} \leq \tfrac{\tli{4}{3} (t+1)^2}{\eta},
    \end{align*}
thus    \begin{align*}
        \epsilon_t = \tfrac{\tli{5\revision{\Omega}{\newOmega}}{2}  q_t}{\tli{2}{} m_t} \leq \tfrac{\tli{10 \revision{\Omega}{\newOmega}}{6} (t+1)^2}{\eta m} .
    \end{align*}
    Given the above inequality,  for $t \geq 3$ we have
    \begin{align*}
        \theta_t \tau_t + \cL \epsilon_{t-1} = t \,\tfrac{t-1}{3} + \cL \epsilon_{t-1 } \leq \tfrac{t (t-1)}{3} + \tfrac{\tli{10 \revision{\Omega}{\newOmega}}{6}t^2 \cL}{\eta m} \leq \tfrac{t(t-1)}{3} + \tfrac{t-1}{3} = \tfrac{(t+1)(t-1)}{3},
    \end{align*}
    where the second inequality follows from $\eta \geq \tfrac{\tli{30 \revision{\Omega}{\newOmega}}{18} (k+2) \cL}{m}$.
    Combining the above bound with  $\theta_{t-1} (1 + \tau_{t-1}) = \tfrac{(t-1)(t+1)}{3}$, we arrive at
    \begin{align*}
        \theta_{t-1}(1+ \tau_{t-1}) - (\theta_t \tau_t + \cL \epsilon_{t-1})  \geq 0, \qquad t \geq 2,
    \end{align*}
    which is \eqref{eq:batchsizeBt_adaptive}. We conclude that the bound \rf{eq:GEM_theorem_result_general} of Theorem \rf{thm:main} holds.

On the other hand, when $\eta\geq \sqrt{\tfrac{\tli{10 \revision{\Omega}{\newOmega}}{2} (k+1)^3 \sigma_*^2}{m D ^2}}$ we have
    \begin{align*}
        \sum_{t=0}^{k-1} \epsilon_t \sigma_*^2  \leq \sum_{t=0}^{k-1} \tfrac{\tli{10 \revision{\Omega}{\newOmega}}{6} (t+1)^2\sigma_*^2}{\eta m} \leq \tfrac{\tli{10 \revision{\Omega}{\newOmega}}{2} (k+1)^3\sigma_*^2}{\eta m} \leq \sqrt{\tfrac{\tli{10 \revision{\Omega}{\newOmega}}{2}(k+1)^3 \sigma_*^2 D ^2}{m }}.
    \end{align*}
When substituting the above bound into \rf{eq:GEM_theorem_result_general} and noticing that $ \theta_k (1 + \tau_k) = \tfrac{k(k+2)}{3}$ and $\cL \epsilon_0 \leq \tfrac{1}{3}$, we conclude that
    \begin{align*}
         \bbe[ f(x_k) - f^*] \leq \tfrac{\tli{91}{73} L D ^2}{k(k+2)} + \tfrac{\tli{90 \revision{\Omega}{\newOmega}}{54} \cL D ^2}{m k} +  \sqrt{\tfrac{\tli{120\revision{\Omega}{\newOmega}}{72} \sigma_*^2 D ^2}{m k }}
    \end{align*}
    which completes the proof of \eqref{stepsize_non_strong_convex_2_4}.
\paragraph{2$^o$.} The ``Furthermore'' part of the statement immediately follows from \rf{bound_bregman} and the fact that
% Now let us prove Ineq.~\eqref{stepsize_non_strong_convex_2_5}. The only difference is that
when $\eta\geq \sqrt{\tfrac{\tli{20}{4}(k+1)^3\tli{\revision{}{\newOmega}}{} \sigma_*^2}{\tli{3}{}m \revision{}{\Omega} R^2}}$ one has
    \[%begin{align*}
        \sum_{t=0}^{k-1} \epsilon_t \sigma_*^2  \leq \tfrac{\tli{10 \revision{\Omega}{\newOmega}}{2}(k+1)^3\sigma_*^2}{\tli{3}{}\eta m} \leq \sqrt{\tfrac{\tli{5\revision{\Omega^2}{\Omega\newOmega} }{\Omega}(k+1)^3 \sigma_*^2 R^2}{\tli{3}{}m }}.\eqno{\mbox{\qed}}
    \]%end{align*}
%    Recalling Ineq.~\eqref{bound_bregman}, we have
%    \begin{align*}
%         \bbe[ f(x_k) - f^*] \leq \tfrac{91 L \Omega R_0^2}{2k(k+2)} + \tfrac{45 \Omega^2 \cL R_0^2}{m k} +  \sqrt{\tfrac{60 \Omega^2 \sigma_*^2 R_0^2}{m k }},
%    \end{align*}
%    and we complete the proof.\qed

\subsection{Proofs of Corollaries \ref{coro_strongly_convex} and  Corollary \ref{coro_SR}}
\paragraph{Proof of  Corollary \ref{coro_strongly_convex}.}
Note that bound \eqref{stepsize_non_strong_convex_2_5} of Corollary~\ref{coro_non_strongly_convex} implies that whenever size $m$ of the mini-batch satisfies
\[m \geq \max \left\{1, \tfrac{{3}(k+2)\cL}{L}, \tfrac{{16} N (N+2)^2 \sigma_*^2}{{9\Omega}L^2 R^2}\right\},\] we have for the approximate solution $x_N$ by SGE after $N$ iterations,
\begin{align}\label{conv_result}
	    \bbe[ f(x_N) - f^*] \leq  \tfrac{73 L \Omega R^2}{2N(N+2)} + \tfrac{27 \Omega \cL R^2}{m k} +  6\sqrt{\tfrac{\Omega \sigma_*^2 R^2}{m k }}\leq  \tfrac{{73} L \Omega R^2}{2N(N+2)} + \tfrac{{9} L \Omega R^2}{N(N+2)} + \tfrac{9 L \Omega R^2}{\aic{}{2}N(N+2)}=\tfrac{50 L \Omega R^2}{N(N+2)}
	\end{align}
where $R$ is an upper bound for the ``initial distance to $x^*$.''

Note that $\|\stx^0-x^*\|\leq R_0$. Let us now assume that for $1\leq k\leq K$, $\|\stx^{k-1}-x^*\|\leq R_k$, so that at the beginning of the $k$th stage
$\|x_0-x^*\|\leq R_k$.
Based on the above inequality, by the choice of $N$,
	\begin{align*}
	    \bbe[ f(x_N) - f^*]
	\le  \tfrac{50 L \Omega R_{k-1}^2}{N(N + 2)}\leq \tfrac{\mu R_{k-1}^2}{4}=2^{-k-1}\mu R_{0}^2= \tfrac{1}{2} \mu R_{k}^2.
	\end{align*}
Due to \rf{eq:general_quadractic_growth} this also means that
	\[
	    \bbe[\|x_N-x^*\|^2] \leq {R_k^2}= 2^{-k}R_{0}^2.\eqno{\mbox{\qed}}
	\]
\paragraph{Proof of Corollary \ref{coro_SR}.}
We have $\|\sstx^0-x^*\|\leq R_0$. Let us assume that $\|\sstx^{k-1}-x^*\|\leq R_{k-1}$ for some $1\leq k\leq K$.
From \eqref{conv_result} we conclude that
\begin{align*}
	    \bbe[ f(\stx^k) - f^*] \leq
	 \tfrac{50 L \Omega R_{k-1}^2 }{N(N + 2)}
	\leq \tfrac{\underline \kappa R_k^2}{16 s}={2^{-k-4}s^{-1}\underline \kappa R_0^2}
	\end{align*}
by the choice of $N$.
Furthermore, recalling that $\sstx^k=\mathsf{sparse}(\stx^k)$, we have
\beq\label{lemma_norm_tran}
    \|\sstx^k - x^*\|\leq \sqrt{2s}\|\sstx^k -x^*\|_2 \leq 2\sqrt{2s}\|\stx^k-x^*\|_2.\eeq
Indeed,
given that both $\sstx^k$ and $x^*$ are $s$-sparse, we conclude that $\sstx^k - x^*$ is $2s$-sparse, thus
$   \|\sstx^k - x^*\|\leq \sqrt{2s}\|\sstx^k -x^*\|_2.$ On the other hand, by the optimality of $\sstx^k$ for \rf{eq:sparsex},
\[
    \|\sstx^k - x^*\|_2\leq \|\sstx^k - \stx^k\|_2 + \|\stx^k - x^*\|_2 \leq 2 \|\stx^k - x^*\|_2.
\]
We conclude that
\begin{align*}
    \bbe[\|\sstx^k-x^*\|^2] \leq 8s\| \stx^k-x^*\|^2_2\leq 16s\underline \kappa^{-1}\bbe[ f(\stx^k) - f^*] \leq R_k^2= 2^{-k}R_0^2
\end{align*}
which completes the proof.\qed
% \newpage

%\section*{Funding and Conflicts of interests}
%Sasila Ilandarideva and Anatoli Juditsky were partially supported by Multidisciplinary Institute in Artificial intelligence MIAI {@} Grenoble Alpes ANR-19-P3IA-0003. Guanghui Lan and Tianjiao Li were partially supported by Division of Mathematical Sciences grant DMS-1953199 and Air Force Office of Scientific Research grant FA9550-22-1-0447. The authors have no other conflicts of interests to disclose.

\renewcommand\refname{Reference}

\bibliographystyle{abbrv}
\bibliography{paper_revised.bib}

\end{document}